\documentclass[a4paper,11pt]{amsart}
 \usepackage{graphicx,amsfonts,amsmath,amssymb,latexsym,cancel}
 \theoremstyle{plain}
 \newtheorem {hypo}{\bf\hspace{-\parindent}Hypothesis}
 
 \newtheorem {thm}{Theorem}
 \newtheorem {prop}[hypo]{Proposition}%[section]
 
 \newtheorem {lem}[hypo]{Lemma}%[section]

 \theoremstyle{definition}
 \newtheorem {defn}[hypo]{Definition}%[section]
 \newtheorem {eg}[hypo]{Example}%[section]

 \theoremstyle{remark}
 \newtheorem {rmk}[hypo]{Remark}%[section]

 \newcommand\omi{\omega_{4}}
 \newcommand\omx{\omega_X}
 \newcommand\omy{\omega_Y}
 \newcommand\omz{\omega_Z}

 \newcommand\Zb{\mathbb{Z}}
 
 \newcommand\Cb{\mathbb{C}}
 \newcommand\Pb{\mathbb{P}}
 \newcommand\Qb{\mathbb{Q}}

 \newcommand\ben{\begin{equation*}}
 \newcommand\ebn{\end{equation*}}
 \newcommand\be{\begin{equation}}
 \newcommand\eb{\end{equation}}
 \newcommand\ds{\displaystyle}

 \newcommand{\pf}{\begin{bpf}}

 \newcommand{\pfms}{\begin{bpfms}}
 \newcommand{\epf}{\end{bpf}\hfill$\square$\vspace{0.1cm}}
 \newcommand{\epfms}{\end{bpfms}\hfill$\square$\\ }

 \title{Algebraic solutions  of the sixth  Painlev\'e equation}
 \author{Oleg Lisovyy and Yuriy Tykhyy}
 \address{Laboratoire de Math\'ematiques et Physique Th\'eorique CNRS/UMR 6083,
 Universit\'e de Tours, Parc de Grandmont, 37200 Tours, France}
 \email{lisovyi@lmpt.univ-tours.fr}
 \address{Bogolyubov Institute for Theoretical Physics, 03680 Kyiv, Ukraine}
 \email{tykhyy@bitp.kiev.ua}
 \date{}

\begin{document}

 \begin{abstract}
 We describe all  finite orbits of an action of the extended modular group~$\bar{\Lambda}$ on
 conjugacy classes of $SL_2(\Cb)$-triples. The result is used to classify all algebraic
 solutions of the general Painlev\'e~VI equation up to parameter equivalence.
 \end{abstract}

 \maketitle

 \tableofcontents
 \newcounter{rem}
 \setcounter{rem}{1}
 \section{Introduction}
 Modular group $\Gamma=PSL_2(\Zb)$ consists of $2\times2$
 matrices with integer entries and unit determinant, considered up
 to overall sign. It has a presentation $\Gamma=\langle
 s,t\,|\,s^3=t^2=1\rangle$, and is known to be isomorphic to the quotient of
 3-braid group $\mathcal{B}_3$ by its center $\mathcal{Z}\cong \Zb$. The kernel of the
 canonical homomorphism $\Gamma\rightarrow PSL_2(\Zb_2)\cong S_3$
 defines a congruence subgroup $\Lambda\subset \Gamma$, also known as
 $\Gamma(2)$:
 \ben
 \Lambda=\left\{\left(\begin{array}{cc}a & b \\ c & d\end{array}\right)\in
 SL_2(\Zb)\,|\,a,d\text{ odd};\;b,c\text{ even}\right\}/\{\pm 1\}.
 \ebn
 There are isomorphisms $\Lambda\cong \mathcal{P}_3/\mathcal{Z}\cong
 \mathcal{F}_2$, where $ \mathcal{P}_3$ denotes the group of pure
 3-braids and $\mathcal{F}_2$ is the free group with 2 generators.

 Extended modular groups $\bar{\Gamma}$ and
 $\bar{\Lambda}$ are obtained by replacing the unit determinant
 condition with $ad-bc=\pm1$. These groups have the following
 presentations:
 \be\label{gammaext}
 \bar{\Gamma}=\langle r,s,t\,|\,r^2=s^3=t^2=(tr)^2=(sr)^2=1\rangle,
 \eb
 \be\label{lambdaext}
 \bar{\Lambda}=\langle x,y,z\,|\,x^2=y^2=z^2=1\rangle\cong C_2\ast
 C_2\ast C_2\,,
 \eb
 where
 \ben
 t=\left(\begin{array}{cr}0 & -1 \\ 1 & 0\end{array}\right),\quad
 s=\left(\begin{array}{cr}0 & -1 \\ 1 & 1\end{array}\right),\quad
 r=\left(\begin{array}{cr}0 & 1 \\ 1 & 0\end{array}\right),\ebn
 \ben
 x=rsts=\left(\begin{array}{rr}-1 & -2 \\ 0 &
 1\end{array}\right),\quad y=rt=\left(\begin{array}{cr}1 & 0 \\ 0 &
 -1\end{array}\right),\quad
 z=stsr=\left(\begin{array}{rr}1 & 0 \\ -2 &
 -1\end{array}\right).
 \ebn
 Note that $\Lambda$ is isomorphic to the subgroup (of
 index 2) of $\bar{\Lambda}$ containing words of even length in $x$, $y$,
 $z$. Hence, given a  $\bar{\Lambda}$ action on a set $U$ and a point $u\in U$, the
 orbits $\bar{\Lambda}(u)$ and $\Lambda(u)$ are simultaneously
 finite or infinite.

 In this paper the last observation is used to classify
 algebraic solutions of the sixth Painlev\'e equation (see \cite{ince}):
% (actually discovered by Gambier \cite{gambier}):
 \begin{align}
 \label{pvi}\frac{d^2w}{dt^2}=\;\frac{1}{2}\left(\frac{1}{w}+\frac{1}{w-1}
 +\frac{1}{w-t}\right)\left(\frac{dw}{dt}\right)^2 -
 \left(\frac{1}{t}+\frac{1}{t-1}+\frac{1}{w-t}\right)\frac{dw}{dt}\,+\tag{PVI}\\
 \nonumber +\,
 \frac{w(w-1)(w-t)}{2t^2(t-1)^2}\left((\theta_{\infty}-1)^2-\frac{\theta_x^2 t}{w^2}+
 \frac{\theta_y^2(t-1)}{(w-1)^2}+\frac{(1-\theta_z^2)t(t-1)}{(w-t)^2}\right).
 \end{align}
 This is the most general ODE of the form $w''=F(t,w,w')$, with $F$ rational
 in $w$, $w'$ and $t$, whose general solution has no movable branch points
 and essential singularities. It can therefore be analytically continued to
 a meromorphic function on the universal covering of $\Pb^1\backslash\{0,1,\infty\}$.

 A result from Watanabe \cite{watanabe} suggests that, roughly speaking, any solution
 of PVI is either a)~algebraic or b)~solves a Riccati equation or
 c)~cannot be expressed via classical functions. Known examples of algebraic solutions \cite{boalch_SL}
 turn out to be related to various mathematical structures, including e.g.
 Frobenius manifolds \cite{dubrovin2}, symmetry groups of regular polyhedra \cite{dubrovin,hitchin2},
 complex reflections \cite{boalch_klein},  Grothendieck's dessins d'enfants and
 their deformations \cite{andreev_kitaev,kitaev,kitaev2}. A few families of non-classical solutions have also been constructed
 in terms of Fredholm determinants, see \cite{borodin,lisovyy}.

 In the case $\theta_x=\theta_y=\theta_z=0$ a full classification of algebraic solutions
 has been obtained by Dubrovin and Mazzocco \cite{dubrovin}.
 Their approach, followed to some extent in the present work, is based on
 the description of PVI as the equation of monodromy
 preserving deformation of Fuchsian systems of the form
 \be\label{fuchsian}
 \frac{d\Phi}{d\lambda}=\left(\frac{A_x}{\lambda-u_x}+\frac{A_y}{\lambda-u_y}+
 \frac{A_z}{\lambda-u_z}\right)\Phi,\qquad A_{\nu}\in\mathfrak{sl}_2(\Cb),
 \eb
 where the poles $u_{\nu}$ are pairwise distinct,
 $A_{\nu}$ are $2\times 2$ matrices independent of $\lambda$ with eigenvalues
 $\pm \theta_{\nu}/2$  and
 \ben
 A_x+A_y+A_z=\left(\begin{array}{cc}-\theta_{\infty}/2 & 0 \\ 0 &
 \theta_{\infty}/2\end{array}\right), \qquad\theta_{\infty}\neq0.
 \ebn
  \begin{figure}[!h]
 \begin{center}
 \resizebox{6cm}{!}{
 \includegraphics{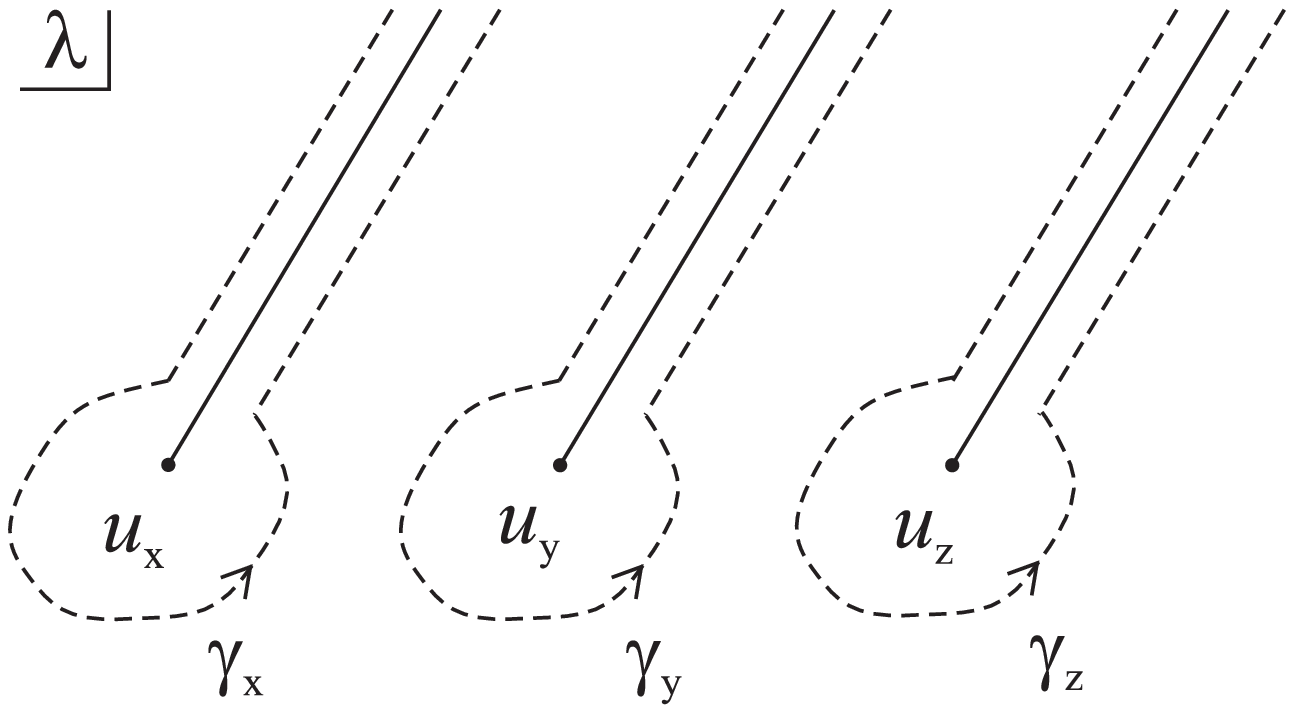}} \\   Fig. 1: Branch cuts and loops $\gamma_{x,y,z}$
 \end{center}
 \end{figure}

 The fundamental matrix  $\Phi(\lambda)$ is a multivalued analytic
 function on $\Cb\backslash\{u_x,u_y,u_z\}$. Fix a basis of loops
 and branch cuts in $\pi_1(\Pb^1\backslash\{u_x,u_y,u_z,\infty\},\infty)$ as shown in
 Fig.~1. To each branch of a solution of the PVI equation corresponds
 a unique (up to conjugation) triple of  monodromy matrices $(M_x,M_y,M_z)\in G^3$, $G= SL_2(\Cb)$
 of $\Phi(\lambda)$ w.r.t. the loops $\gamma_x$,  $\gamma_y$,  $\gamma_z$.
 One consequence of isomonodromy is that analytic continuation of  solutions of PVI induces an action of the
 pure braid group on 3 strings on the space of conjugacy classes of such triples (i.e.
 on the quotient $\mathcal{M}=G^3/G$ of three copies of $G$ by diagonal conjugation by $G$).
 It extends to the standard Hurwitz action of the braid group
 $\displaystyle \mathcal{B}_3=\langle \beta_x,\beta_z\,|\,\beta_x\beta_z\beta_x=\beta_z\beta_x\beta_z\rangle$
 on $G^3$. Explicitly,
 \begin{align}
 \nonumber&\beta_x:\left(M_x,M_y,M_z\right)\mapsto\left(M_x,M_z,M_zM_yM_z^{-1}\right), \\
 \nonumber&\beta_z:\left(M_x,M_y,M_z\right)\mapsto\left(M_y,M_yM_xM_y^{-1},M_z\right).
 \end{align}

 Observe that $\beta_z\beta_x$ acts on a representative triple
 $\left(M_x,M_y,M_z\right)\in\mathcal{M}$ by a cyclic
 permutation. The center $\mathcal{Z}$ of $\mathcal{B}_3$ is generated
 by $\left(\beta_z\beta_x\right)^3$ and therefore it acts on $\mathcal{M}$ trivially. This leads to
 an action of the modular group  $\Gamma\cong\mathcal{B}_3/\mathcal{Z}$ on $\mathcal{M}$, with
 \begin{align}
 \label{gammas} & s:\left(M_x,M_y,M_z\right)\mapsto \left(M_z,M_x,M_y\right),\\
 \label{gammat} & t:\left(M_x,M_y,M_z\right)\mapsto \left(M_z,M_y,M_yM_xM_y^{-1}\right)
 \end{align}
 in the above notation. The action of $\bar{\Gamma}$ on $\mathcal{M}$  is obtained by adding
 the involution
 \be
 \label{gammar} r:\left(M_x,M_y,M_z\right)\mapsto
 \left(M_z^{-1},M_y^{-1},M_x^{-1}\right).
 \eb
 \begin{lem} The transformations $s,t,r:\mathcal{M}\rightarrow\mathcal{M}$, as given by
 (\ref{gammas})--(\ref{gammar}), satisfy the defining relations (\ref{gammaext})
 of the extended modular group $\bar{\Gamma}$. \end{lem}
 \noindent
 As a corollary, we obtain the  restriction of the $\bar{\Gamma}$  action to its level 2
 subgroup~$\bar{\Lambda}$:
 \begin{lem} The generators $x,y,z\in \bar{\Lambda}$ act on  representative
 triples from $\mathcal{M}$ as follows:
 \begin{align}
 \nonumber & x:\left(M_x,M_y,M_z\right)\mapsto \left(M_x^{-1},M_y^{-1},M_xM_z^{-1}M_x^{-1}\right),\\
 \label{lamb_matr} & y:\left(M_x,M_y,M_z\right)\mapsto \left(M_yM_x^{-1}M_y^{-1},M_y^{-1},M_z^{-1}\right),\\
 \nonumber & z:\left(M_x,M_y,M_z\right)\mapsto \left(M_x^{-1},M_z M_y^{-1} M_z^{-1},M_z^{-1}\right).
 \end{align}
 \end{lem}
 \pf Both lemmas can be proved by direct calculation. \epf

 Let us now describe the last action  in more detail, introducing on $\mathcal{M}$
 a suitable set of coordinates. Following \cite{iwasaki}, to a point $(M_x,M_y,M_z)\in\mathcal{M}$
 we associate a 7-tuple $(p_x,p_y,p_z,p_{\infty},X,Y,Z)\in\Cb^7$ given by
  \be\label{pxyzi}
  p_{x}=\mathrm{Tr}\,M_{x},\quad
    p_{y}=\mathrm{Tr}\,M_{y},\quad
      p_{z}=\mathrm{Tr}\,M_{z},\quad
  p_{\infty}=\mathrm{Tr}\left(M_z M_y
  M_x\right),
  \eb
   \be \label{xyz}
  X=\mathrm{Tr}\left(M_yM_z\right),\qquad Y=\mathrm{Tr}\left(M_zM_x\right),\qquad
  Z=\mathrm{Tr}\left(M_xM_y\right).
  \eb
  Naive dimension of the quotient $\mathcal{M}$ is equal to 6 and thus it is not surprising that the above monodromy ivariants are
  not all independent --- there is a constraint
  \be\label{jfr}
  XYZ+X^2+Y^2+Z^2-\omx X-\omy Y-\omz Z+\omi=4\,,
  \eb
  where
  \be\label{omxyz}
  \omx= p_x p_{\infty}+p_y p_z,\qquad
  \omy= p_y p_{\infty}+p_z p_x,\qquad
  \omz= p_z p_{\infty}+p_x p_y,
  \eb
  \be\label{omi}
  \omi= p_x^2+p_y^2+p_z^2+p_{\infty}^2+p_xp_yp_zp_{\infty}.
  \eb
  \begin{rmk}
  Boalch \cite{boalch_klein} refers to an equation equivalent to (\ref{jfr}) as `Fricke relation'.
  In the context of Painlev\'e~VI, it was first obtained by Jimbo
  in \cite{jimbo},  p.1140.
  \end{rmk}
  \begin{rmk} Four quantities (\ref{pxyzi}) are related to PVI parameters
  by
  \be\label{pnu}
  \qquad\qquad p_{\nu}=2\cos\pi\theta_{\nu},\qquad\qquad\nu=x,y,z,\infty.
  \eb
  Remaining three parameters $X$, $Y$, $Z$ satisfying Jimbo-Fricke
  relation (\ref{jfr}) can  be generically thought of as giving two
  PVI integration constants.
  \end{rmk}
  The $\bar{\Lambda}$ action (\ref{lamb_matr}) is defined for any group
  $G$. That $G=SL_2(\Cb)$ in our case leads to important
  simplifications, in particular  $\mathrm{Tr}\,M=\mathrm{Tr}\,M^{-1}$ for any $M\in
  G$. Monodromy
  parameters (\ref{pxyzi}) are then fixed by the induced  action of
  $\bar{\Lambda}$, and quadratic functions~(\ref{xyz}) transform according to the following:
  \begin{lem}\label{lemxyz} The induced action of the generators $x,y,z\in\bar{\Lambda}$ on
  the parameters (\ref{xyz}) is
  \begin{align}
  \nonumber x(X,Y,Z)&=\left(\omx-X-YZ,\,Y,\,Z\right),\\
  \label{lxyz} y(X,Y,Z)&=\left(X,\,\omy-Y-ZX,\,Z\right),\\
  \nonumber z(X,Y,Z)&=\left(X,\,Y,\,\omz-Z-XY\right).
  \end{align}
  \end{lem}
  \pf
  Using again that for $M\in SL_2(\Cb)$ one has $\mathrm{Tr}\,M=\mathrm{Tr}\,M^{-1}$
  and also $\ds M+M^{-1}=\mathrm{Tr}\,M\cdot \mathbf{1}$ we find
  for example
  \ben
  x(X)=\mathrm{Tr}\left(M_y^{-1}M_xM_z^{-1}M_x^{-1}\right)=
  \mathrm{Tr}\left(M_yM_xM_zM_x^{-1}\right)=p_xp_{\infty}-
  \mathrm{Tr}\left(M_yM_xM_zM_x\right)=
  \ebn
  \ben
  =p_xp_{\infty}-YZ+\mathrm{Tr}\left(M_yM_z^{-1}\right)=p_xp_{\infty}+p_yp_{z}-X-YZ.
  \ebn
  Proof of the other relations follows in a similar manner.
 \epf
 \begin{rmk}
 After this work has been completed, we became aware of two recent papers \cite{cantat,iwasaki_uehara}, where
 the group $\bar{\Lambda}$ was introduced into Painlev\'e~VI context in a way similar to ours and in
 particular its action (\ref{lxyz}) on monodromy invariants has been computed (cf. relations
 (2.10)--(2.12) in  \cite{cantat} and formula (37) in \cite{iwasaki_uehara}).
 We also note another interesting recent preprint \cite{iwasaki3} on algebraic PVI solutions.
 \end{rmk}

  \noindent\textbf{Idea of classification}. Finite branch (in particular, algebraic) solutions of
  Painlev\'e~VI necessarily lead to finite orbits of the
  $\mathcal{P}_3/\mathcal{Z}\cong \Lambda$ action on the space $\mathcal{M}$ of
  conjugacy classes of monodromy. Classification of such orbits is
  equivalent to finding all finite orbits of the action
  (\ref{lamb_matr}) of the extended modular group  $\bar{\Lambda}$. Finally, the orbit
  $\bar{\Lambda}(m)$, $m\in\mathcal{M}$ can be finite only if the corresponding orbit
  of the induced $\bar{\Lambda}$~action (\ref{lxyz}) on $\Cb^3$ is
  finite.\vspace{0.1cm}

  \begin{rmk}
  One usually obtains explicit algebraic solution
  curves from monodromy by applying Jimbo's asymptotic
  formula \cite{jimbo} (or an appropriate modification of it) and computing sufficiently many terms in
  the Puiseux expansions of solutions near
  singular points. Another extremely useful tool, especially for solutions of high degree,
  are Kitaev's quadratic transformations \cite{kitaev_LMP,kitaev_vidunas}.
  \end{rmk}
  In the next
  section, we classify all finite orbits of the action
  (\ref{lxyz}) (Theorem~1). It then turns out that the
  resulting 7-tuples of monodromy invariants completely determine
  ${\Lambda}$-orbits in  $\mathcal{M}$ except in the case when
  $M_{x,y,z}$ can be simultaneously transformed into upper  triangular form.
  In Section~3, we give a complete (up to parameter equivalence) list of  Painlev\'e~VI solutions
  with finite branching. All of them are algebraic with one possible exception of Picard solutions;  in that way our explicit computation confirms
  a recent result by Iwasaki \cite{iwasaki2}.

  Somewhat unexpectedly for the authors, the solutions corresponding to all possible
  finite $\Lambda$-orbits have already appeared in  various papers \cite{andreev_kitaev,boalch_klein,boalch52,boalch_RH,boalch_HG,
  dubrovin2,dubrovin,hitchin,hitchin2,kitaev,kitaev2}. However, four of them
  (solutions~13, 24, 43 and 44 below) were published with
  misprints, which are fixed
  in the present paper.
  \vspace{0.2cm}\\
  {\small
  \textbf{Acknowledgements}. O.L. is grateful to Dublin Institute
  for Advanced Studies where the project of classification of algebraic Painlev\'e~VI
  solutions was first conceived (supported by the Schroedin\-ger
  Fellowship).  The work of Yu.T. is supported by
  Eiffel PhD Fellowship of French government.}
  \normalsize

 \section{Finite orbits of $\bar{\Lambda}$}
  \subsection{Orbit graphs} Our main subject in this section is the $\bar{\Lambda}$-action
  (\ref{lxyz}) which we consider as an action on $\Cb^3$ by fixing
  the parameters $\boldsymbol{\omega}=(\omega_X,\omega_Y,\omega_Z)$. To any orbit ${O}$ of this action
  we associate a 3-colored (pseudo)graph $\Sigma(O)$ as follows:
  \begin{itemize}
  \item the vertices of $\Sigma(O)$ represent distinct points $\mathbf{r}=(X,Y,Z)\in O$,
  \item two vertices $a,b\in\Sigma(O)$ are connected by an undirected
  edge of colour $x$, $y$ or $z$ if $x(a)=b$
  (resp. $y(a)=b$ or $z(a)=b$),
  \item if a point $a\in \Sigma(O)$ is fixed by the transformation $x$, $y$ or $z$,
  we assign to it a self-loop of the corresponding color.
  \end{itemize}
  In fact $\Sigma(O)$ is a Schreier coset graph as its
  vertices can be identified with the cosets of the stabilizer of
  any point in $O$. Also observe that the structure of (\ref{lxyz}) imposes a number of restrictions on
  $\Sigma(O)$, in particular it forbids multiple
  edges and simple cycles with only one edge of a given color.
  \begin{eg}\label{orbitexample}
  Set $\boldsymbol{\omega}=(0,1,1)$ and consider the orbit of the
  point $\mathbf{r}=(-1,1,1)$. It consists of 5 points with coordinates given
 below along with the orbit graph.
  \begin{center}\begin{tabular}{c|ccc}
    point & $X$ & $Y$ & $Z$ \\ \hline
    1 & $-1$ & 1 & $1$ \\
    2 & 0 & 1 & 1 \\
    3 & 0 & 1 & 0 \\
    4 & 0 & 0 & 0 \\
    5 & 0 & 0 & 1 \\
  \end{tabular}
  $\qquad\begin{array}{c} \includegraphics[height=3.0cm]{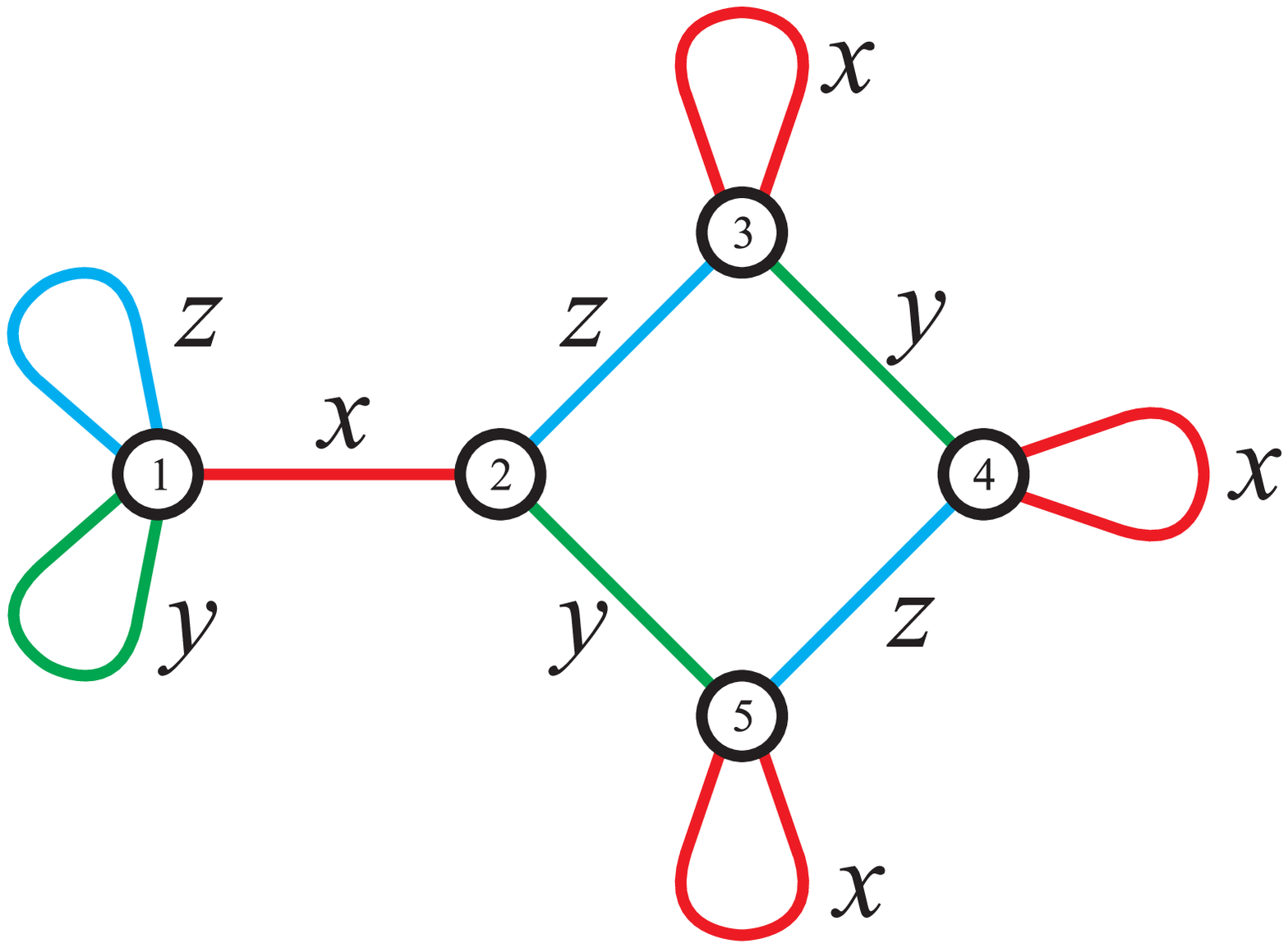}\end{array}$\end{center}
  This orbit does not split under the action of non-extended
  modular group $\Lambda$. The same result is immediate for any
  $\bar{\Lambda}$-orbit whose graph contains at least one self-loop (recall that $\Lambda$ consists of even-length words in $x$, $y$,
  $z$).
  \end{eg}
  \subsection{Symmetries} Before we move on to the classification,
  it is useful to look at the symmetries of the space of orbits
  and their relation to B\"acklund transformations for
  Painlev\'e~VI.

 Let $T:\mathcal{M}\rightarrow\mathcal{M}$ be an invertible map and
  let $\mathcal{O}\in\mathcal{M}$ be an orbit of the $\bar{\Lambda}$-action
  (\ref{lamb_matr}). If there exists an automorphism
  $\varphi\in\mathrm{Aut}\,\bar{\Lambda}$ compatible
  with $T$ (i.e.
  $\lambda\left(T(u)\right)=T\left(\varphi(\lambda)(u)\right)$ for any $\lambda\in\bar{\Lambda},
  u\in\mathcal{M}$), then $T(\mathcal{O})$ is also an orbit, and we will say that $\mathcal{O}$ and $T(\mathcal{O})$ are equivalent.
  The symmetries
  to be considered below are generated by
  \begin{itemize}
  \item permutations: $T:(M_x,M_y,M_z)\mapsto P(M_x,M_y,M_z)$, $\varphi:(x,y,z)\mapsto P(x,y,z)$ with some $P\in
  S_3$, where permutations act on $(x,y,z)$ in the standard way, and on the triples $(M_x,M_y,M_z)$ as follows:
  \ben
  \begin{array}{lll}
  (123)&: \bigl(M_x,M_y,M_z\bigr)\mapsto\bigl(M_z,M_x,M_y\bigr), \vspace{0.1cm}\\
  (12)(3)&: \bigl(M_x,M_y,M_z\bigr)\mapsto\bigl(M_y^{-1},M_x^{-1},M_z^{-1}\bigr).
  \end{array}
  \ebn
  \item sign flips: $T:(M_x,M_y,M_z)\mapsto (\varepsilon_x M_x,\varepsilon_y M_y,\varepsilon_z
  M_z)$, $\varepsilon_{x,y,z}=\pm1 $,
  $\varphi=id$.
  \end{itemize}
  To any orbit $O$ of the induced $\bar{\Lambda}$ action (\ref{lxyz})
  with parameters $\boldsymbol{\omega}\in\Cb^3$ therefore corresponds
  a number of equivalent orbits whose parameter triples are obtained from
  $\boldsymbol{\omega}$ by permutations and the action of the Klein four-group $K_4$
  (by sign changes of two coordinates). By virtue of (\ref{jfr}), all these orbits are characterized
  by the same value of $\omi$.
  To deal with nonequivalent orbits, we
  quotient the parameter space $\Cb^3$ by $K_4\rtimes S_3$, although it is convenient not to fix the fundamental domain
  explicitly.

  B\"acklund transformations (BTs) map solutions of a given
  Painlev\'e~VI equation to solutions of the same equation with different values of
  parameters $\theta_{x,y,z,\infty}$. The list of fundamental BTs
  for PVI is given in the table below, cf.
  \cite{noumi}:

  \begin{center}\begin{tabular}{c||c|c|c|c||c|c||c|c|c|c|}
   & $\theta_x$ &   $\theta_y$  & $\theta_z$  & $\theta_{\infty}$  & $w$ & $t$ & $\omega_X$ &  $\omega_Y$ & $\omega_Z$ & $\omi$\\
   \hline\hline
  $s_x$ & $-\theta_x$ &  $\theta_y$  & $\theta_z$  & $\theta_{\infty}$  & $w$ & $t$ &  $\omega_X$ &  $\omega_Y$ & $\omega_Z$ & $\omi$\\
  \hline
  $s_y$ & $\theta_x$ &  $-\theta_y$  & $\theta_z$  & $\theta_{\infty}$ & $w$ & $t$ &  $\omega_X$ &  $\omega_Y$ & $\omega_Z$ & $\omi$ \\
  \hline
  $s_z$ & $\theta_x$ &  $\theta_y$  & $-\theta_z$  & $\theta_{\infty}$ & $w$ & $t$  &  $\omega_X$ &  $\omega_Y$ & $\omega_Z$ & $\omi$\\
  \hline
  $s_{\infty}$ & $\theta_x$ &  $\theta_y$  & $\theta_z$  & $2-\theta_{\infty}$ & $w$ & $t$  &  $\omega_X$ &  $\omega_Y$ & $\omega_Z$ & $\omi$\\
   \hline
  $s_{\delta}$ & $\theta_x-\delta$ &  $\theta_y-\delta$  & $\theta_z-\delta$  & $\theta_{\infty}-\delta$ & $w+\frac{\delta}{p}$ & $t$  &  $\omega_X$ &  $\omega_Y$ &
  $\omega_Z$ & $\omi$ \\
  \hline\hline
  $r_x$ & $\theta_{\infty}-1$ & $\theta_z$ & $\theta_y$ & $\theta_{x}+1$ & $t/w$ & $t$  & $\omega_X$ & $-\omega_Y$ & $-\omega_Z$ & $\omi $\\
  \hline
  $r_y$ & $\theta_z$ & $\theta_{\infty}-1$ & $\theta_x$ & $\theta_y+1$& $\frac{w-t}{w-1}$ & $t$  & $-\omega_X$ & $\omega_Y$ & $-\omega_Z$ & $\omi $ \\
  \hline
  $r_z$ & $\theta_y$ & $\theta_x$ & $\theta_{\infty}-1$ & $\theta_z+1$ & $\frac{t(w-1)}{w-t}$ & $t$  & $-\omega_X$ & $-\omega_Y$ & $\omega_Z$ & $\omi $  \\
   \hline\hline
  $P_{xy}$ & $\theta_y$ & $\theta_x$ & $\theta_z$ & $\theta_{\infty}$ & $1-w$ & $1-t$  & $\omega_Y$ & $\omega_X$ & $\omega_Z$ & $\omi$  \\
  \hline
  $P_{yz}$ & $\theta_x$ & $\theta_z$ & $\theta_y$ & $\theta_{\infty}$ & $w/t$ & $1/t$  & $\omega_X$ & $\omega_Z$ & $\omega_Y$ & $\omi$  \\
  \hline\end{tabular}\vspace{0.2cm}\\
  \nopagebreak[2] Table 1: B\"acklund transformations for Painlev\'e
  VI\vspace{0.2cm}\end{center}
  Here we use the standard notation
  $\ds\delta=\frac{\theta_x+\theta_y+\theta_z+\theta_{\infty}}{2}$ and
  \ben
  2p=\frac{t(t-1)w'}{w(w-1)(w-t)}-\left(\frac{\theta_x}{w}+\frac{\theta_y}{w-1}+\frac{\theta_z+1}{w-t}\right).
  \ebn

  \begin{rmk}\label{shifts}Five transformations $s_{\nu}$ ($\nu=x,y,z,\infty,\delta$)
  generate affine Weyl group of type~$D_4$. Using these
  transformations,
  one can construct shift operators
  \begin{eqnarray*}
  &t_{x}=s_x s_{\delta} \left(s_y s_z  s_{\infty} s_{\delta}\right)^2,\qquad
  &t_{y}=s_y s_{\delta} \left(s_x s_z  s_{\infty} s_{\delta}\right)^2,\\
  &t_{z}=s_z s_{\delta} \left(s_x s_y  s_{\infty} s_{\delta}\right)^2,\qquad
  &t_{\infty}=s_{\infty} s_{\delta} \left(s_{x} s_y s_z s_{\delta}\right)^2,
  \end{eqnarray*}
  acting on the parameter space by simple translations:

  \begin{center}\begin{tabular}{c||c|c|c|c|}
  & $\theta_x$ & $\theta_y$ & $\theta_z$ & $\theta_{\infty}$ \\
  \hline\hline
  $t_x$ & $\theta_x+2$ & $\theta_y$ & $\theta_z$ &
  $\theta_{\infty}$\\
  \hline
   $t_y$ & $\theta_x$ & $\theta_y+2$ & $\theta_z$ &
  $\theta_{\infty}$\\
  \hline
   $t_z$ & $\theta_x$ & $\theta_y$ & $\theta_z+2$ &
  $\theta_{\infty}$\\
  \hline
 $t_{\infty}$ & $\theta_x$ & $\theta_y$ & $\theta_z$ &
  $\theta_{\infty}+2$\\
  \hline
  \end{tabular}\vspace{0.2cm}\end{center}
  Enlarging affine $D_4$ by the Klein four-group $K_4\cong\langle
  r_x,r_y,r_z\rangle$ gives extended affine Weyl group $D_4$. Full Okamoto affine $F_4$ action
  involves additional generators $P_{xy}$, $P_{yz}$ changing PVI independent variable $t$ by M\"obius transformations of
  $\Pb^1$ permuting 0, 1 and $\infty$.
  \end{rmk}

  Last four columns of Table~1 describe the action of BTs on parameters $\omega_{X,Y,Z,4}$ defined by (\ref{omxyz})--(\ref{pnu}).
  Observe that all BTs lead to equivalent points in the
  parameter space of orbits of the induced $\bar{\Lambda}$ action (\ref{lxyz}). We now want to prove a converse statement:
  \begin{prop}\label{inj} Given $\omega_X,\omega_Y,\omz,\omi\in\Cb$,
  consider (\ref{omxyz})--(\ref{pnu}) as a system of equations for
  unknown $\theta_{x,y,z,\infty}$. Any two solutions of this system are related by the
  affine $D_4$ transformations introduced above.
  \end{prop}
  \pf Choose an arbitrary solution $\{\theta_{\nu}^{0}\}$ ($\nu=x,y,z,\infty$) and denote $p_{\nu}^{0}=2\cos\pi\theta_{\nu}^{0}$.
   Introduce the auxiliary variable $\xi=p_x^2+p_y^2+p_z^2+p_{\infty}^2$. It satisfies the
  cubic equation
  \be\label{cubicxi}
  \xi^3-a(\omega)\xi^2+b(\omega)\xi-c(\omega)=0,
  \eb
  where
  \begin{eqnarray*}
  &a(\omega)=\omi+16,\qquad b(\omega)=\omx\omy\omz-4(\omx^2+\omy^2+\omz^2)+32\omi, \\
  &c(\omega)=\omx^2\omy^2+\omx^2\omz^2+\omy^2\omz^2-4\omi(\omx^2+\omy^2+\omz^2)+16\omi^2.
  \end{eqnarray*}
  Write $\omega_{X,Y,Z,4}$ in terms of $\{p_{\nu}^{0}\}$, then three roots of
  (\ref{cubicxi}) are
  \begin{eqnarray*}
  \xi_0&=&\left(p_x^0\right)^2+\left(p_y^0\right)^2+\left(p_z^0\right)^2+\left(p_{\infty}^0\right)^2,\\
  \xi_{\pm}&=&
 8\left(1\;\;+\!\!\!\prod\limits_{\nu=x,y,z,\infty}\cos\pi\theta_{\nu}^{0}\;\;\pm\!\!\!\prod\limits_{\nu=x,y,z,\infty}\sin\pi\theta_{\nu}^{0}\right).
  \end{eqnarray*}
  Applying $s_{\delta}$ (or $s_{\delta}s_x$) to initial solution
  $\{\theta_{\nu}^{0}\}$ gives a solution with $\xi=\xi_-$
  (resp. $\xi=\xi_+$). Therefore it is sufficient to prove the Proposition for
  solutions of (\ref{omxyz})--(\ref{pnu}) with $\xi=\xi_0$.

  Assume that at least two of three numbers $\omx^2,\omy^2,\omz^2\in\Cb$ are distinct,
  say $\omy^2\neq\omz^2$. Substituting  $\xi=\xi_0$ into easily verified relations
  \ben
  (p_x\pm p_{\infty})^4-(\xi\pm 2\omx)(p_x\pm p_{\infty})^2+(\omy\pm\omz)^2=0
  \ebn
  we find  $\ds\left(p_x+ p_{\infty}\right)^2=\left(p^0_x+ p^0_{\infty}\right)^2\text{ or
  }\left(p^0_y+  p^0_{z}\right)^2$,  $\ds\left(p_x- p_{\infty}\right)^2=\left(p^0_x- p^0_{\infty}\right)^2\text{ or
  }\left(p^0_y-  p^0_{z}\right)^2$.
  Also if $\xi=\xi_0$  then
  \ben
  p_xp_yp_zp_{\infty}=\omi-\xi= p_x^0p_y^0p_z^0p_{\infty}^0,\qquad\qquad
  p_xp_{\infty}+p_yp_z=\omx=p^0_xp^0_{\infty}+p^0_yp^0_z,
  \ebn
  so that $p_xp_{\infty}=p^0_xp^0_{\infty}\text{ or
  }p^0_yp^0_{z}$. But now if e.g. $\ds\left(p_x+ p_{\infty}\right)^2=\left(p^0_x+
  p^0_{\infty}\right)^2$,
  $\ds\left(p_x-p_{\infty}\right)^2=\left(p^0_y-
  p^0_{z}\right)^2$, combining with the latter result we find $\left(p^0_x+
  p^0_{\infty}\right)^2=\left(p^0_y+
  p^0_{z}\right)^2$ (for $p_xp_{\infty}=p^0_yp^0_{z}$) or $\left(p^0_x-
  p^0_{\infty}\right)^2=\left(p^0_y-
  p^0_{z}\right)^2$ (for $p_xp_{\infty}=p^0_xp^0_{\infty}$).
  Therefore we necessarily have
  \be\label{pxpi}
  \begin{cases}
  \left(p_x+ p_{\infty}\right)^2=\left(p^0_x+
  p^0_{\infty}\right)^2,\\
  \left(p_x-p_{\infty}\right)^2=\left(p^0_x-
  p^0_{\infty}\right)^2,
  \end{cases}\qquad \text{or}\qquad
  \begin{cases}  \left(p_x+ p_{\infty}\right)^2=\left(p^0_y+
  p^0_{z}\right)^2,\\
  \left(p_x-p_{\infty}\right)^2=\left(p^0_y-
  p^0_{z}\right)^2.
  \end{cases}
  \eb
   Choose a solution of (\ref{pxpi}) for $p_x$ and $p_{\infty}$, then $p_y$ and $p_z$ are unambigously fixed by
  \ben
  (p_{\infty}\pm p_x)(p_y\pm p_z)=\omy\pm\omz=(p^0_{\infty}\pm p^0_{x})(p^0_y\pm p^0_z)
  \ebn
  (here we used that $\omy^2\neq\omz^2$). Hence there are
  8 possible solutions for $(p_x,p_y,p_z,p_{\infty})$, namely
  \be\label{pxyzid4}
  \begin{array}{rr}
  (\pm p_x^0,\pm p_y^0,\pm p_z^0, \pm p_{\infty}^0),\quad
  &(\pm p_y^0,\pm p_x^0,\pm p_{\infty}^0, \pm p_{z}^0),\vspace{0.1cm}\\
  (\pm p_z^0,\pm p_{\infty}^0,\pm p_{x}^0, \pm p_{y}^0),\quad
  &(\pm p_{\infty}^0,\pm p_z^0,\pm p_y^0, \pm p_{x}^0).
  \end{array}
  \eb
   All of them
  can be obtained from $\left\{p^0_{\nu}\right\}$ using three
  affine $D_4$ transformations
  $\left(s_xs_ys_zs_{\infty}s_{\delta}\right)^2$, $s_{\delta}s_{x}s_y s_{\delta} s_z s_{\infty}$
  and $s_{\delta}s_{x}s_z s_{\delta} s_y s_{\infty}$. Now given $\left\{p_{\nu}\right\}$,
  all possible solutions for $\left\{\theta_{\nu}\right\}$ are clearly related by
  the transformations $\{s_{\nu}\}$, $\{t_{\nu}\}$, see Remark~\ref{shifts}.

  Now let $\omx^2=\omy^2=\omz^2$. We can set for
  definiteness
  $\omx=\omy=\omz$, then three out of four
  $p_{\nu}$ are equal. Denote this common value by $p$ and let
  $\tilde{p}$ be the fourth variable. Then
  \be\label{aux01}
  \omx=p(p+\tilde{p}),\qquad \omi=3p^2+\tilde{p}^2+p^3\tilde{p}.
  \eb
  Choose a solution $(p^0,\tilde{p}^0)$ of (\ref{aux01}). If $\omx\neq0$ then the only other
  solution such that $3p^2+\tilde{p}^2=\xi_0=3\left(p^0\right)^2+\left(\tilde{p}^0\right)^2$
  is given by $p=-p^0$, $\tilde{p}=-\tilde{p}^0$. Thus  $(p_x,p_y,p_z,p_{\infty})$ can only be a permutation
  of $(p^0,p^0,p^0,\tilde{p}^0)$ or
  $(-p^0,-p^0,-p^0,-\tilde{p}^0)$, which yields at most 8 distinct solutions. As above, all these 4-tuples are
  related by $\left(s_xs_ys_zs_{\infty}s_{\delta}\right)^2$, $s_{\delta}s_{x}s_y s_{\delta} s_z s_{\infty}$
  and $s_{\delta}s_{x}s_z s_{\delta} s_y s_{\infty}$. Now if
  $\omx=0$ there are 2 possibilities: 1) $p^0=0$, then the only other
  solution of (\ref{aux01}) with the same value of $\xi$ has the
  form $p=0$, $\tilde{p}=-\tilde{p}^0$; 2) $\tilde{p}^0=-p^0$,
  then the only such solution is $p=-p^0$, $\tilde{p}=p^0$.
  Clearly in both cases possible 4-tuples
  $(p_x,p_y,p_z,p_{\infty})$ are related by the affine $D_4$
  transformations.
  \epf
  \begin{rmk}
  We have just shown that the map
  \be\label{inject}
  \rho:{\begin{array}{c}\text{\rm parameter} \\ \text{\rm space of
  PVI}\end{array}}\Bigl/\text{\rm affine
  }D_4\;\;\rightarrow\;\;\Cb^4,\qquad
  [\theta_x,\theta_y,\theta_z,\theta_{\infty}]\mapsto(\omx,\omy,\omz,\omi)
  \eb
  is injective. Direct calculation shows that $\rho$ is in fact a
  bijection. Moreover the same result holds true if we replace in
  (\ref{inject}) affine $D_4$ by the full affine $F_4$ action and quotient the
  set of all triples $(\omx,\omy,\omz)$ by $K_4\rtimes S_3$ as described above.
  \end{rmk}
  \begin{rmk}
  It is more delicate to establish the equivalence of actual PVI
  solutions as BTs may become singular ($w(t)=0,1,t$ or $p=0$)
  in the way of transforming a given solution into another one with
  equivalent parameters.
%  We will come back to this question in  Section~3.
  \end{rmk}
  \subsection{2-colored suborbits}\label{ss2c} Take a point
  $\mathbf{r}=(X,Y,Z)\in\Cb^3$, fix
  $\boldsymbol{\omega}\in\Cb^3$ and consider the
  suborbit $O_{yz}(\mathbf{r})$ of the $\bar{\Lambda}$ action (\ref{lxyz}), generated from $\mathbf{r}$ by two transformations $y$
  and $z$. Clearly all points of $O_{yz}(\mathbf{r})$ have the same first
  coordinate $X$. We set $Y_0=Y$, $Z_0=Z$ and label remaining
  coordinates as shown on the suborbit graph below.
  \ben\includegraphics[height=1.3cm]{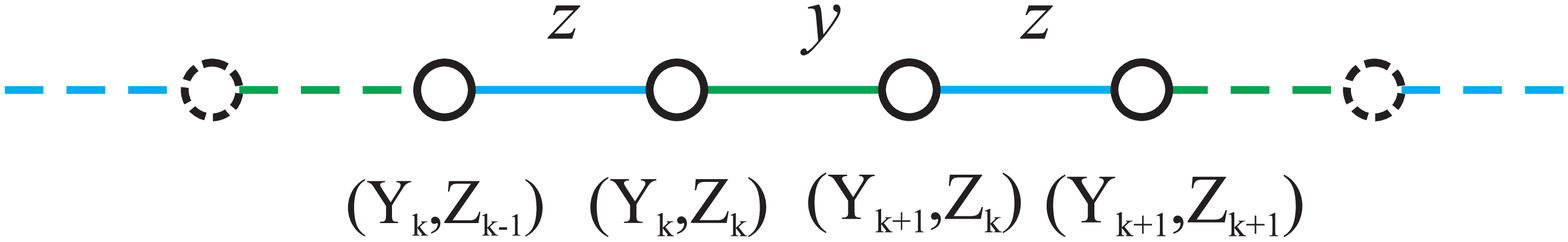}\ebn

  From (\ref{lxyz}) one finds a first order linear inhomogeneous difference equation
  \be\label{yzrecursion}
  \left(\begin{array}{c} Y_{k+1} \\ Z_{k+1}\end{array}\right)=
  \left(\begin{array}{cc} -1 & -X \\ X & X^2-1 \end{array}\right)
  \left(\begin{array}{c} Y_{k} \\ Z_{k}\end{array}\right)+
  \left(\begin{array}{c} \omy \\ \omz -  X\omy\end{array}\right).
  \eb
  A straightforward computation gives
  \begin{lem}\label{yzsolutions} If $X\neq \pm 2$, then the general solution of (\ref{yzrecursion}) is
  \be\label{yzsolution1}
   \left(\begin{array}{c} Y_{k} \\
  Z_{k}\end{array}\right)=\frac{1}{\sin{\lambda}/{2}}
  \left(\begin{array}{cc} \sin\frac{(1-2k)\lambda}{2} & -\sin k\lambda \\
  \sin k\lambda &  \sin\frac{(1+2k)\lambda}{2} \end{array}\right)
  \left(\begin{array}{c} \alpha \\
  \beta\end{array}\right)+\frac{1}{4-X^2}
  \left(\begin{array}{c} 2\omy-X\omz \\
  2\omz-X\omy\end{array}\right),
  \eb
   where $\alpha$, $\beta$ are arbitrary constants and
  $X=2\cos\lambda/2$. For $X=\pm2$ we have
  \be\label{yzsolution2}
  \left(\begin{array}{c} Y_{k} \\
  Z_{k}\end{array}\right)=\left(\begin{array}{cc}
  1-2k & \mp 2k \\ \pm 2k & 1+2k \end{array}\right)
    \left(\begin{array}{c} \alpha \\ \beta\end{array}\right)+
  \left(\begin{array}{c} \frac{\omy\pm\omz}{8}-\frac{(\omy\mp\omz)k}{2}+(\omy\mp\omz)k^2 \\
  \frac{\omz\pm\omy}{8}+\frac{(\omz\mp\omy)k}{2}+(\omz\mp\omy)k^2\end{array}\right).
  \eb
  \end{lem}

  Now assume that $O_{yz}(\mathbf{r})$ is finite.
  We call  the length of $O_{yz}(\mathbf{r})$
  the smallest  positive integer $N$ such that $Y_{k+N}=Y_k$, $Z_{k+N}=Z_k$.
  Since $x$, $y$, $z$ are involutions, the graph of
  any 2-colored finite suborbit can only be a simple cycle (as the length 2 $yz$-suborbit 2-3-4-5 in
  Example~\ref{orbitexample}) or a line with a self-loop at each
  of its ends (as e.g. the length 3 $xz$-suborbit 1-2-3 or the length 2 $xy$-suborbit 3-4 of
  the same  example).
  \begin{lem}\label{orbcosines}
  Let $N$ be the length of $O_{yz}(\mathbf{r})$.
  If $N>1$, then $X=2\cos \pi{n_X}/{N}$, where $n_X$ is an integer relatively prime to $N$
  satisfying  $0<n_X<N$.
  \end{lem}
  \pf Let $X\neq\pm 2$ and impose $Y_{k+N}=Y_k$, $Z_{k+N}=Z_k$ in
  (\ref{yzsolution1}). This gives $\sin\frac{N\lambda}{2}=0$, otherwise $\alpha=\beta=0$
  and hence
  $N=1$. Therefore $\lambda=2\pi n_X/N$, $n_X\in\Zb$ and
  we can choose $0<n_X<N$. Clearly $n_X$ and $N$ are coprime; otherwise $N$ is not the smallest period
  of  (\ref{yzsolution1}).

  Now if $X=\pm2$, then substituting $Y_{k+N}=Y_k$, $Z_{k+N}=Z_k$
  into (\ref{yzsolution2}) we find two conditions: 1) $\omy\mp\omz=0$ and 2) $\alpha\pm\beta=0$.
  This in turn implies that $Y_k=\mathrm{const}$, $Z_k =\mathrm{const}$,
  i.e.  $O_{yz}(\mathbf{r})$ consists of a single point. \epf
  \begin{defn}
  Let $O\subset \Cb^3$ be an orbit of the $\bar{\Lambda}$ action (\ref{lxyz}). A point
  $\mathbf{r}\in O$ is called \textit{good} if it is not fixed
  by at least two of three transformations $x$, $y$, $z$;
  otherwise we say that $\mathbf{r}$ is a \textit{bad} point.
  \end{defn}
  The case when the whole orbit consists of a single point is
  trivial. Hence below by a bad point we most often mean a point fixed by
  two transformations. The orbit graph has then two self-loops at the corresponding
  vertex.
  \begin{eg} The point 1 in
  Example~\ref{orbitexample} is bad, and the others are good.
  \end{eg}
  \begin{lem}\label{xyzcosineslemma}
  Let $O\subset \Cb^3$ be a {finite} orbit of (\ref{lxyz}). If
  $\mathbf{r}=(X,Y,Z)\in O$ is a good point, then
  \be\label{xyzcosines}
  X=2\cos \pi r_X,\qquad Y=2\cos \pi r_Y,\qquad Z=2\cos \pi r_Z,
  \eb
  where $r_{X,Y,Z}\in\Qb$ and $0<r_{X,Y,Z}<1$. If $\mathbf{r}\in O$ is a bad point, fixed by $y$ and $z$
  but not by $x$, then (\ref{xyzcosines}) still
  holds for $Y$ and $Z$.
  \end{lem}
  \pf If $\mathbf{r}$ is not fixed by $x$, then the lengths of $xz$- and
  $xy$-suborbit of $\mathbf{r}$ are strictly greater than 1. If
  $\mathbf{r}$ is good the same is true for each of three 2-colored suborbits of
  $\mathbf{r}$. Both statements then follow from Lemma~\ref{orbcosines}. \epf

  \subsection{Main technical lemma}
 This subsection is devoted to a technical result
  to be extensively used later. Namely, we want to find all rational solutions of the equation
 \be\label{sumcos}\sum_{j=1}^n\cos2\pi\varphi_j=0\eb
 with $n\leq6$. Without loss of generality we assume
 that $0\leq\varphi_j<1$ and consider the $n$-tuples $(\varphi_1,\ldots,\varphi_n)$ related by
 permutations, transformations $\varphi_j\rightarrow 1-\varphi_j$ and by the
 simultaneous change  $\varphi_j\rightarrow 1/2-\varphi_j$ as
 equivalent.
 \begin{defn} A rational $n$-tuple $(\varphi_1,\ldots,\varphi_n)$
 is called irreducible if it satisfies (\ref{sumcos}) and
 $\sum_{j\in E} \cos2\pi \varphi_j\neq 0$ for any proper
 subset $E\subset\{1,\ldots,n\}$.
 \end{defn}
 It then suffices to classify irreducible $n$-tuples
 $(\varphi_1,\ldots,\varphi_n)$ with $n\leq 6$.
  We first prove an auxiliary result concerning rational solutions of
  the equation
   \be\label{eqexp}
 \sum_{j=1}^n e^{2\pi i {\varphi}_j}=0.
 \eb
 Again we can assume that $0\leq \varphi_j<1$ and consider the solution $n$-tuples up to
 permutations. Also note that the shift of all $\varphi_j$ by a
 common phase $\varphi\in\Qb$ yields another solution.
 \begin{lem}\label{lemexp}
 All inequivalent irreducible (in the sense that $\displaystyle \sum_{j\in E}e^{2\pi i {\varphi}_j}\neq0 $
 for any proper subset $E\subset\{1,\ldots,n\}$) rational
 $n$-tuples with $n\leq 6$ solving (\ref{eqexp}) are given by
 \begin{itemize}
 \item  the 6-tuple
 \be\label{exp6tuple}
 \left(\varphi-\frac16,\varphi+\frac16,\varphi+\frac15,\varphi+\frac25,
 \varphi+\frac35,\varphi+\frac45\right),
 \eb
 \item the 5-tuple
 \be\label{exp5tuple} \left(\varphi,\varphi+\frac15,\varphi+\frac25,
 \varphi+\frac35,\varphi+\frac45\right),
 \eb
 \item the triple
 $\ds\left(\varphi,\varphi+\frac13,\varphi+\frac23\right)$ and
 the pair $\ds\left(\varphi,\varphi+\frac12\right)$,
 \end{itemize}
 with $\varphi\in\Qb$.
 \end{lem}
 \pf
  First part of the proof follows
 \cite{crosby, dubrovin}. Write ${\varphi}_k=\frac{n_k}{d_k}$,
 where $k=1,\ldots,n$ ($1<n\leq 6$) and $d_k$, $n_k$ are either positive coprime integers with
 $d_k>n_k$ or $n_k=0$. Let $p$ be a prime which is a divisor of at
 least one of $d_1,\ldots,d_n$, and denote by $\delta_k$, $l_k$, $c_k$,
 $\nu_k$  the integers such that
 \ben
 d_k=\delta_k p^{l_k},\qquad n_k=c_k\delta_k+\nu_k p^{l_k},
 \ebn
 where $\delta_k$ is prime to $p$, $0\leq c_k<p^{l_k}$; $c_k$
 is prime to $p$ for $l_k\neq 0$, otherwise  $c_k=0$. Then
 \ben
 \varphi_k=f_k+\frac{c_k}{p^{l_k}},\qquad
 f_k=\frac{\nu_k}{\delta_k}\,.
 \ebn
 Reordering $\varphi_1,\ldots,\varphi_n$ so that $l_1\geq l_2\geq\ldots\geq
 l_n$, we define the function
 \ben
 g_k(x)=\begin{cases}
 e^{2\pi i f_k} x^{c_k p^{l_1-l_k}} & \text{if}\; c_k\neq0, \\
 \qquad e^{2\pi i\varphi_k} & \text{if}\; c_k=0,
 \end{cases}
 \ebn
 and the polynomial
 \be\label{polu1}
 U(x)=\sum_{k=1}^n g_k(x).
 \eb
 By construction $g_k\left(\exp\left(\frac{2\pi
 i}{p^{l_1}}\right)\right)=e^{2\pi i\varphi_k}$, and
 (\ref{eqexp}) then implies that $U\left(\exp\left(\frac{2\pi
 i}{p^{l_1}}\right)\right)=0$.

 It is known since 1854  \cite{kronecker}
  that the polynomial
 \ben
 P(x)=1+x^{p^{l_1-1}}+x^{2p^{l_1-1}}+\ldots+x^{(p-1)p^{l_1-1}}
 \ebn
 is irreducible in the ring of polynomials with coefficients in
 any extension of the form $\mathbb{Q}(\zeta_1,\ldots,\zeta_m)$,
 where $\zeta_j$ is a root of unity of the order coprime with~$p$.
 Since $P\left(\exp\left(\frac{2\pi i}{p^{l_1}}\right)\right)=0$,
 then either (a) $U(x)\equiv
 0$ or (b) $U(x)\not\equiv 0$ is divisible by $P(x)$.\vspace{0.1cm}\\
 \underline{Case (a)}. The powers $c_k p^{l_1-l_k}$, appearing in the
 functions $g_k(x)$, are all equal. Otherwise one could write
 $U(x)$ as a sum of at least two polynomials equal to $0$, and the
 irreducibility condition fails. Therefore $l_k=l_1$,
 $c_k=c_1$. Now it is sufficient to subtract common phase $\frac{c_1}{p^{l_1}}$
 from all $\varphi_k$ to eliminate $p$ from all denominators.
 \vspace{0.1cm}\\
 \underline{Case (b)}. Write $U(x)=P(x)Q(x)$. The degree of
 $U(x)$ is at most $p^{l_1}-1$, hence the degree of $Q(x)$ is at most
 $p^{l_1-1}-1$. Then the numbers $N_U$ and $N_Q$ of different
 powers of $x$ in $U(x)$ and $Q(x)$ must be related by $N_U=p
 N_Q$. In particular, since in our case $N_U\leq 6$, the prime $p$ can only be equal to 2, 3 or 5.

 The powers $c_k p^{l_1-l_k}$ are all equal modulo $p^{l_1-1}$ to
 $s$, where $s$ is some integer independent of $k$, $0\leq s<p^{l_1-1}$.
 Otherwise one could collect powers corresponding to different $s$ and write $U(x)$ as
 a sum of at least two polynomials, each of them either divisible by $P(x)$ or
 vanishing identically. Corresponding $n$-tuple is then reducible,
 therefore we can only have $N_Q=1$, $Q(x)=\alpha x^s$.

  Suppose that $l_1\geq2$. Since $c_1$ is prime to $p$,  $s$ is also prime to $p$
  and all $n$ powers of $x$ that appear in the functions
 $g_k(x)$ are not divisible by $p^{l_1-1}$ and by $p$; in particular, all
 $c_k$ are non-zero. This in turn implies that $l_k=l_1$ for any $k$. Now $c_k=s+N_k
 p^{l_1-1}$ and subtracting from all $\varphi_k$ the common phase
 $\frac{s}{p^{l_1}}$ eliminates all higher (greater than
 1) powers of $p$  from the denominators.

 It remains to consider $l_1=1$, $p=2,3$ or $5$:\\
 (b.1) Let $l_1=1$, $p=5$, then $n=5$ or $6$. If $n=6$, then from
 $U(x)=\alpha x^s P(x)$ four out of six phases are equal, say
 $f_1=f_2=f_3=f_4$, and the remaining two satisfy
 $e^{2\pi i f_5}+e^{2\pi i f_6}=e^{2\pi i f_1}$. Setting $f_1=0$ gives $\ds f_5=\frac{1}{6}$, $\ds
 f_6=-\frac{1}{6}$, then $(c_1,c_2,c_3,c_4,c_5=c_6)$ is a permutation of $(0,1,2,3,4)$ and
 we obtain the 6-tuple  (\ref{exp6tuple}).

 If $n=5$, then  $f_1=f_2=f_3=f_4=f_5$, $(c_1,c_2,c_3,c_4,c_5)$
 is a permutation of $(0,1,2,3,4)$, which leads to the 5-tuple
 (\ref{exp5tuple}).\\
 (b.2) Now every $\varphi_k$ can only be equal to
 $\ds 0$, $\ds\frac12$, $\ds\pm\frac13$ or $\ds\pm\frac16$. Direct check shows
  that the only irreducible $n$-tuples with $n\leq6$ that can be built from such numbers are
 (equivalent to) the triple $\ds\left(0,\frac13,-\frac13\right)$ and
 the pair $\ds\left(0,\frac12\right)$.
 \epf

  We now establish a similar classification of rational solutions
  of  (\ref{sumcos}):
 \begin{lem}\label{lemcos}
   Inequivalent irreducible rational $n$-tuples solving (\ref{sumcos})  with $1<n\leq6$
 fall into one of the following classes:
 \begin{itemize}
 \item 13 nontrivial irreducible 6-tuples
 \begin{align*}
 \label{vi1}&\displaystyle\left(\frac{1}{11},\frac{2}{11},\frac{3}{11},\frac{4}{11},\frac{5}{11},\frac{1}{6}\right),&\tag{$\text{VI}_{\text{1}}$}\\
 \label{vi2} &\displaystyle\left(\frac{L}{7}+\frac{1}{6},\frac{L}{7}-\frac{1}{6},\frac{2L}{7},\frac{3L}{7},0,\frac{1}{3}\right),\qquad
 &L=1,2,3,\tag{$\text{VI}_{\text{2}}$}\\
 \label{vi3}
&\displaystyle\left(\frac{L}{7}+\frac{1}{6},\frac{L}{7}-\frac{1}{6},\frac{2L}{7},\frac{3L}{7},\frac{1}{10},\frac{3}{10}\right),\qquad
 &L=1,2,3,\tag{$\text{VI}_{\text{3}}$}\\
 \label{vi4} &\displaystyle\left(\frac{L}{7}+\frac{1}{6},\frac{L}{7}-\frac{1}{6},\frac{2L}{7}+\frac{1}{6},\frac{2L}{7}-\frac{1}{6},\frac{3L}{7},\frac{1}{6}\right),\qquad
 &L=1,2,3,\tag{$\text{VI}_{\text{4}}$}
 \end{align*}
 \ben\label{vi5}
 \left(\frac{1}{7},\frac{2}{7},\frac{3}{7},0,\frac15,\frac25\right),
 \;%\bcancel
 {%\cancel
 {\left(\frac{1}{7},\frac{2}{7},\frac{3}{7},\frac{1}{15},\frac{4}{15},\frac{3}{10}\right)}},
 \;%\bcancel
 {%\cancel
 {\left(\frac{1}{7},\frac{2}{7},\frac{3}{7},\frac{1}{10},\frac{2}{15},\frac{7}{15}\right)}},\tag{$\text{VI}_{\text{5}}$}
  \ebn
 and an infinite family of the form
 \ben\label{via}
 \left(\varphi+\frac16,\varphi-\frac16,\varphi+\frac15,\varphi+\frac25,\varphi+\frac35,\varphi+\frac45\right),\qquad
 \varphi\in\mathbb{Q},\tag{$\text{VI}_{\varphi}$}
 \ebn
 \item 7 nontrivial irreducible 5-tuples
  \ben\label{v1}
  \left(0,\frac{1}{30},\frac{1}{3},\frac{11}{30},\frac{2}{5}\right),\qquad
  \left(0,\frac{1}{5},\frac{7}{30},\frac{1}{3},\frac{13}{30}\right),\tag{$\text{V}_{\text{1}}$}\ebn
  \ben\label{v2}
  \left(\frac{L}{7}+\frac{1}{6},\frac{L}{7}-\frac{1}{6},\frac{2L}{7},\frac{3L}{7},\frac{1}{6}\right),\qquad
  L=1,2,3,\tag{$\text{V}_{\text{2}}$}
  \ebn
  \ben\label{v3}
  \left(\frac{1}{7},\frac{2}{7},\frac{3}{7},0,\frac{1}{3}\right),\qquad
  %\bcancel
  {%\cancel
  {\left(\frac{1}{7},\frac{2}{7},\frac{3}{7},\frac{1}{10},\frac{3}{10}\right)}},\tag{$\text{V}_{\text{3}}$}
  \ebn
 and an infinite family of the form
 \ben\label{va}
 \left(\varphi,\varphi+\frac15,\varphi+\frac25,\varphi+\frac35,\varphi+\frac45\right),\qquad
 \varphi\in\mathbb{Q},\tag{$\text{V}_{\varphi}$}
 \ebn
 \item 4 nontrivial irreducible quadruples
 \ben\label{iv}
 \left(0,\frac{1}{5},\frac{1}{3},\frac{2}{5}\right),\;\left(\frac{1}{30},\frac{1}{6},\frac{11}{30},\frac{2}{5}\right),\;
 \left(\frac{1}{15},\frac{4}{15},\frac{3}{10},\frac{1}{3}\right),\;
 \left(\frac{1}{7},\frac{2}{7},\frac{3}{7},\frac{1}{6}\right),\tag{IV}
 \ebn
 \item 1 nontrivial irreducible triple
 \ben\label{iii1}
 \left(\frac{1}{10},\frac{3}{10},\frac{1}{3}\right)\tag{$\text{III}_{\text{1}}$}
 \ebn
 and an infinite family of the form
 \ben\label{iiia}
 \left(\varphi,\varphi+\frac{1}{3},\varphi-\frac{1}{3}\right),\qquad
 \varphi\in\mathbb{Q},\tag{$\text{III}_{\varphi}$}
 \ebn
 \item an infinite family of pairs of the form
 \ben\label{iia}
 \left(\varphi,\frac{1}{2}-\varphi\right),\qquad\varphi\in\mathbb{Q}.\tag{$\text{II}_{\varphi}$}
 \ebn
 \end{itemize}
 \end{lem}
 \pf
 We use the same ideas, notations and conventions as in  the proof
 of Lemma~\ref{lemexp}. One modification concerns the functions $g_k(x)$ which are now defined by
  \ben
 g_k(x)=\begin{cases}
 \frac12\left[e^{2\pi i f_k} x^{c_k p^{l_1-l_k}}+e^{-2\pi i f_k} x^{p^{l_1}-c_k
 p^{l_1-l_k}}\right] & \text{if}\; c_k\neq0, \\
 \qquad \cos 2\pi \varphi_k & \text{if}\; c_k=0,
 \end{cases}
 \ebn
 As $g_k\left(\exp\left(\frac{2\pi
 i}{p^{l_1}}\right)\right)=\cos2\pi\varphi_k$, one has again $U\left(\exp\left(\frac{2\pi
 i}{p^{l_1}}\right)\right)=0$, so that either (a) $U(x)\equiv
 0$ or (b) $U(x)\not\equiv 0$ is divisible by $P(x)$.\vspace{0.1cm} \\
 \underline{Case (a)}. All $2n$ powers $c_kp^{l_1-l_k}$,
 $p^{l_1}-c_kp^{l_1-l_k}$, appearing in the functions
 $g_k(x)$,
 are simultaneously divisible or non-divisible by $p$ unless we
 have a reducible $n$-tuple. Since $c_1$ is prime to $p$, they
 are actually non-divisible, which in turn gives $l_k=l_1$ for any
 $k$. Irreducibility then implies that
 $c_k$ can only be equal to $c_1$ or $p^{l_1}-c_1$.  In fact we
 can assume that $c_k=c_1$, as the transformation $\varphi_k\mapsto 1-\varphi_k$ maps $f_k\mapsto -f_k$,
 $c_k\mapsto p^{l_k}-c_k$. Now one has
 \ben
 U(x)=\frac12\, x^{c_1}\sum_{k=1}^n e^{2\pi i f_k}+\frac12\, x^{p^{l_1}-c_1}\sum_{k=1}^n e^{-2\pi i
 f_k}=0,
 \ebn
 and, since  $c_1\neq p^{l_1}-c_1$ except in
 the trivial case $p=2$, $l_1=1$, the problem is reduced to
 the classification of rational solutions of the equation (\ref{eqexp}), given
 by Lemma~\ref{lemexp}.
 \vspace{0.1cm}\\
 \underline{Case (b)}. Set $U(x)=P(x)Q(x)$, then by the same reasoning as above $N_U=pN_Q$.
 However, here $N_U\leq 12$, therefore $p$ can be equal to 2, 3, 5, 7 or
 11.

 $2n$ powers $c_k p^{l_1-l_k}$, $p^{l_1}-c_k
 p^{l_1-l_k}$ are all equal modulo $p^{l_1-1}$ to
 $s$ or $p^{l_1-1}-s$, where the integer $s$ does not depend on $k$, $0\leq s<p^{l_1-1}$.
 Otherwise one could collect powers corresponding to different $s$ and write $U(x)$ as
 a sum of at least two polynomials, each of them either divisible by $P(x)$ or
 vanishing. Since $p^{l_1}-c_k
 p^{l_1-l_k}=-c_k p^{l_1-l_k}\;\mathrm{mod}\; p^{l_1-1}$, both terms coming
 from a given $g_k(x)$ will appear in the same polynomial, and then the
 corresponding $n$-tuple is reducible. Hence
 $N_Q$ can only be equal to 1 or 2.

 If $l_1\geq2$, then all $2n$ powers of $x$ that appear in the functions
 $g_k(x)$ are not divisible by $p$ and therefore $l_k=l_1$ for any $k$.

 Two powers $c_k$ and
 $p^{l_1}-c_k$ are distinct modulo $p^{l_1-1}$ for all but a finite number of values of $l_1$ and $p$.
 Indeed, if they are the same, one has $2c_k=0\;\mathrm{mod}\; p^{l_1-1}$. However, this
 is impossible  for $p\geq3$, $l_1\geq2$ and for $p=2$, $l_1\geq3$, since all $c_k$ are prime to
 $p$. Let us now consider separately two cases:
 \begin{itemize}
 \item[(b.1)]  $p\geq3$, $l_1\geq2$ or $p=2$, $l_1\geq3$;
 \item[(b.2)] $p=3,5,7,11$, $l_1=1$ or $p=2$, $l_1=1, 2$.
 \end{itemize}

 (b.1) When $c_k\neq p^{l_1}-c_k\;\mathrm{mod}\;p^{l_1-1}
 $, we use $N_Q\leq 2$ to write the relation
 $U(x)=P(x)Q(x)$ as two distinct equations containing different ($\mathrm{mod}\; p^{l_1-1}$) powers of $x$. Replacing
 $\varphi_k\mapsto 1-\varphi_k$ if necessary, one finds that both
 equations are equivalent to the following one:
 \be\label{condition1}
 \sum_{j=1}^n e^{2\pi i f_k} x^{c_k}=\alpha x^s P(x),\qquad
 \alpha\neq0.
 \eb
 Assume that $n=6$. It is impossible to satisfy (\ref{condition1}) if
 $p=7,11$. For $p=5$ four out of six phases are equal, say
 $f_1=f_2=f_3=f_4$, and the remaining two satisfy
 \begin{itemize}
 \item[(b.1.1)]
 $\displaystyle e^{2\pi i f_5}+e^{2\pi i f_6}=e^{2\pi if_1}$.
 \end{itemize}
 In addition we have $c_k=s+N_k \cdot 5^{l_1-1}$, where
 $(N_1,N_2,N_3,N_4,N_5=N_6)$ is a permutation of $(0,1,2,3,4)$.
 Now applying Lemma~\ref{lemexp} to find rational solutions of
 (b.1.1) we see that resulting 6-tuples are of type (\ref{via}).

 For $p=3$, up to permutations there are only three possibilities:
 \begin{itemize}
 \item[(b.1.2)]
 $\displaystyle e^{2\pi
 if_1}=e^{2\pi
 if_2}= e^{2\pi i f_3}+e^{2\pi i f_4}+e^{2\pi i f_5}+e^{2\pi i f_6}$,
 \item[(b.1.3)] $\displaystyle  e^{2\pi i f_1}=e^{2\pi i f_2}+e^{2\pi i f_3}=
 e^{2\pi i f_4}+e^{2\pi i f_5}+e^{2\pi i f_6}$,
 \item[(b.1.4)] $\displaystyle  e^{2\pi i f_1}+e^{2\pi i f_2}=e^{2\pi i
 f_3}+
 e^{2\pi i f_4}=e^{2\pi i f_5}+e^{2\pi i f_6}\neq0$.
 \end{itemize}
 Finally, for $p=2$ one should have one of the following:
 \begin{itemize}
 \item[(b.1.5)]
 $\displaystyle e^{2\pi
 if_1}=e^{2\pi
 if_2}+ e^{2\pi i f_3}+e^{2\pi i f_4}+e^{2\pi i f_5}+e^{2\pi i f_6}$,
 \item[(b.1.6)] $\displaystyle  e^{2\pi i f_1}+e^{2\pi i f_2}=e^{2\pi i
 f_3}+
 e^{2\pi i f_4}+e^{2\pi i f_5}+e^{2\pi i f_6}\neq0$,
 \item[(b.1.7)] $\displaystyle  e^{2\pi i f_1}+e^{2\pi i f_2}+e^{2\pi i
 f_3}=
 e^{2\pi i f_4}+e^{2\pi i f_5}+e^{2\pi i f_6}\neq0$.
 \end{itemize}
 In each of these cases the problem is reduced to  Lemma~\ref{lemexp}. The 6-tuples we obtain at the end
 turn out to be reducible or belong to the family (\ref{via}).

 Other possibilities ($n=3,4,5$) can be treated in a similar
 manner. They lead to  5-tuples of type (\ref{va}) and triples of
 type (\ref{iiia}).\vspace{0.2cm}

 (b.2) We first consider the case when the denominator of every $\varphi_k$  ($k=1,\ldots,n$) is not
 divisible by $7$ and $11$:
   \begin{lem}\label{lemma60}
 Inequivalent irreducible $n$-tuples solving (\ref{sumcos}) with $3\leq n\leq6$ such that
 every $d_k$ ($k=1,\ldots,n$) is a divisor of $2^2\cdot3\cdot5=60$ are given by
 \begin{itemize}
 \item 6-tuples:
 \ben
 \left(0,\frac{1}{30},\frac{1}{5},\frac{11}{30},\frac{2}{5},\frac{2}{5}\right),\quad
 \left(0,\frac{1}{30},\frac{7}{30},\frac{1}{3},\frac{11}{30},\frac{13}{30}\right),\quad
 \left(0,\frac{1}{5},\frac{1}{5},\frac{7}{30},\frac{2}{5},\frac{13}{30}\right),
 \ebn
 \ben
 \left(\frac{1}{60},\frac{1}{60},\frac{13}{60},\frac{7}{20},\frac{23}{60},\frac{5}{12}\right),\;
 \left(\frac{1}{60},\frac{1}{20},\frac{11}{60},\frac{23}{60},\frac{23}{60},\frac{5}{12}\right),\;
  \left(\frac{1}{60},\frac{11}{60},\frac{13}{60},\frac{13}{60},\frac{5}{12},\frac{9}{20}\right),
 \ebn
 \ben
 \left(\frac{1}{12},\frac{7}{60},\frac{17}{60},\frac{19}{60},\frac{19}{60},\frac{7}{20}\right).
 \ebn
 \item 5-tuples:
 \ben
  \left(0,\frac{1}{5},\frac{1}{5},\frac{2}{5},\frac{2}{5}\right),\qquad
  \left(\frac{1}{60},\frac{11}{60},\frac{13}{60},\frac{23}{60},\frac{5}{12}\right),\qquad
  \left(\frac{1}{30},\frac{1}{6},\frac{7}{30},\frac{11}{30},\frac{13}{30}\right),
 \ebn
 \ben
 \left(0,\frac{1}{30},\frac{1}{3},\frac{11}{30},\frac{2}{5}\right),\qquad
 \left(0,\frac{1}{5},\frac{7}{30},\frac{1}{3},\frac{13}{30}\right).
 \ebn
 \item quadruples:
 \ben
 \left(0,\frac{1}{5},\frac{1}{3},\frac{2}{5}\right),\qquad
 \left(\frac{1}{30},\frac{1}{6},\frac{11}{30},\frac{2}{5}\right),\qquad
 \left(\frac{1}{15},\frac{4}{15},\frac{3}{10},\frac{1}{3}\right).
 \ebn
 \item triples:
 \ben
 \left(0,\frac{1}{3},\frac{1}{3}\right),\quad
 \left(\frac{1}{60},\frac{19}{60},\frac{7}{20}\right),\quad
 \left(\frac{1}{30},\frac{3}{10},\frac{11}{30}\right),\quad
 \left(\frac{1}{20},\frac{17}{60},\frac{23}{60}\right),\ebn
 \ben
 \left(\frac{1}{15},\frac{4}{15},\frac{2}{5}\right),\quad
 \left(\frac{1}{10},\frac{3}{10},\frac{1}{3}\right).
 \ebn
 \end{itemize}
 \end{lem}
 \pf Direct (e.g., Mathematica) computation. Notice that all obtained $6$-tuples, first three $5$-tuples
 and all but the last triple belong to the infinite families (\ref{via}), (\ref{va}) and (\ref{iiia}), respectively.
 \epf\vspace{0.1cm}

 The case $p=11$, $l_1=1$ is possible only for $n=6$. We
 have $N_{Q}=1$, $\mathrm{deg}\,Q=0$, hence $Q(x)=\alpha$,
 $N_U=11$, $\mathrm{deg}\,U=10$ and, consequently, one can choose
 $l_1=\ldots=l_5=1$, $l_6=0$, $c_k=k$ ($k=1,\ldots,5$),
 $c_6=0$. This gives the irreducible $6$-tuple (\ref{vi1}).

 Remaining case $p=7$, $l_1=1$ is possible only for $n=4,5,6$. Similarly
 to the above, $N_Q=1$, $\mathrm{deg}\,Q=0$, $Q(x)=\alpha$,
 $N_U=7$, $\mathrm{deg}\,U=6$, and in addition for all $k=2,\ldots, n$
 either $l_k=1$ or $c_k=0$. For $n=6$ one then has four
 possibilities:
 \begin{itemize}
  \item $(c_1,c_2,c_3=c_4)$ is a permutation of $(1,2,3)$,
 $c_5=c_6=0$; this gives $f_1=f_2=0$ and
 \be\label{7a}
 e^{2\pi if_3}+e^{2\pi i f_4}=2\cos2\pi f_5+2\cos2\pi f_6=1.
 \eb
 Recall that $f_1,\ldots,f_6$ are rational numbers with
 denominator which is a divisor of $60$. Using
 Lemma~\ref{lemma60}  to classify the appropriate solutions
 of (\ref{7a}), one finds that the only irreducible $6$-tuples
 obtained in this way are given by (\ref{vi2}) and (\ref{vi3}).
  \item $(c_1,c_2=c_3,c_4=c_5)$ is a permutation of $(1,2,3)$,
 $c_6=0$; then
 \ben
 f_1=0,\qquad e^{2\pi if_2}+e^{2\pi i f_3}=e^{2\pi if_4}+e^{2\pi i
 f_5}=2\cos2\pi f_6=1,
 \ebn
 which leads to the family of irreducible $6$-tuples (\ref{vi4}).
 \item $(c_1,c_2,c_3=c_4=c_5)$ is a permutation of $(1,2,3)$,
 $c_6=0$; then $f_1=f_2=0$ and
 \ben
 e^{2\pi if_3}+e^{2\pi i f_4}+e^{2\pi i f_5}=2\cos2\pi f_6=1.
 \ebn
 All $6$-tuples arising here turn out to be reducible.
 \item $(c_1,c_2,c_3)=(1,2,3)$, $c_4=c_5=c_6=0$, which implies
 $f_1=f_2=f_3=0$ and
 \be\label{7b}
 2\cos2\pi f_4+2\cos2\pi f_5+2\cos2\pi f_6=1.
 \eb
 Using again Lemma~\ref{lemma60} to find irreducible solutions of
 (\ref{7b}), we obtain 3 irreducible 6-tuples (\ref{vi5}).
 \end{itemize}
 For $n=5$, there are two possibilities:
 \begin{itemize}
 \item $(c_1,c_2,c_3=c_4)$ is a permutation of $(1,2,3)$, $c_5=0$;
 this implies $f_1=f_2=0$ and
 \ben
 e^{2\pi if_3}+e^{2\pi i f_4}=2\cos2\pi f_5=1,
 \ebn
 so that we find 3 irreducible $5$-tuples (\ref{v2}).
 \item $(c_1,c_2,c_3)=(1,2,3)$, $c_4=c_5=0$,
 hence $f_1=f_2=f_3=0$ and
 \ben
 2\cos2\pi f_4+2\cos2\pi f_5=1.
 \ebn
 This gives 2 irreducible 5-tuples (\ref{v3}).
 \end{itemize}
 Finally, for $n=4$ we should have $(c_1,c_2,c_3)=(1,2,3)$, $c_4=0$ and,
 therefore, $f_1=f_2=f_3=0$, $2\cos2\pi f_4=1$, which leads to the
 fourth irreducible quadruple in (\ref{iv}). This concludes the
 proof of Lemma~\ref{lemcos}.
 \epf
 \begin{rmk} The classification of irreducible rational solutions
 of (\ref{sumcos}) with $n\leq 4$ is essentially equivalent to
 Lemma~1.13 in \cite{dubrovin}. In fact we will see shortly that
 this partial result is already sufficient to find all finite $\bar{\Lambda}$~orbits
 with $\omx^2\neq\omy^2\neq\omz^2$. Its extension to $n=5,6$
 is needed to treat the case when $\boldsymbol{\omega}\in\Cb^3$ is
 fixed by some of the $K_4\rtimes S_3$ transformations.
 \end{rmk}
 \subsection{Bounds on suborbit lengths} Let $O\subset\Cb^3$ be a
 finite orbit of the induced $\bar{\Lambda}$ action~(\ref{lxyz}).
 We choose an arbitrary 2-colored  suborbit $O_{yz}\subset O$ (i.e.
 the suborbit generated from a given point by two
 transformations $y$ and $z$), denote its length by $N$ and label the points of $O_{yz}$ as in
 Subsection~\ref{ss2c}.

 Throughout this subsection we assume that $N>1$. Denote $X=2\cos\lambda/2$, then by Lemma~\ref{orbcosines}
 one has $\lambda=2\pi r_X$, $r_X=n_X/N$, where $n_X\in\Zb$ is prime to $N$ and
 we choose $0<n_X<N$. Lemma~\ref{yzsolutions} implies in addition
 that $Y_k$, $Z_k$ ($k=0,1,\ldots,N-1$) are given by (\ref{yzsolution1}).

 When the graph of $O_{yz}$ is a simple cycle, it
 contains $2N$ points and all of them are good. Then by
 Lemma~\ref{xyzcosineslemma}
 for $k=0,\ldots,N-1$ we have
 \be\label{yzgood}
 Y_k=2\cos\pi r_{Y_k},\qquad  Z_k=2\cos\pi r_{Z_k},\qquad\qquad
 r_{Y_k}, r_{Z_k}\in\Qb,\quad 0<r_{Y_k}, r_{Z_k}<1.
 \eb
 If $\Sigma(O_{yz})$ is a line with self-loops at the
 ends, then there are $N$ distinct points. While two endpoints can in principle be bad,
 the other $N-2$ points are good so that their coordinates satisfy (\ref{yzgood}).
 \begin{lem}\label{distinct}
 Two distinct vertices of $\Sigma(O_{yz})$ characterized by the same coordinate $Y$ (or $Z$) are
 necessarily connected by an edge of color $z$ (resp. $y$).
 \end{lem}
 \pf Let $(X,Y,Z)$ be an arbitrary point in $O$. Since $\omega_{X,Y,Z,4}$ are fixed by the $\bar{\Lambda}$ action, the quantity
 \ben
 XYZ+X^2+Y^2+Z^2-\omx X-\omy Y-\omz Z =\mathrm{const}= 4-\omi
 \ebn
 is an orbit invariant. Computing this invariant for two distinct points $(X,Y,Z)$, $(X,Y,Z')$ in $O_{yz}$
 we find $Z'=\omz-Z-XY=z(Z)$. \epf
 \begin{rmk}
 In the simple cycle case, Lemma~\ref{distinct} implies that
 $Y_k\neq Y_{k'}$, $Z_k\neq Z_{k'}$ for $k\neq k'$ where
 $k,k'=0,\ldots,N-1$. Similarly, in the line case for any
 $k$ there exists at most one $k'\neq k$ such that $Y_k=Y_{k'}$
 (or $Z_k=Z_{k'}$).
 \end{rmk}
 \begin{lem}
 The coordinates $\{Y_k\}$, $\{Z_k\}$ satisfy the following
 identities:
 \begin{align}
 \label{evenodd1}
  &\text{for $N$ even, $n_X$ odd:}
  && \begin{cases}
  Y_k+Y_{k+N/2}=p_++p_-\,, \\
  Z_k+Z_{k+N/2}=p_+-p_-\,,
  \end{cases}  \\
  \label{evenodd2}
  &\text{for $N$ odd, $n_X$ even:} &&
 \quad Y_k+Z_{k+(N-1)/2}=p_+\, , \\
 \label{evenodd3}
  &\text{for $N$ odd, $n_X$ odd:} &&
 \quad Y_k-Z_{k+(N-1)/2}=p_-\, ,
 \end{align}
 where $k=0,\ldots,N-1$ and $\ds p_{\pm}=\frac{\omy\pm\omz}{2\pm X}$.
 \end{lem}
 \pf Straightforward substitution of (\ref{yzsolution1}) into (\ref{evenodd1})--(\ref{evenodd3}). \epf
 \begin{prop}\label{evenbound}
 If $N$ is even and at least one of two parameters
 $\omy$, $\omz$ is different from $0$, then $N\leq 10$.
 \end{prop}
 \pf When at least one of $\omy$, $\omz$ differs from $0$, at least one of $p_+\pm p_-$ is also non-zero.
 Assume for definiteness that $p_++p_-\neq0$ and consider the
 first equation in (\ref{evenodd1}). It implies that for any
 $k,k'=0,\ldots,N-1$ one has
 \be\label{yneven}
 Y_k+Y_{k+N/2}=Y_{k'}+Y_{k'+N/2}\neq0.
 \eb

 First assume that the graph of
 $O_{yz}$ is a simple cycle. All $Y_k$ are then distinct and have
 the form (\ref{yzgood}). Hence  (\ref{yneven}) reduces to an equation
 of type (\ref{sumcos}) with $n=4$, whose rational solutions have been classified
 in Lemma~\ref{lemcos}. We now consider
 different types of solutions to maximize the number $N^c$ of possible unordered couples
 $\bigl(Y_k,Y_{k+N/2}\bigr)$ of the form
 (\ref{yzgood}), characterized by the same value of $Y_k+Y_{k+N/2}$:
 \begin{itemize}
 \item Splitting of the rational solution quadruple into two (not necessarily irreducible) pairs
 is possible only for $k'=k$ or $k'=k+N/2$, therefore one should not
 take such solutions into account when computing $N^c$ (here we used that
 $p_++p_-\neq0$~!).
 \item Assume that $Y_{k_0}=0$ for some $k_0$, then for any
 $k$ one has $Y_k+Y_{k+N/2}=Y_{k_0+N/2}$. This is an equation of
 type (\ref{sumcos}) with $n=3$. By Lemma~\ref{lemcos},
 if $Y_{k_0+N/2}\neq \pm1$, $\pm 2\cos \pi/5$, $\pm 2\cos2\pi/5$, the only possible couple different from
 $\bigl(0,Y_{k_0+N/2}\bigr)$ is
 \ben
 \bigl(2\cos\pi(r_{Y_{k_0+N/2}}+1/3),2\cos\pi(r_{Y_{k_0+N/2}}-1/3)\bigr)
 \ebn
 and therefore $N^c=2$. When $Y_{k_0+N/2}= \pm 1$, the only compatible couple
 is $\bigl(\pm 2\cos\pi/5,\mp2\cos 2\pi/5\bigr)$ so that again
 $N^c=2$.

 Finally, for (a) $Y_{k_0+N/2}= \pm 2\cos\pi/5$ and (b) $Y_{k_0+N/2}=\pm
 2\cos2\pi/5$ one has $N^c=3$ as in both cases we have three compatible couples:\\
 (a) $\bigl(0,\pm2\cos\pi/5\bigr)$, $\bigl(\pm1,\pm2\cos2\pi/5\bigr)$,
 $\bigl(\pm 2\cos2\pi/15,\pm2\cos8\pi/15\bigr)$;\\
 (b)  $\bigl(0,\pm2\cos2\pi/5\bigr)$, $\bigl(\mp1,\pm2\cos\pi/5\bigr)$, $\bigl(\pm 2\cos\pi/15,\pm2\cos11\pi/15\bigr)$.
 \item If there is no $Y_k$ equal to zero, the solution quadruple
 can only be equivalent to one of the last 3 quadruples in
 (\ref{iv}) (first quadruple is excluded because $Y_k\neq\pm2$).
 Direct check then shows that for any choice of  $\bigl(Y_k,Y_{k+N/2}\bigr)$ there is only one
 compatible couple, i.e. $N^c=2$.
 \end{itemize}
 Since the maximal possible value of $N^c$ is $3$, even length $N$
 of the simple cycle cannot exceed~$6$.

 When the graph of $O_{yz}$ is a line, the same reasoning shows that $N\leq14$, otherwise
 the number of distinct compatible couples
 $\bigl(Y_k,Y_{k+N/2}\bigr)$ satisfying (\ref{yzgood}) is greater than
 3. We now want to improve this bound to $N\leq10$ using that for $N=12,14$ the number of such
 couples is $3$ and therefore $Y$-coordinates of good points should give (a) or (b) above.
 \begin{figure}[!h]
 \begin{center}
 \resizebox{13cm}{!}{
 \includegraphics{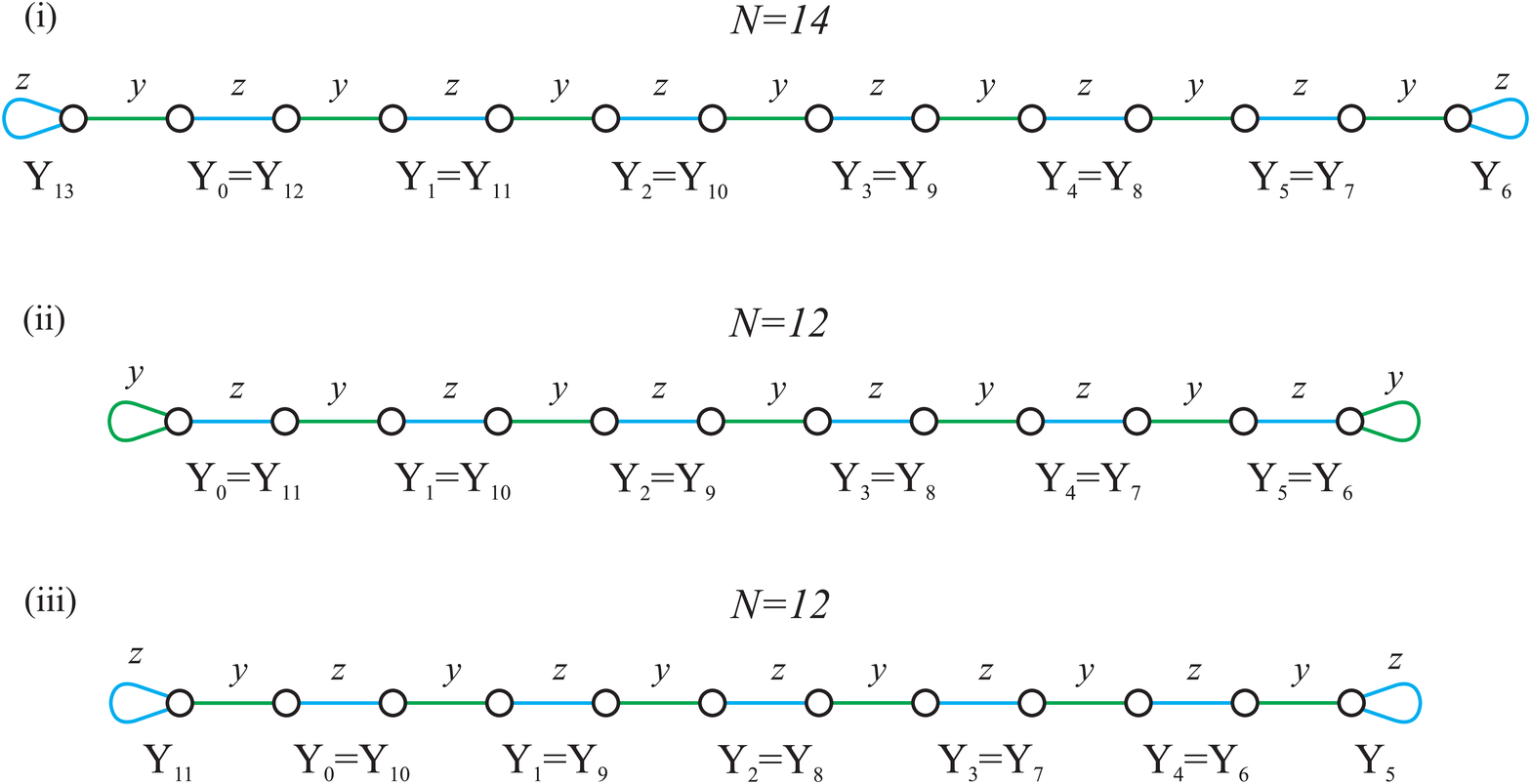}} \\   Fig. 2: Three possible graphs for $N=12,14$
 \end{center}
 \end{figure}

 In Fig.~2 we show three possible graphs and label each vertex by
 its $Y$-coordinate. Third diagram (iii) can in fact be immediately
 excluded, since in this case $2Y_2=Y_1+Y_3=Y_0+Y_4$ but no couple in (a) or
 (b) contains two equal cosines. To exclude the remaining two cases,
 use that from (\ref{yzrecursion}) follows a 2nd order
 difference equation for $\{Y_k\}$:
 \ben
 Y_{k+2}+(2-X^2)Y_{k+1}+Y_k=2\omy-X\omz.
 \ebn
 It implies in particular that for both (i) and (ii) we should
 have
 \be\label{x1314}
 X^2-1=\frac{Y_4-Y_1}{Y_3-Y_2}.
 \eb
 Since $(Y_1,Y_4)$ and $(Y_2,Y_3)$ are necessarily given by two couples
 from (a) or (b), the RHS of (\ref{x1314}) can only take one of 12 values
 \ben
 \varepsilon_1 (\sqrt{5}+2\varepsilon_2),\qquad \varepsilon_1(15+6\varepsilon_2
 \sqrt{5})^{\varepsilon_3/2},\qquad\qquad
 \varepsilon_{1,2,3}=\pm1.
 \ebn
 Possible values of the LHS also belong to an explicitly defined finite set: recall that $X=2\cos\pi
 n_X/N$, where $n_X=1,3,5,9,11$ or $13$ for $N=14$ and $n_X=1,5,7$ or $11$ for
 $N=12$. Now it is easy to check that the LHS and the RHS of (\ref{x1314})
 never match, and thus the lengths $N=12,14$ are forbidden.
  \epf
 \begin{prop}\label{oddbound1}
 If $N$ is odd and $\omy^2\neq\omz^2$, then $N\leq 9$.
 \end{prop}
 \pf The condition $\omy^2\neq\omz^2$ guarantees that both $p_+$ and $p_-$ are non-zero.
 Assuming for definiteness that $n_X$ is odd, one finds from (\ref{evenodd3})
 \ben
 Y_k-Z_{k+(N-1)/2}=Y_{k'}-Z_{k'+(N-1)/2}\neq 0.
 \ebn
 We can now use the same approach as in the previous proof. One difference is that here we
 maximize the number of \textit{ordered} couples $\bigl(Y_k,Z_{k+(N-1)/2}\bigr)$
 of the form (\ref{yzgood})
 characterized by the same value of $Y_k-Z_{k+(N-1)/2}$. This maximal
 number is equal to~$6$ (twice the maximal $N^c$),
 therefore by Lemma~\ref{distinct} simple cycles of length
 $N\geq7$ and the lines of length $N\geq15$ are forbidden.

 The lengths $N=11,13$ are excluded similarly to the above, since in
 this case $Y$- and $Z$-coordinates of good points take only a finite
 number of explicitly defined values. Straightforward computation shows that
 possible values of $X$ determined from (\ref{yzrecursion}) never
 match $X=2\cos \pi n_X/N$.
  \epf
 \begin{rmk} In the proof of Proposition~\ref{oddbound1}
 we used only that $p_-\neq0$. Therefore the bound ``$\text{odd }N\leq9$''  also holds for
 $\omy=\omz\neq0$ when $n_X$ is even and for $\omy=-\omz\neq0$ when $n_X$ is odd.
 \end{rmk}

 Next we study the case $\omy=\omz$, $n_X$ odd, where the relation (\ref{evenodd3}) gives
 just $Y_k=Z_{k+(N-1)/2}$. For $\omy=-\omz$,
 $n_X$ even the upper bound for $N$ is the same
 by symmetry; recall that e.g. the transformation $\omx\mapsto -\omx$, $\omy\mapsto -\omy$,
 $(X,Y,Z)\mapsto(-X,-Y,Z)$ for all
 $(X,Y,Z)\in O$ yields an orbit equivalent to $O$.
 \begin{prop}\label{oddbound2}
 Let $N$ and $n_X$ be odd and let $\omy=\omz\neq0$. If the graph $\Sigma(O_{yz})$ is
 a line, then the only possible values of $N$ are $3, 5, 7, 9, 11,
 15$ and $21$.
 \end{prop}
 \pf The suborbit graph for odd $N$ is presented in Fig.~3. Each vertex
 is labeled by its coordinates $(Y,Z)$. For $\omy=\omz$ one has $p_-=0$, hence (\ref{evenodd3})
 implies in particular that for the center point  $Z=Y$.
 \begin{figure}[!h]
 \begin{center}
 \resizebox{12cm}{!}{
 \includegraphics{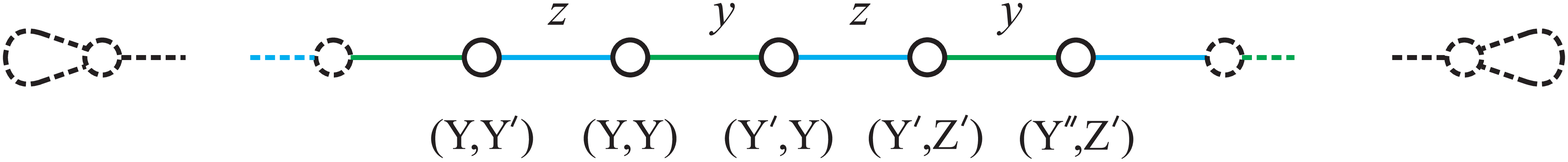}} \\  Fig. 3: Line of odd length, $\omy=\omz$
 \end{center}
 \end{figure}

 Let us denote $\omy=\omz=\omega$ and
 $X=2\cos\pi r_X$, $Y=2\cos\pi r_{Y}$, $Y'=2\cos\pi r_{Y'}$ etc.
 From the relations
 \ben
 Y+Y'+XY=\omega=Y+Z'+XY'
 \ebn
 one finds an equation of type (\ref{sumcos}) with $n=6$:
 \begin{align}
 \nonumber &\cos\pi r_{Y'}+\cos\pi(r_X-r_Y)+\cos\pi(r_X+r_Y)=\\
 \label{oddboundaux1} &\qquad\qquad=\cos\pi r_{Z'}+
 \cos\pi(r_X-r_{Y'})+\cos\pi(r_X+r_{Y'}).
 \end{align}
 We assume that $N\geq7$, then $r_{X,Y,Y',Z'}\in\Qb$ by Lemma~\ref{xyzcosineslemma}.

 General idea of the proof
 is to obtain the restrictions on $r_X$ from Lemma~\ref{lemcos}.
 Not all solutions listed in Lemma~\ref{lemcos} are of interest
 here because the arguments of cosines in (\ref{oddboundaux1}) are
 not all independent. Five
 entries in the solution $6$-tuple, say
 $\varphi_1\ldots\varphi_5$, should satisfy\\
 (a) $\varepsilon_1\varphi_1+\varepsilon_2\varphi_2+\varepsilon_3\varphi_3+\varepsilon_4\varphi_4\in\Zb$
 for some choice of $\varepsilon_{1,2,3,4}=\pm1$.\\
 (b) $\varepsilon_3\varphi_3-\varepsilon_4\varphi_4=2\varepsilon_5\varphi_5\;(\mathrm{mod}\;\Zb)$
 for the same $\varepsilon_{3,4}$ and some
 $\varepsilon_5=\pm1$.
 \begin{rmk} In many cases below, the number of possible solutions for $r_X$ is rather large
 and their complete description becomes too cumbersome.
 However, since $r_X=n_X/N$ and $N$ is odd, in practice it is easy to
 determine admissible values of $N$
 by simply looking at odd integers that can appear
 in the denominator of~$r_X$. The reader should keep in mind that probably not all
 such admissible values do actually occur. For clarity, the values $N=3,5$
 (not satisfying the above assumption $N\geq7$) will not be omitted
 in the course of this shortcut computation.
 \end{rmk}

 First assume that the solution of (\ref{oddboundaux1})
 is equivalent to one of the $6$-tuples
 (\ref{vi1})--(\ref{via}):\vspace{0.1cm}\\
 (\ref{vi1}) In this 6-tuple, $1/6$ clearly corresponds to $r_{Z'}$ in
 (\ref{oddboundaux1}), otherwise conditions (a) and (b) cannot be
 simultaneously satisfied. Hence the only possible odd denominator of
 $r_X$ is $11$.\vspace{0.1cm}\\
 (\ref{vi2}) Considering the sum and the difference
 of any two elements in (\ref{vi2}), one readily concludes that
 the only possible odd denominators of $r_X$ are $3, 7$ and
 21.\vspace{0.1cm}\\
 (\ref{vi3}) Condition~(a) fails unless $1/10$ and $3/10$ correspond to $r_{Y'}$ and
 $r_{Z'}$ or vice versa. In both cases, however, (b) is violated.\vspace{0.1cm}\\
 (\ref{vi4}) Possible $N$ are $3, 7, 21$ by the same argument as
 in (\ref{vi2}).
 \vspace{0.1cm}\\
 (\ref{vi5}) With the second and the third 6-tuple condition (a) always
 fails. With the first 6-tuple it can be satisfied only if $1/5$
 and $2/5$ correspond to $r_{Y'}$ and $r_{Z'}$ or vice versa, but
 then (b) is violated.\vspace{0.1cm}\\
 (\ref{via}) Taking the sum and the difference of any two elements
 (meant to be $r_X\pm r_{Y'}$) we see that that  odd
 divisors of the denominator of either $r_X$ or $r_{Y'}$ can only
 be $3, 5, 15$. However, in the second case $\varphi$ becomes fixed
 so that admissible $N$ are again $3,5,15$.\vspace{0.1cm}

 Reducible $6$-tuples consisting of one 5-tuple from
 ($\text{V}_{\text{1}}$)--($\text{V}_{\varphi}$) and one zero cosine
 (we will say that the solution is of type
 ``$\text{V}_{\text{1,2,3},\varphi}+\text{I}$'') can be treated in a
 completely similar manner, leading to $N=3,5,7,15,21$. These values of $N$
 are also the only admissible ones for the solutions of type
 ``$\text{IV}+\text{II}_{\varphi}$'', where the solution 6-tuple
 splits into one of the irreducible quadruples (\ref{iv}) and a pair
 of the form (\ref{iia}). Solutions of type
 ``$\text{III}_{\text{1}}+\text{III}_{\text{1}}$'' and
 ``$\text{III}_{\text{1}}+\text{II}_{\varphi}+\text{I}$'' lead to
 $N=3,5,15$, and those of type
 ``$\text{III}_{\text{1}}+\text{III}_{\varphi}$'' to
 $N=3,5,9,15$. There remain three types of possible rational solution
 6-tuples:

 (1)~``$\text{III}_{\varphi}+\text{II}_{\psi}+\text{I}$'';

 (2)~``$\text{II}_{\varphi}+\text{II}_{\psi}+\text{II}_{\mu}$'';

 (3)~``$\text{III}_{\varphi}+\text{III}_{\psi}$''.\\
 \underline{Case (1)}. We first study the case when (\ref{oddboundaux1})
 contains at least one zero cosine (in particular, this includes (1)). There
 are four inequivalent possibilities:\vspace{0.1cm}\\
 (1.1) Set {$Y'=0$}, then from (\ref{oddboundaux1})
 follows  $XY=Z'$. This equation clearly reduces to
 (\ref{sumcos}) with $n=3$ and $\varphi_{1,2,3}\in\Qb$, hence its
 solutions are described by Lemma~\ref{lemcos}. Solutions equivalent to
 (\ref{iii1}) can lead only to $N=3,5,15$, and it remains to consider
 solutions of type ``\ref{iiia}'' and
 ``$\text{II}_{\varphi}+\text{I}$''.\\
 (1.1.1) Solution of $XY=Z'$ has the form (\ref{iiia}) only if
 $X=\pm1$ (i.e. $N=3$) or $Y=\pm1$. In the latter case $Z'=\pm X$
 and $\omega=\pm(1+X)$. Now computing $Y''=\omega-Y'-XZ'$ we
 find $\cos\pi r_{Y''}=\pm\left(\cos\pi r_X-\cos2\pi
 r_X-\cos{\pi}/{3}\right)$. By virtue of Lemma~\ref{xyzcosineslemma}, for $N\geq9$ one
 has $r_{Y''}\in\Qb$. We can thus apply Lemma~\ref{lemcos} to the
 last relation. Irreducible quadruples (\ref{iv}) lead to
 $N=3,5,7,15$, solutions of type
 ``$\text{III}_{\text{1}}+\text{I}$'' to $N=5$, and solutions of
 type ``$\text{III}_{\varphi}+\text{I}$'' and
 ``$\text{II}_{\varphi}+\text{II}_{\psi}$'' to $N=3$.\\
 (1.1.2) Now consider  solutions of $XY=Z'$ containing at least one zero cosine.
 Note that $Z'\neq0$ for $N\geq7$, since
 by (\ref{evenodd3}) $Y_k=Z_{k+(N-1)/2}$ and we have already put $Y'=0$.
 One can therefore assume that $r_Y=r_X\pm1/2\;(\mathrm{mod}\;2\Zb)$,
 $Z'=2\cos\pi(2r_X\pm1/2)$. Computation of $Y''$ then gives
 \be\label{oddboundaux1a}
 \cos\pi r_{Y''}=\cos\pi (2r_X\pm1/2)-\cos\pi (3r_X\pm1/2).
 \eb
 If $N\geq9$, one can apply to (\ref{oddboundaux1a}) Lemma~\ref{lemcos}.
 Solutions (\ref{iii1}) and (\ref{iiia}) can lead only to $N=3,5$ and
 $N=3,5,15$ correspondingly. Since $Y''\neq0$, the only possible $N$ for solutions of type
 ``$\text{II}_{\varphi}+\text{I}$'' is 3.
 As a consequence, from now on we can assume that
 $Y'\neq0$.\vspace{0.1cm}\\
 (1.2) Suppose that $Z'=0$. Here we will use two relations of the form (\ref{sumcos}). The
 first one, with $n=5$, is merely (\ref{oddboundaux1}) with
 $Z'=0$:
 \be\label{oddboundaux3}
 Y'+XY=XY'.
 \eb
 Recall that we can restrict our attention to solutions of (\ref{oddboundaux3})
 of 2 types: ``$\text{II}_{\varphi}+\text{II}_{\psi}+\text{I}$'' and
 ``$\text{III}_{\varphi}+\text{II}_{\psi}$''.
 The second equation, with $n=4$, comes from the computation of
 $Y''$,
 \be\label{oddboundaux4}
 Y''=Y+XY.
 \eb
 Assume that $N\geq9$ to guarantee $r_{Y''}\in\Qb$ and consider rational solutions
 of (\ref{oddboundaux4}) given by Lemma~\ref{lemcos}. It is easy to check that the quadruples
 equivalent to (\ref{iv}) can only lead to $N=3,5,7,15,21$, while for solutions
 of type ``$\text{III}_{\text{1}}+\text{I}$'' one has $N=3,5,15$.

 Next we examine solutions of (\ref{oddboundaux4}) of type
 ``$\text{III}_{\varphi}+\text{I}$''. Since $Y,Y''\neq 0$ it
 can be assumed that $r_Y=r_X\pm1/2\;(\mathrm{mod}\;2\Zb)$ and then
 the triple (\ref{iiia}) becomes
 \ben
 \cos\pi r_{Y''}=\cos\pi (r_X\pm1/2)+\cos\pi (2r_X\pm1/2),
 \ebn
 giving $N=3,9$. Finally, for solutions of type
 ``$\text{II}_{\varphi}+\text{II}_{\psi}$'', since $Y\neq Y''$, we
 may write
 \ben
 Y''=2\cos\pi (r_Y+r_X),\qquad Y+2\cos\pi(r_Y-r_X)=0.
 \ebn
 Second relation implies that $r_X=2r_Y+1\;(\mathrm{mod}\;2\Zb)$ (remember that
 $X\neq\pm2$). Substituting this into (\ref{oddboundaux3}),
 one finds
 \be\label{oddboundaux5}
 \cos\pi r_{Y'}-\cos\pi r_Y-\cos3\pi
 r_Y+\cos\pi(2r_Y+r_{Y'})+\cos\pi(2r_Y-r_{Y'})=0.
 \eb
 (1.2.1) Now consider solutions of (\ref{oddboundaux5}) of type
 ``$\text{II}_{\varphi}+\text{II}_{\psi}+\text{I}$''. Note that
  $Y,Y'\neq0$. Furthermore $\cos3\pi
 r_Y=0$ implies $N=3$, therefore it may be assumed that
 $\cos\pi(2r_Y-r_{Y'})=0$, i.e. $r_{Y'}=2r_Y\pm1/2\;(\mathrm{mod}\;2\Zb)$. Then (\ref{oddboundaux5})
 transforms into
 \ben
 \cos\pi (2r_Y\pm1/2)-\cos\pi r_Y-\cos3\pi
 r_Y+\cos\pi(4r_Y\pm1/2)=0.
 \ebn
 We are looking for rational solutions of the last relation that
 have type ``$\text{II}_{\varphi}+\text{II}_{\psi}$'', hence the
 only admissible values of $N$ are 3 and 5.\\
 (1.2.2) Consider a solution of (\ref{oddboundaux5}) of type
 ``$\text{III}_{\varphi}+\text{II}_{\psi}$'' and take into account
 the following comments:
 \begin{itemize}
 \item $\cos\pi r_Y$
 and $\cos3\pi r_Y$ cannot belong simultaneously to
 ($\text{II}_{\psi}$) because then the denominators of $r_Y$ and $r_X$ would not
 have odd divisors. They can neither belong simultaneously to
 ($\text{III}_{\varphi}$) unless $N=3$. Therefore we may assume that $\cos\pi r_Y$
 and $\cos3\pi r_Y$ are divided between ($\text{III}_{\varphi}$)
 and ($\text{II}_{\psi}$).
 \item $\cos\pi(2r_Y\pm r_{Y'})$ cannot belong
 simultaneously to ($\text{II}_{\psi}$) as there is no enough place.
 It they are both in ($\text{III}_{\varphi}$) then either $N=3$ or
 $Y'=\pm1$. In the latter case, since $Y'$ belongs to
 ($\text{II}_{\psi}$), one can only have $N=3,9$. Hence it may be
 assumed that  $\cos\pi(2r_Y\pm r_{Y'})$ are divided between ($\text{III}_{\varphi}$)
 and ($\text{II}_{\psi}$), and in particular $Y'$ belongs to
 ($\text{III}_{\varphi}$).
 \end{itemize}
 Then we are left with two inequivalent possibilities:
 \begin{align}
 \begin{cases}
 \cos\pi r_{Y'}-\cos\pi r_{Y}+\cos\pi (2r_{Y}-r_{Y'})=0 & \qquad (\text{III}_{\varphi}) \\
 \cos3\pi r_{Y}=\cos\pi (2r_{Y}+r_{Y'}) & \qquad (\text{II}_{\psi})
 \end{cases} \tag{1.2.2.1}
 \end{align}
 From the second equation one finds either ${Y'}=Y$ (forbidden) or
 $r_{Y'}=-5r_Y\;(\mathrm{mod}\;2\Zb)$. But then the first equation
 transforms into $\cos5\pi r_Y+\cos7\pi r_Y-\cos\pi r_Y=0$, which
 implies $N=3,9$.
  \begin{align}
 \begin{cases}
 \cos\pi r_{Y'}-\cos3\pi r_{Y}+\cos\pi (2r_{Y}+r_{Y'})=0 & \qquad (\text{III}_{\varphi}) \\
 \cos\pi r_{Y}=\cos\pi (2r_{Y}-r_{Y'}) & \qquad (\text{II}_{\psi})
 \end{cases} \tag{1.2.2.2}
 \end{align}
 Again from the second equation follows either $Y'=Y$ or
 $r_{Y'}=3r_Y\;(\mathrm{mod}\;2\Zb)$. In the latter case the
 substitution into the first equation gives $\cos 5\pi r_Y=0$, hence
 the only admissible $N$ is~5.\\
 (1.3) Set $\cos\pi(r_X-r_Y)=0$. This implies
 $r_Y=r_X+\varepsilon_1/2\;(\mathrm{mod}\;2\Zb)$, $\varepsilon_1=\pm1$ and our initial equation
 (\ref{oddboundaux1}) transforms into
 \be\label{oddboundaux6}
 \cos\pi r_{Y'}+\cos\pi(2r_X+\varepsilon_1/2)=\cos\pi r_{Z'}+
 \cos\pi(r_X-r_{Y'})+\cos\pi(r_X+r_{Y'}).
 \eb
 (1.3.1) We first study solutions of (\ref{oddboundaux6}) of type
 ``$\text{II}_{\varphi}+\text{II}_{\psi}+\text{I}$''. All cases
 when $Y'=0$ or $Z'=0$ have been considered above. Moreover
 $\cos\pi(2r_X+\varepsilon_1/2)=0$ would lead only to even $N$,
 therefore it can be assumed that $\cos\pi(r_X-r_{Y'})=0$, i.e.
 $r_{Y'}=r_X+\varepsilon_2/2\;(\mathrm{mod}\;2\Zb)$,
 $\varepsilon_2=\pm1$. Now $Y\neq Y'$ implies that
 $\varepsilon_2=-\varepsilon_1$. Setting e.g.
 $r_Y=r_X+1/2$,
 $r_{Y'}=r_X-1/2$ in (\ref{oddboundaux6}) one
 finds
 \ben
 \cos\pi (r_X-1/2)+\cos\pi(2r_X+1/2)=\cos\pi r_{Z'}+\cos\pi(2r_X-1/2).
 \ebn
 Since we are looking for solutions of type
 ``$\text{II}_{\varphi}+\text{II}_{\psi}$'' of this equation and
 since $Y'\neq Z'$, the only possible $N$ is equal to 3.\\
 (1.3.2) Next consider solutions of type
 ``$\text{III}_{\varphi}+\text{II}_{\psi}$''.
 It can be assumed that $ \cos\pi(r_X\pm r_{Y'})$  do
 not belong simultaneously to
 ($\text{II}_{\psi}$), as this would lead to $X=0$ ($N=2$) or
 $Y'=0$ (case studied above).

 We may further assume that they are
 not simultaneously in ($\text{III}_{\varphi}$), because one would
 then have $N=3$ or $Y'=\varepsilon_2$, where
 $\varepsilon_2=\pm1$. In the latter case (\ref{oddboundaux6}) would transform
 into
 \ben
 \varepsilon_2\cos\pi/3+\cos\pi(2r_X+\varepsilon_1/2)=\cos\pi r_{Z'}+
 \varepsilon_2\cos\pi r_{X}.
 \ebn
 Since solutions of this equation should have type ``$\text{II}_{\varphi}+\text{II}_{\psi}$''
 and since $Y'\neq Z'$, one concludes that $N=3$. \\
 (1.3.2.1) Let $\cos\pi r_{Y'}$ be in ($\text{II}_{\psi}$), then
 we may write (\ref{oddboundaux6}) as
 \ben
 \begin{cases}
 \cos\pi(2r_X+\varepsilon_1/2)=\cos\pi r_{Z'}+
 \cos\pi(r_X+r_{Y'}), & \qquad (\text{III}_{\varphi}) \\
 \cos\pi r_{Y'}=\cos\pi (r_X-r_{Y'}).& \qquad (\text{II}_{\psi})
 \end{cases}
 \ebn
 Second equation implies that
 $r_{X}=2r_{Y'}\;(\mathrm{mod}\;2\Zb)$.
 Substituting this into the first equation one finds $ \cos\pi(4r_{Y'}+\varepsilon_1/2)=\cos\pi r_{Z'}+
 \cos3\pi r_{Y'}$, therefore $N$ can only be equal to $3,7,21$.\\
 (1.3.2.2) Let $\cos\pi r_{Y'}$ be in ($\text{III}_{\varphi}$) and
 let $\cos\pi(2r_X+\varepsilon_1/2)$ be in ($\text{II}_{\psi}$).
 Then one can write
  \ben
 \begin{cases}
 \cos\pi r_{Y'}=\cos\pi r_{Z'}+
 \cos\pi(r_X-r_{Y'}), & \qquad (\text{III}_{\varphi}) \\
 \cos\pi(2r_X+\varepsilon_1/2)=\cos\pi (r_X+r_{Y'}),& \qquad (\text{II}_{\psi})
 \end{cases}
 \ebn
 and it follows that possible values of $N$ are $3,7,21$.
 Similarly if both $\cos\pi r_{Y'}$ and $\cos\pi(2r_X+\varepsilon_1/2)$
 are in ($\text{III}_{\varphi}$), one finds $N=3,5,9,15$.\\
 (1.4) Finally suppose that $\cos\pi(r_X-r_{Y'})=0$.
 Then $r_{Y'}=r_X+\varepsilon_1/2\;(\mathrm{mod}\;2\Zb)$, $\varepsilon_1=\pm1$
 and from (\ref{oddboundaux1}) follows the relation
 \be\label{oddboundaux7}
 \cos\pi (r_X+\varepsilon_1/2)+\cos\pi(r_X-r_Y)+\cos\pi(r_X+r_Y)=
 \cos\pi r_{Z'}+ \cos\pi(2r_X+\varepsilon_1/2).
 \eb
 It is not necessary to examine solutions of (\ref{oddboundaux7}) of type
 ``$\text{II}_{\varphi}+\text{II}_{\psi}+\text{I}$'' because all
 cases when $Y'=0$, $Z'=0$ or $\cos\pi(r_X\pm r_Y)=0$ have already
 been considered above, and $\cos\pi(2r_X\pm1/2)=0$ gives $N=2$.
 Hence we may restrict our attention to solutions of type
 ``$\text{III}_{\varphi}+\text{II}_{\psi}$''.
 \begin{itemize}
 \item $\cos\pi (r_X+\varepsilon_1/2)$ and
 $\cos\pi(2r_X+\varepsilon_1/2)$ cannot be simultaneously in ($\text{II}_{\psi}$)
 unless $N=3$ and in ($\text{III}_{\varphi}$) unless $N=3,9$.
 Therefore one can assume that they are divided between ($\text{III}_{\varphi}$)
 and ($\text{II}_{\psi}$).
 \item If both $\cos\pi(r_X\pm r_{Y})$ belong to
 ($\text{III}_{\varphi}$), then either $N=3$ or $Y=\varepsilon_2$,
 $\varepsilon_2=\pm1$, but in the latter case (\ref{oddboundaux7})
 becomes
 \ben
 \cos\pi (r_X+\varepsilon_1/2)+\varepsilon_2\cos\pi r_X=
 \cos\pi r_{Z'}+ \cos\pi(2r_X+\varepsilon_1/2).
 \ebn
 Solution of this equation should be of type
 ``$\text{II}_{\varphi}+\text{II}_{\psi}$''.
 Since $Y'\neq Z'$ and by the above assumption $\cos\pi (r_X+\varepsilon_1/2)$ and
 $\cos\pi(2r_X+\varepsilon_1/2)$ are not in the same pair,
 this can happen only
 if $\cos\pi (r_X+\varepsilon_1/2)+\varepsilon_2\cos\pi r_X=0$,
 i.e. odd $N$ are impossible. Thus we can assume that $\cos\pi(r_X\pm
 r_{Y})$ in (\ref{oddboundaux7}) are also divided between ($\text{III}_{\varphi}$)
 and ($\text{II}_{\psi}$) and in particular $\cos\pi r_{Z'}$
 belongs to ($\text{III}_{\varphi}$).
 \end{itemize}
 We then have two inequivalent possibilities:
 \ben
 \begin{cases}
 \cos\pi(r_X+r_Y)=\cos\pi r_{Z'}+
 \cos\pi(2r_X+\varepsilon_1/2), & \qquad (\text{III}_{\varphi}) \\
 \cos\pi(r_X+\varepsilon_1/2)+\cos\pi (r_X-r_{Y})=0.& \qquad (\text{II}_{\psi})
 \end{cases}\tag{1.4.1}
 \ebn
 From the second equation one finds that either $Y=0$ or
 $r_Y=2r_X+\varepsilon_1/2+1\;(\mathrm{mod}\;2\Zb)$.
 In the former case, substitution into the first equation
 gives admissible values $N=3,9$, while for the latter $N=3,5,15$.
 \ben
 \begin{cases}
 \cos\pi(r_X+\varepsilon_1/2)+\cos\pi(r_X+r_Y)=\cos\pi r_{Z'}, & \qquad (\text{III}_{\varphi}) \\
 \cos\pi (r_X-r_{Y})=\cos\pi(2r_X+\varepsilon_1/2).& \qquad (\text{II}_{\psi})
 \end{cases}\tag{1.4.2}
 \ebn
 Here from $(\text{II}_{\psi})$ follows that either
 $r_Y=-r_X-\varepsilon_1/2\;(\mathrm{mod}\;2\Zb)$ (forbidden because then
 $Y=Y'$) or $r_Y=3r_X+\varepsilon_1/2\;(\mathrm{mod}\;2\Zb)$. In
 the latter case first equation implies that $N=3,5,9,15$.\\
 \underline{Case (2)}. Now we come back to the initial equation
  (\ref{oddboundaux1}) and consider its solutions of type
 ``$\text{II}_{\varphi}+\text{II}_{\psi}+\text{II}_{\mu}$''.

 It can be assumed that $\cos\pi(r_X\pm r_{Y'})$ are not in the
 same pair, as otherwise $X=0$ ($N=2$) or $Y'=0$ (already
 considered). Similarly, if both $\cos\pi(r_X\pm r_{Y})$ are in
 the same pair, then $Y=0$ and one can write
 \ben
 \begin{cases}
 \cos\pi r_{Y'}=\cos\pi(r_X-r_{Y'}),& \qquad (\text{II}_{\varphi}) \\
 \cos\pi r_{Z'}+\cos\pi(r_X+r_{Y'})=0.& \qquad (\text{II}_{\psi})
 \end{cases}
 \ebn
 Since $X\neq\pm2$, from ($\text{II}_{\varphi}$) follows that
 $r_X=2r_{Y'}\;(\mathrm{mod}\;2\Zb)$ and then $Z'=-2\cos3\pi r_{Y'}$.
 Moreover $Y=0$ implies that $\omega=Y'$, therefore $Y''=-XZ'$,
 i.e.
 \ben
 \cos\pi r_{Y''}=\cos\pi r_{Y'}+\cos5\pi r_{Y'}.
 \ebn
 For $N\geq 9$ we can apply Lemma~\ref{lemcos} to the last
 relation. Its solutions of type (\ref{iii1}) and (\ref{iiia})
 lead to $N=3,5$ and $N=3,9$ correspondingly. Since $Y',Y''\neq0$
 (because we already have $Y=0$), solutions of type ``$\text{II}_{\varphi}+\text{I}$''
 are possible only if $N=5$.

 Hence we can assume that $\cos\pi(r_X\pm r_{Y})$ are divided
 between two different pairs. These cannot be the same as
 for $\cos\pi(r_X\pm r_{Y'})$, otherwise the third pair would
 give $Y'=Z'$. Therefore we may assume
 one of the pairs in (\ref{oddboundaux1}) to be
 \ben
 \cos\pi(r_X-r_Y)=\cos\pi(r_X-r_{Y'}).\qquad{(\text{II}_{\varphi})}
 \ebn
 Since $Y\neq Y'$, the last relation gives
 $r_Y=2r_X-r_{Y'}\;(\mathrm{mod}\;2\Zb)$. Now for the remaining
 two pairs there are two inequivalent possibilities:\\
 (2.1) If $\cos\pi r_{Y'}$ and $\cos\pi(r_X+r_{Y})$ are in the
 same pair, then
 \ben
 \begin{cases}
 \cos\pi r_{Y'}+\cos\pi(3r_X-r_{Y'})=0,&\qquad (\text{II}_{\psi})\\
 \cos\pi r_{Z'}+\cos\pi(r_X+r_{Y'})=0.&\qquad (\text{II}_{\mu})
 \end{cases}
 \ebn
 From ($\text{II}_{\psi}$) one finds that either $N=3$ or
 $\ds\cos{\pi(3r_X-2r_{Y'})}/{2}=0$. In the latter case, compute
 $\omega$:
 \ben
 \ds\omega=Y+Y'+XY=4\,\cos{\pi
 r_X}/{2}\,\cos{\pi(3r_X-2r_{Y'})}/{2}=0,
 \ebn
 i.e. the initial assumption $\omega\neq0$ does not hold.\\
 (2.2)  If $\cos\pi r_{Y'}$ and $\cos\pi(r_X+r_{Y'})$ are in the
 same pair, then
  \ben
 \begin{cases}
 \cos\pi r_{Y'}=\cos\pi(r_X+r_{Y'}),&\qquad (\text{II}_{\psi})\\
 \cos\pi r_{Z'}=\cos\pi(3r_X-r_{Y'}).&\qquad (\text{II}_{\mu})
 \end{cases}
 \ebn
 First equation implies that $r_X=-2r_{Y'}\;(\mathrm{mod}\;2\Zb)$.
 Therefore $X=2\cos2\pi r_{Y'}$, $Y=2\cos5\pi r_{Y'}$,
 $Z'=2\cos7\pi r_{Y'}$. Let us compute $\omega=Y+Y'+XY$:
 \ben
 \ds\omega=2\cos\pi r_{Y'}+2\cos3\pi r_{Y'}+2\cos5\pi r_{Y'}+
 2\cos7\pi r_{Y'}.
 \ebn
 The computation of $Y''=\omega-Y'-XZ'$ now gives
 \be\label{oddboundaux8}
 \cos\pi r_{Y''}=\cos3\pi r_{Y'}+\cos7\pi r_{Y'}-\cos9\pi r_{Y'}.
 \eb
 For $N\geq9$, we can apply to (\ref{oddboundaux8})
 Lemma~\ref{lemcos}. Solutions of type (\ref{iv}),
 ``$\text{III}_{\text{1}}+\text{I}$'' and
 ``$\text{III}_{\varphi}+\text{I}$'' can lead only to
 $N=3,5,7,9,15,21$. Since $Y''\neq Z'$,
 solutions of type ``$\text{II}_{\varphi}+\text{II}_{\psi}$''
 are possible only if $N=5$.\\
 \underline{Case (3)}. It remains to consider solutions of
 (\ref{oddboundaux1}) of type
 ``$\text{III}_{\varphi}+\text{III}_{\psi}$''.\\
 (3.1) If both $\cos\pi(r_X\pm r_{Y'})$ appear in
 the same triple, then $N=3$ or $Y'=\pm1$. In
 the latter case, (\ref{oddboundaux1}) transforms into
 \be\label{oddboundaux9}
 \pm \cos\pi/3+\cos\pi(r_X+r_Y)+\cos\pi(r_X-r_Y)=
 \cos\pi r_{Z'}\pm\cos\pi r_X.
 \eb
 The solution of (\ref{oddboundaux9}) should have type
 ``$\text{III}_{\varphi}+\text{II}_{\psi}$'', and moreover
 $\cos\pi r_X$ belongs to ($\text{II}_{\psi}$). If the second
 cosine in ($\text{II}_{\psi}$) is $\cos\pi/3$, then $N=3$. If
 $\cos\pi r_{Z'}\pm\cos\pi r_X=0$, then from (\ref{iiia}) again
 follows $N=3$. Therefore it can be assumed that
 \ben
 \begin{cases}
 \pm\cos\pi /3+\cos\pi(r_X+r_Y)=\cos\pi r_{Z'},&\qquad
 (\text{III}_{\varphi})\\
 \cos\pi(r_X-r_{Y})=\pm\cos\pi r_X.&\qquad (\text{II}_{\psi})
 \end{cases}
 \ebn
 Since $Y\neq \pm2$, second equation implies that
 $r_Y=2r_X+1/2\mp1/2$, but then from the first equation follows $N=3,9$.
 Hence from now on we assume that $\cos\pi(r_X\pm r_{Y'})$ belong to
 different triples.\\
 (3.2) If $\cos\pi(r_X\pm r_{Y})$ are in the same triple, then
 $N=3$ or $Y=\pm1$. In the latter case
 (\ref{oddboundaux1}) can be rewritten as
  \ben
 \begin{cases}
 \cos\pi r_{Y'}=\cos\pi r_{Z'}+\cos\pi(r_X+r_{Y'}),&\qquad
 (\text{III}_{\varphi})\\
 \pm\cos\pi r_X=\cos\pi (r_X-r_{Y'}).&\qquad (\text{II}_{\psi})
 \end{cases}
 \ebn
 Again from ($\text{II}_{\psi}$) follows  $r_{Y'}=2r_X+1/2\mp1/2$,
 and (\ref{iiia}) then implies that $N=3,5,15$. Therefore
 we assume in the following that $\cos\pi(r_X\pm r_{Y})$,
 as well as $\cos\pi r_{Y'}$ and $\cos\pi r_{Z'}$,
 are divided between the two triples.\\
 (3.3) Without loss of generality we can now write
 (\ref{oddboundaux1}) as
 \be\label{oddboundaux11}
 \begin{cases}
 \cos\pi r_{Y'}+\cos\pi(r_X-r_Y)-\cos\pi(r_X-r_{Y'})=0,&\qquad
 (\text{III}_{\varphi})\\
 \cos\pi r_{Z'}+\cos\pi (r_X+r_{Y'})-\cos\pi(r_X+r_Y)=0,&\qquad (\text{III}_{\psi})
 \end{cases}
 \eb
 or, in another form,
    \ben
 \begin{cases}
 \ds\cos\pi r_{Y'}+2\sin\frac{\pi(2r_X-r_Y-r_{Y'})}{2}
 \sin\frac{\pi(r_Y-r_{Y'})}{2}=0,&\qquad
 (\text{III}_{\varphi})\\
 \ds\cos\pi r_{Z'}+2\sin\frac{\pi(2r_X+r_Y+r_{Y'})}{2}
 \sin\frac{\pi(r_Y-r_{Y'})}{2}=0.&\qquad (\text{III}_{\psi})
 \end{cases}
 \ebn
 If $\ds\sin\frac{\pi(r_Y-r_{Y'})}{2}\neq\pm\frac{1}{2}$, then one
 should simultaneously have
 \be\label{oddboundaux10}
 \sin\frac{\pi(2r_X-r_Y-r_{Y'})}{2}=\frac{\varepsilon_1}{2},\qquad
 \sin\frac{\pi(2r_X+r_Y+r_{Y'})}{2}=\frac{\varepsilon_2}{2},
 \eb
 where $\varepsilon_{1,2}=\pm 1$ (in fact $\varepsilon_2=-\varepsilon_1$,
 otherwise $Y'=Z'$). Equations (\ref{oddboundaux10}) lead to
  $N=3$, therefore we can  assume that
 \be\label{oddboundaux12}
 \ds\sin\frac{\pi(r_Y-r_{Y'})}{2}=\frac{\varepsilon_3}{2},\qquad
 \varepsilon_3=\pm1.
 \eb
 Let us compute $Y''=Y+XY-XZ'$ using (\ref{oddboundaux11}) and
 (\ref{oddboundaux12}). After some
 simplifications one finds
 \be\label{oddboundaux13}
 \ds\cos\pi r_{Y''}=\cos\pi(r_X+r_Y)+\cos\pi (r_X-r_{Y'})+
 \varepsilon_3\sin\frac{\pi(4r_X+r_Y+r_{Y'})}{2}\,.
 \eb
 Relation (\ref{oddboundaux12}) implies that
 $r_Y=r_{Y'}+\varepsilon_3/3\;(\mathrm{mod}\;4\Zb)$ or
 $r_Y=r_{Y'}+5\varepsilon_3/3\;(\mathrm{mod}\;4\Zb)$.
 Similarly, first relation in (\ref{oddboundaux11}) gives either $N=3$
 or $r_X=2r_{Y'}+\varepsilon_4/3\;(\mathrm{mod}\;2\Zb)$,
 $\varepsilon_4=\pm1$. We now substitute this into
 (\ref{oddboundaux13}) and apply Lemma~\ref{lemcos} (for $N\geq9$).
 Solutions of type (\ref{iv}) and ``$\text{III}_{\text{1}}+\text{I}$''
 then lead to admissible values $N=3,5,7,9,15,21$ and $N=3,5,15$ correspondingly,
 while solutions of type ``$\text{III}_{\varphi}+\text{I}$'' and
 ``$\text{II}_{\varphi}+\text{II}_{\psi}$''  give $N=3,5,9,15$ and
 $N=3,9$. This concludes the proof of Proposition~\ref{oddbound2}.
   \epf
 \begin{lem}\label{oddbound3}
 Let $N$ and $n_X$ be odd and let $\omy=\omz\neq0$. If the graph $\Sigma(O_{yz})$ is
 a simple cycle and $O_{yz}$ contains a point with coordinate $Z$ (or $Y$)
 equal to $0$,
 then the only possible values of $N$ are $3,5,7,9,15,21$.
 \end{lem}
 \pf Analogously to the previous proof, let us label
 the vertices of $\Sigma(O_{yz})$ by their coordinates $(Y,Z)$,
 as shown in Fig.~4. Because of
 simple cycle assumption all points of $O_{yz}$ are good,
 therefore all $\{Y_k\}$ and $\{Z_k\}$
 have the form (\ref{yzgood}). It will be assumed that $N>3$, then by
 Lemma~\ref{distinct}
 four numbers $Y,Y',Y'',Z'$ are distinct and
 non-zero (recall that $Y_k=Z_{k+(N-1)/2}$).
  \begin{figure}[!h]
 \begin{center}
 \resizebox{7.5cm}{!}{
 \includegraphics{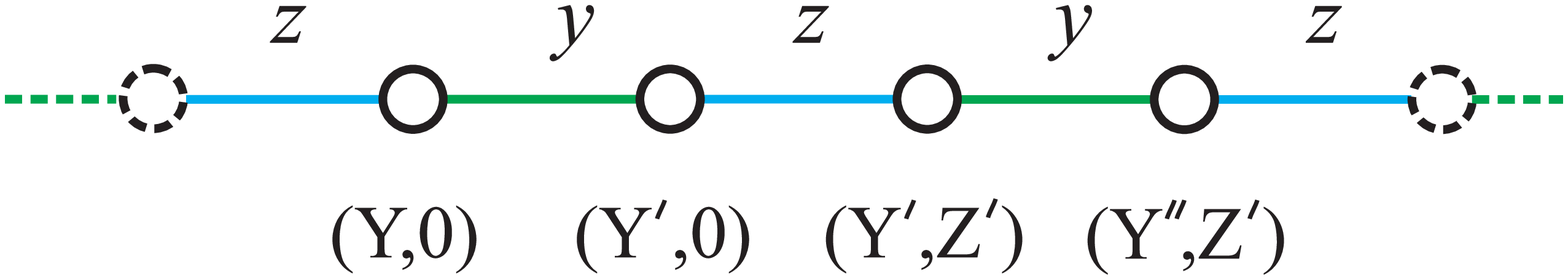}} \\   Fig. 4
 \end{center}
 \end{figure}

 We now apply Lemma~\ref{lemcos} to the relation
 \be\label{oddboundaux14}
 \cos\pi r_Y+\cos\pi r_{Y'}=\cos\pi r_{Z'}+\cos\pi(r_X+r_{Y'})
 +\cos\pi(r_X-r_{Y'}).
 \eb

 Its solutions of type (\ref{va}), (\ref{v1})--(\ref{v3}),
 ``$\text{IV}+\text{I}$'' can lead only to $N=3,5,7,15,21$.

 The solutions of type
 ``$\text{II}_{\varphi}+\text{II}_{\psi}+\text{I}$'' are
 forbidden. Indeed, since $Y,Y',Z'\neq0$,
 in this case one could write
 $\cos\pi(r_X-r_{Y'})=0$, but then one of the pairs
 ($\text{II}_{\varphi}$), ($\text{II}_{\psi}$) would give
 $Y=Z'$ or $Y'=Z'$ (impossible) or $Y+Y'=0$ (excluded because
 then $\omega=0$).

 Next we consider solutions of type
 ``$\text{III}_{\text{1}}+\text{II}_{\varphi}$''.
 Since $Y'\neq0$, two cosines
 $\cos\pi(r_X\pm r_{Y'})$ cannot belong both to
 (\ref{iia}). They can neither be simultaneously
 in (\ref{iii1}), as (\ref{iia}) would then give $Y=Z'$ or $Y'=Z'$
 or $Y+Y'=0$. Therefore it can be assumed that $\cos\pi(r_X- r_{Y'})$
 belongs to (\ref{iia}) and $\cos\pi(r_X+ r_{Y'})$ is in (\ref{iii1}).
 Now if $\cos\pi r_{Y'}$ is in (\ref{iii1}), then admissible values of $N$
 are $3,5,15$. If
 $\cos\pi r_{Y'}$ belongs to (\ref{iia}), then
 $r_X=2r_{Y'}\;(\mathrm{mod}\;2\Zb)$. Substituting this into
 (\ref{iii1}), we obtain $N=5,9,15$.

 It remains to consider solutions of (\ref{oddboundaux14}) of type
 ``$\text{III}_{\varphi}+\text{II}_{\psi}$''. By the same argument
 as above we can assume that $\cos\pi(r_X- r_{Y'})$
 is in ($\text{II}_{\psi}$) and $\cos\pi(r_X+ r_{Y'})$ is in
 (\ref{iiia}).

  Assume that $\cos\pi r_{Y'}$ is in
 ($\text{II}_{\psi}$). Then $r_X=2r_{Y'}\;(\mathrm{mod}\;2\Zb)$
 and the triple (\ref{iiia}) becomes
 \ben
 \cos\pi r_Y=\cos\pi r_{Z'}+\cos3\pi r_{Y'}.
 \ebn
 Therefore we can assume that
 $r_Y=3r_{Y'}\pm1/3\;(\mathrm{mod}\;2\Zb)$,
 $r_{Z'}=3r_{Y'}\pm2/3\;(\mathrm{mod}\;2\Zb)$.
 Let us substitute these expressions into an easily verified
 relation
 \be\label{oddboundaux15}
 \cos\pi r_{Y''}= \cos\pi
 r_Y-\cos\pi(r_X-r_{Z'})-\cos\pi(r_X+r_{Z'}).
 \eb
 Its solutions of type (\ref{iv}) and ``$\text{III}_{\text{1}}+\text{I}$''
 lead to admissible values $N=3,5,7,15,21$ and $N=3,5,15$
 correspondingly (in fact this conclusion does not depend on any of
 our previous assumptions). Since $Y,Y''\neq 0$,
 solutions of type ``$\text{III}_{\varphi}+\text{I}$''
 give $N=3,5,15$. Finally, since $Y\neq Y''$, solutions of
 type ``$\text{II}_{\varphi}+\text{II}_{\psi}$'' are only possible for $N=3$.

 On the other hand, if $\cos\pi r_{Y'}$ belongs to the triple (\ref{iiia})
 then, since $\cos\pi(r_X+r_{Y'})$ is also in (\ref{iiia}), we can
 set $r_X=-2r_{Y'}+\varepsilon/3\;(\mathrm{mod}\;2\Zb)$,
 $\varepsilon=\pm1$,
 otherwise $N=3$. Hence \\
 (1) If  $\cos\pi r_{Y}$ is the third cosine in (\ref{iiia}), then
 \ben
 Y=-2\cos\pi(r_{Y'}+\varepsilon/3),\qquad
 Z'=-2\cos\pi(3r_{Y'}-\varepsilon/3).
 \ebn
 Now let us look at the equation
 (\ref{oddboundaux15}). When its solution has type
 ``$\text{III}_{\varphi}+\text{I}$'',
 it can be assumed that $\cos\pi(r_X-r_{Z'})=0$, but then $N=3,5,15$.
 For solutions of type ``$\text{II}_{\varphi}+\text{II}_{\psi}$''
 we can write $\cos\pi r_{Y}=\cos\pi(r_X-r_{Z'})$, which leads to
 $N=3,9$.\\
 (2) If  $\cos\pi r_{Y}$ belongs to ($\text{II}_{\psi}$), then one
 finds
  \ben
 Y=2\cos\pi(3r_{Y'}-\varepsilon/3),\qquad
 Z'=2\cos\pi(r_{Y'}+\varepsilon/3).
 \ebn
 In this case, solutions of (\ref{oddboundaux15})
 of type ``$\text{III}_{\varphi}+\text{I}$'' and
 ``$\text{II}_{\varphi}+\text{II}_{\psi}$''
 lead to admissible values $N=3,9$.
 \epf
 \begin{prop}\label{auuu}
 Let $N$ and $n_X$ be odd and let $\omy=\omz\neq0$. If the
 graph $\Sigma(O_{yz})$ is a simple cycle
 then the only possible values of $N$ are $3,5,7,9,11,15,21$.
 \end{prop}
 \pf Let us start with the obvious relation $Y+XZ=Y''+XZ'$ (see Fig.~5), written as
 \begin{align}
 \label{oddboundaux16}
 &\cos\pi r_{Y}+\cos\pi(r_X+r_{Z})+\cos\pi(r_X-r_{Z})=\qquad \\
 \nonumber&\qquad=\cos\pi r_{Y''}+\cos\pi(r_X+r_{Z'})+\cos\pi(r_X-r_{Z'}).
 \end{align}
 We can assume that this relation does not contain zero cosines. Indeed, the case when
 $Y=0$ or $Y''=0$ is completely described by Lemma~\ref{oddbound3}. If $\cos\pi(r_X\pm r_{Z})=0$
 or $\cos\pi(r_X\pm r_{Z'})=0$, then $Z$ or $Z'$ is equal to $\pm\sqrt{4-X^2}$.
 Now recall that by Lemma~\ref{distinct} in a simple cycle all $\{Z_k\}$ are distinct, therefore
 already for $N\geq5$ it will be possible to find a pair $(Z_k,Z_{k+1})$ which does not contain prescribed
 two values $\pm\sqrt{4-X^2}$ (Assumption~1).
   \begin{figure}[!h]
 \begin{center}
 \resizebox{10cm}{!}{
 \includegraphics{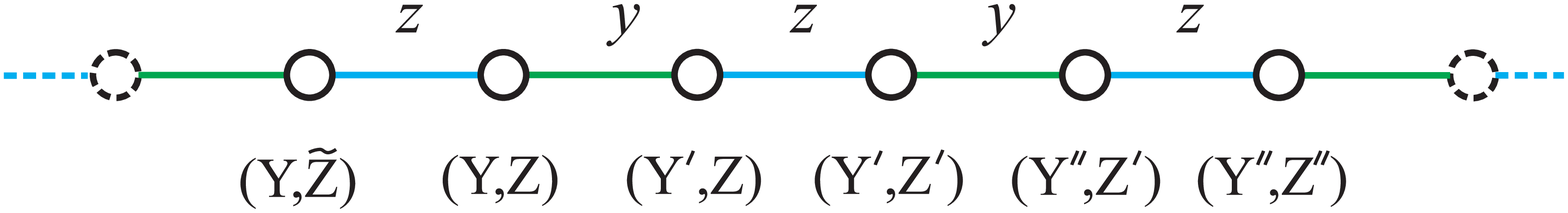}} \\   Fig. 5\\
 \end{center}
 \end{figure}

 Next we exclude solutions of type (\ref{vi1})--(\ref{vi5}), ``$\text{IV}+\text{II}_{\varphi}$'',
 ``$\text{III}_{\text{1}}+\text{III}_{\text{1}}$'',  ``$\text{III}_{\text{1}}+\text{III}_{\varphi}$'',
 as they can lead only to $N=3,5,7,9,11,15,21$ (note that solutions of (\ref{oddboundaux16})
 satisfy a condition similar to (a) in the proof of Proposition~\ref{oddbound2}). Then there remain
 three types of possible solution 6-tuples:

 (1) ``$\text{II}_{\varphi}+\text{II}_{\psi}+\text{II}_{\mu}$'';

 (2) ``$\text{III}_{\varphi}+\text{III}_{\psi}$'';

 (3) ``$\text{VI}_{\varphi}$''.\\
 \underline{Case (1)}. It can be assumed that
 two cosines $\cos\pi(r_X\pm r_Z)$ (and $\cos\pi(r_X\pm r_{Z'})$)
 are divided between two different pairs. Otherwise $Z=0$ (resp. $Z'=0$) and
 one obtains restrictions on $N$ from Lemma~\ref{oddbound3}. The pairs cannot be the same
 in both cases because then  $Y=Y''$.  Therefore we can set one of the pairs to be
 \be
 \cos\pi(r_X-r_Z)=\cos\pi(r_X-r_{Z'}).\qquad (\text{II}_{\varphi})
 \eb
 Since $Z\neq Z'$, one has $r_{Z'}=2r_X-r_Z\;(\mathrm{mod}\;2\Zb)$. For the remaining two
 pairs, there are two inequivalent possibilities:
 \ben
 \begin{cases}
 \cos\pi r_Y+\cos\pi(r_X+r_Z)=0,\qquad & (\text{II}_{\psi})\\
 \cos\pi r_{Y''}+\cos\pi(3r_X-r_{Z})=0.\qquad & (\text{II}_{\mu})
 \end{cases}\tag{1.1}
 \ebn
 Here from $Y+Y'+XZ=Z+Z'+XY'$ follows that either $N=3$ or
 $Y'=-2\cos\pi(r_X-r_Z)$. In the latter case, however, computing
 $\omega=Y+Y'+XZ$ we find forbidden value $\omega=0$.
  \ben
 \begin{cases}
 \cos\pi r_Y=\cos\pi(3r_X-r_{Z}),\qquad & (\text{II}_{\psi})\\
 \cos\pi r_{Y''}=\cos\pi(r_X+r_{Z}).\qquad & (\text{II}_{\mu})
 \end{cases}\tag{1.2}
 \ebn
 Substituting these relations into  $\tilde{Z}+XY=Z'+XY'$ and $Z''+XY''=Z+XY'$,
 one obtains
 \begin{eqnarray}
 \label{oddboundaux17}
 \cos\pi r_{\tilde{Z}}+\cos\pi(4r_X-r_Z)&=&\cos\pi (r_X-r_{Y'})+\cos\pi (r_X+r_{Y'}),\\
 \label{oddboundaux18}
 \cos\pi r_{Z''}+\cos\pi(2r_X+r_Z)&=&\cos\pi (r_X-r_{Y'})+\cos\pi (r_X+r_{Y'}).
 \end{eqnarray}
 Solutions of (\ref{oddboundaux17}), (\ref{oddboundaux18}) of type (\ref{iv}) and
 ``$\text{III}_1+\text{I}$'' can lead only to $N=3,5,7,15,21$, therefore we can restrict our
 attention to solutions of type ``$\text{III}_{\varphi}+\text{I}$'' and
 ``$\text{II}_{\varphi}+\text{II}_{\psi}$''.\\
 (1.2.1) Suppose that the solution of (\ref{oddboundaux17}) is of type
 ``$\text{III}_{\varphi}+\text{I}$''. If $\cos\pi(4r_X-r_Z)=0$ and $\cos\pi (r_X\pm r_{Y'})$
 are in (\ref{iiia}), then $N=3$ or $r_Z=4r_X+\varepsilon_1/2\;(\mathrm{mod}\;2\Zb)$, $Y'=\varepsilon_2$,
 $\varepsilon_{1,2}=\pm1$. In the second case (\ref{oddboundaux18}) transforms into
 \ben
 \cos\pi r_{Z''}+\cos\pi(6r_X+\varepsilon_1/2)=\varepsilon_2 \cos\pi r_X.
 \ebn
 Now if the solution of this equation has type (\ref{iiia}) or (\ref{iii1}), then $N=3,5,7,15,21$.
 Since it can be assumed that $Z''\neq0$, type ``$\text{II}_{\varphi}+\text{I}$'' solutions give $N=3$.

 On the other hand, if $\cos\pi(r_X-r_{Y'})=0$, i.e. $r_{Y'}=r_X+\varepsilon_1/2\;(\mathrm{mod}\;2\Zb)$,
 $\varepsilon_1=\pm1$, then the triple (\ref{iiia}) in (\ref{oddboundaux17}) is given by
 \ben
 \cos\pi r_{\tilde{Z}}+\cos\pi(4r_X-r_Z)=\cos\pi (2r_X+\varepsilon_1/2).
 \ebn
 This relation implies that either (a) $r_Z=2r_X-\varepsilon_1/2+\varepsilon_2/3\;(\mathrm{mod}\;2\Zb)$
 or (b) $r_Z=6r_X+\varepsilon_1/2+\varepsilon_2/3\;(\mathrm{mod}\;2\Zb)$. In the case (a)
 equation (\ref{oddboundaux18}) transforms into
 \ben
  \cos\pi r_{Z''}+\cos\pi(4r_X-\varepsilon_1/2+\varepsilon_2/3)=\cos\pi (2r_X+\varepsilon_1/2).
 \ebn
 Its solutions of type (\ref{iiia}) and ``$\text{II}_{\varphi}+\text{I}$'' lead to admissible values
 $N=3,9$. Similarly, in the case (b) relation (\ref{oddboundaux18}) gives $N=3,5,9,15$.\\
 (1.2.2) The case when the solution of (\ref{oddboundaux18}) is of type
 ``$\text{III}_{\varphi}+\text{I}$'' is treated analogously to (1.2.1), hence we can assume
 that solutions of both (\ref{oddboundaux17}) and (\ref{oddboundaux18}) have the form
 ``$\text{II}_{\varphi}+\text{II}_{\psi}$''. Thanks to Lemma~\ref{oddbound3}, it can be assumed
 that $Y'\neq0$ so that $\cos\pi(r_X\pm r_{Y'})$ in (\ref{oddboundaux17}), (\ref{oddboundaux18})
 are divided between the two pairs.
 Since $\tilde{Z}\neq Z''$, we may write without loss of generality
 \ben
 \begin{cases}
 \cos\pi (4r_X-r_Z)=\cos\pi(r_X-r_{Y'}),\\
 \cos\pi (2r_X+r_Z)=\cos\pi(r_X+r_{Y'}).
 \end{cases}
 \ebn
 From the first equation follows either $r_{Y'}=-3r_X+r_Z\;(\mathrm{mod}\;2\Zb)$ (forbidden because then
 $Y=Y'$) or $r_{Y'}=5r_X-r_Z\;(\mathrm{mod}\;2\Zb)$. In the latter case the second equation becomes
 \ben
 \cos\pi (2r_X+r_Z)=\cos\pi(6r_X-r_{Z}),
 \ebn
 and implies that $2r_X-r_Z\in\Zb$. This in turn gives $Z'=\pm2$, which is impossible as all points
 in $O_{yz}$ are good.\\
 \underline{Case (2)}. Suppose that $Z$ and $Z'$ are not
 equal to $\pm1$ (Assumption~2). Clearly for $N\geq9$ one will always be able
 to find in $O_{yz}$ a pair $(Z,Z')$ satisfying Assumptions~1 and~2. Then
 in (\ref{oddboundaux16}) the two cosines $\cos\pi(r_X\pm r_Z)$,
 as well as $\cos\pi(r_X\pm r_{Z'})$, are divided between the two triples (\ref{iiia}) and
 ($\text{III}_{\psi}$),
 otherwise $X=\pm1$ and $N=3$. We can therefore write
 \be\label{oddboundaux19}
 \begin{cases}
 \cos\pi r_Y+\cos\pi(r_X-r_Z)-\cos\pi(r_X-r_{Z'})=0,\qquad & (\text{III}_{\varphi})\\
 \cos\pi r_{Y''}-\cos\pi(r_X+r_{Z})+\cos\pi(r_X+r_{Z'})=0.\qquad & (\text{III}_{\psi})
 \end{cases}
 \eb
 Similarly to the proof of Proposition~\ref{oddbound2}, case (3.3) one can show that
 \ben\sin\frac{\pi(r_Z-r_{Z'})}{2}=\pm\frac{1}{2},\ebn
 i.e.  $r_{Z'}=r_Z+\varepsilon_1/3\;(\mathrm{mod}\;2\Zb)$, $\varepsilon_1=\pm1$.

 From $\omega=Y+Y'+XZ=Z+Z'+XY'$ follows that
 \ben
 (X-1)\,\omega=XY+(X^2-2)Z+Z-Z'.
 \ebn
 Substituting (\ref{oddboundaux19}) into this relation, we find
 \begin{align}
 \nonumber(X-1)\,\omega=2\cos\pi(2r_X+r_Z)+2\cos\pi(2r_X-r_{Z'})=\\
 \nonumber=2\cos\pi(2r_X+r_Z)+2\cos\pi(2r_X-r_{Z}-\varepsilon_1/3).
 \end{align}

 Recall that for a simple cycle of length $N$, one may write $N$
 relations of the form (\ref{oddboundaux16}) which correspond to
 different unordered pairs $(Z,Z')$. Suppose there exists a
 second relation whose solution has the form
 ``$\text{III}_{\varphi}+\text{III}_{\psi}$'', and the associated
 pair $(\bar{Z},\bar{Z'})$ satisfies Assumptions~1 and~2. Then we
 can write
 \be\label{oddboundaux20}
 \cos\pi(2r_X+r_Z)+\cos\pi(2r_X-r_{Z}-\varepsilon_1/3)=\cos\pi(2r_X+r_{\bar{Z}})+\cos\pi(2r_X-r_{\bar{Z}}-\varepsilon_2/3),
 \eb
 where $r_{\bar{Z'}}=r_{\bar{Z}}+\varepsilon_2/3\;(\mathrm{mod}\;2\Zb)$,
 $\varepsilon_2=\pm1$. If
 $\varepsilon_1=\varepsilon_2$, then (\ref{oddboundaux20}) implies that
  either $N=3$ or the pairs $(Z,Z')$
 and $(\bar{Z},\bar{Z'})$ coincide. Let us now set
 $\varepsilon_2=-\varepsilon_1$ and consider rational solutions of
 (\ref{oddboundaux20}).

 Solutions of type (\ref{iv}) and
 ``$\text{III}_{\text{1}}+\text{I}$'' can lead only to
 $N=3,5,7,15,21$ and $N=3,5,15$ correspondingly. Solutions of type
 ``$\text{III}_{\varphi}+\text{I}$'' give $N=3,9$. Finally, since
 $\omega\neq0$ and it may be assumed that $X\neq1$, for solutions
 of type ``$\text{II}_{\varphi}+\text{II}_{\psi}$'' there are
 two possibilities:
 \ben
 \tag{2.1}\begin{cases}
 \cos\pi(2r_X+r_Z)=\cos\pi(2r_X+r_{\bar{Z}}),\\
 \cos\pi(2r_X-r_{Z}-\varepsilon_1/3)=\cos\pi(2r_X-r_{\bar{Z}}+\varepsilon_1/3).
 \end{cases}
 \ebn
 If $r_Z=r_{\bar{Z}}\;(\mathrm{mod}\;2\Zb)$, then the second
 equation implies that $r_Z=2r_X+(1-\varepsilon_3)/2\;
 (\mathrm{mod}\;2\Zb)$, $\varepsilon_3=\pm1$. Assume that
 $N\neq3$, then from the relation $(X-1)(Y'-Z)=Y-Z'$ we find
 \ben
 \cos\pi r_{Y'}=\varepsilon_3\Bigl(\cos2\pi r_X-\cos\pi(r_X+\varepsilon_1/3)-\cos\pi/3\Bigr).
 \ebn
 Rational solutions of this equation lead to admissible values
 $N=3,5,7,9,15$. Now if we take as the solution of the first
 equation in (2.1) $r_{\bar{Z}}=-4r_X-r_Z\;(\mathrm{mod}\;2\Zb)$,
 then from the second equation follows $r_Z=-2r_X-\varepsilon_1/3+(1-\varepsilon_3)/2\;
 (\mathrm{mod}\;2\Zb)$. Computing $Y'$ from $(X-1)(Y'-Z')=Y''-Z$,
 one finds the same values of $N$.
 \ben
 \tag{2.2}\begin{cases}
 \cos\pi(2r_X+r_Z)=\cos\pi(2r_X-r_{\bar{Z}}+\varepsilon_1/3),\\
 \cos\pi(2r_X-r_{Z}-\varepsilon_1/3)=\cos\pi(2r_X+r_{\bar{Z}}).
 \end{cases}
 \ebn
 This case is completely analogous to (2.1). \\
 \underline{Case (3)}.
 Recall that solutions of (\ref{sumcos}) relevant for
 (\ref{oddboundaux16}) should satisfy an additional constraint
 $\varepsilon_1\varphi_1+\varepsilon_2\varphi_2+\varepsilon_3\varphi_3+\varepsilon_4\varphi_4\in
 \Zb$ with some $\varepsilon_{1,2,3,4}=\pm1$.
 This condition implies that $\varphi\pm1/6$ in (\ref{via}) belong or do not belong to
 $\{\varphi_1,\varphi_2,\varphi_3,\varphi_4\}$ simultaneously,
 otherwise admissible  $N$ are $3,5,15$.
 Futhermore if we assume that $N\neq 3,5,15$, the  unordered pairs
 $(r_X+r_Z,r_X-r_Z)$ and $(r_X+r_{Z'},r_X-r_{Z'})$ can only be equivalent
 to the following: \\
 (3.1) $(2\varphi+1/3,2\varphi-1/3)$ and
 $(2\varphi+3/5,2\varphi-3/5)$,\\
 (3.2) $(2\varphi+1/3,2\varphi-1/3)$ and
 $(2\varphi+1/5,2\varphi-1/5)$,\\
 (3.3) $(2\varphi+1/5,2\varphi-1/5)$ and
 $(2\varphi+2/5,2\varphi-2/5)$.\\
 Here $\varphi\in\Qb$ and all entries in (3.1)--(3.3) are
 considered $\mathrm{mod}\;2\Zb$. Now observe that in (3.1) and
 (3.2) either $Z$ or $Z'$ is equal to $\pm1$, therefore such
 6-tuples can be excluded by Assumption~1. In the case (3.3),
 unordered pair $(Z,Z')$ is equal to $(2\cos\pi/5,2\cos 2\pi/5)$
 or $(-2\cos\pi/5,-2\cos 2\pi/5)$.\vspace{0.1cm}

 Let us now summarize the above results. If $N\neq3,5,7,9,11,15,21$, then
 $N$ relations (\ref{oddboundaux16}) can have only the following solutions: \\
 (a) with $Z$ or $Z'$ equal to $\pm1$,
 $\pm\sqrt{4-X^2}$,\\
 (b) solutions of type
 ``$\text{III}_{\varphi}+\text{III}_{\psi}$'' (and ``$\text{VI}_{\varphi}$'')
 satisfying Assumptions~1 and~2; these appear in $O_{yz}$
 at most once (resp. twice).\\
  However, under such restrictions the length of $O_{yz}$ cannot exceed $11$ because of
 Lemma~\ref{distinct} (as all $Z_k$ in  the simple cycle are distinct).
 \epf
 \begin{prop}
 Let $\omega_Y=\omega_Z=0$. Then either $N\leq15$ or the
 suborbit $O_{yz}$ has the form
 \be\label{oddboundaux21}
 \begin{cases}
 X=2\cos\pi r_X,\\
 Y_{k}=-2\cos\pi\bigl[r_X(1+2k_0-2k)+r_Z\bigr],\\
 Z_{k}=2\cos\pi\bigl[2r_X(k_0-k)+r_Z\bigr].
 \end{cases}
 \eb
 where $k_0\in \{0,1,\ldots,N-1\}$ and $r_{X,Z}\in\Qb$.
  \end{prop}
 \pf Let us consider the relation (see Fig.~5)
 \be\label{oddboundaux22}
 \cos\pi r_Y+\cos\pi r_{Y'}+\cos\pi(r_X+r_Z)+\cos\pi(r_X-r_Z)=0.
 \eb
 For $N\geq15$ ($N\geq6$ in the simple cycle case) one will always be able to find in
 $O_{yz}$ a solution with $r_{Y,Y',Z}\in\Qb$ satisfying the restrictions $Y,Y'\neq0$ and $Z\neq 0,\pm\sqrt{4-X^2}$.
 With these requirements, the solution of (\ref{oddboundaux22}) cannot be of type ``$\text{III}_{\text{1}}+\text{I}$''
 or ``$\text{III}_{\varphi}+\text{I}$'' as the relation (\ref{oddboundaux22}) does not contain zero
 cosines. Moreover one cannot have solutions of type ``$\text{II}_{\varphi}+\text{II}_{\psi}$''
 with $Y+Y'=0$ unless $X=0$, i.e. $N=2$.

 For the remaining ``$\text{II}_{\varphi}+\text{II}_{\psi}$'' solutions
 one can write
 \ben
 Y=-2\cos\pi(r_X+r_Z),\qquad Y'=-2\cos\pi(r_X-r_Z).
 \ebn
 Setting $Y=Y_{k_0}$, $Y'=Y_{k_0+1}$ we find that $\alpha$, $\beta$ in (\ref{yzsolution1})
 are given by
 \ben
 \alpha=-2\cos\pi\bigl[r_X(1+2k_0)+r_Z\bigr],\qquad \beta=2\cos\pi\bigl[2k_0 r_X+r_Z\bigr],
 \ebn
 and hence $\{Y_k\}$, $\{Z_k\}$ have the form (\ref{oddboundaux21}).

 Now we can assume that all solutions satisfying the above restrictions are equivalent
 to the quadruples (\ref{iv}). This leads
 to admissible values $N=3,5,7,15,30,42$. However, the lengths $N=30,42$ can be excluded because it is not possible to generate from
 (\ref{iv}) a sufficient number of solutions with the same value of $X$ and different
 $Z$.\\
 \textit{Example}. Checking all the quadruples (\ref{iv}) with
 $\ds X=2\cos{\pi}/{30}$ we find that there are only six possible values of
 $Z$:
 $\pm 2\cos 7\pi/30$, $\pm2\cos 11\pi/30$ and  $\pm2\cos 13\pi/30$.
  \epf

 Assume that $O_{yz}$ has the form (\ref{oddboundaux21}).
 If $\omx=0$, then from (\ref{jfr}) and (\ref{oddboundaux21}) follows that
 $\omi=0$.
 Finite orbits of the induced $\bar{\Lambda}$~action (\ref{lxyz})
 with $\omx=\omy=\omz=\omi=0$ will be
 called \textit{Cayley orbits} because in this case Jimbo-Fricke relation
 (\ref{jfr}) reduces to Cayley cubic
 \be\label{cayley}
 XYZ+X^2+Y^2+Z^2-4=0.
 \eb

  Cayley orbits admit a simple characterization, though their size can be
  arbitrarily large. To each of these orbits one can
  assign in a non-unique way a pair of rational numbers.
  Indeed, consider an arbitrary
  point $\mathbf{r}=(X,Y,Z)\in O$. It is not fixed by at least one transformation, say $x$
  (we assume that $O$ consists of more than one point).
  Lemma~\ref{xyzcosineslemma} then implies that $Y=2\cos\pi r_Y$, $Z=2\cos\pi
  r_Z$ with $r_{Y,Z}\in\Qb$. The relation (\ref{cayley}) can be rewritten as
  \ben
  \bigl(X+2\cos\pi(r_Y+r_Z)\bigr)
  \bigl(X+2\cos\pi(r_Y-r_Z)\bigr)=0,
  \ebn
  hence we may assume that $X=-2\cos\pi(r_Y+r_Z)$ (if $X=-2\cos\pi(r_Y-r_Z)$,
  start from $x(X,Y,Z)$).
  Now making one step from $(X,Y,Z)$ by $x$, $y$ and $z$ one finds
  \ben
  \begin{cases}
  X(x(\mathbf{r}))=-2\cos\pi(r_Y-r_Z),\\
  Y(y(\mathbf{r}))=2\cos\pi(r_Y+2 r_Z),\\
  Z(z(\mathbf{r}))=2\cos\pi(2r_Y+r_Z).
  \end{cases}
  \ebn
  Continuing by induction we see that for any other point
  $(X',Y',Z')\in O$ one has $X'=2\cos\pi r_{X'}$, $Y'=2\cos\pi
  r_{Y'}$, $Z'=2\cos\pi r_{Z'}$, where $r_{X',Y',Z'}\in\Qb$ and
  the denominators of  $r_{X',Y',Z'}$ are divisors of the common
  denominator of $r_Y$ and $r_Z$. Lemma~\ref{distinct} then guarantees
  that $O$ is finite.
  \begin{prop}
  Let $\omy=\omz=0$. If $O_{yz}$ has the form
  (\ref{oddboundaux21}) and $\omx\neq0$, then $N\leq12$.
  \end{prop}
  \pf
   Let us make one step by $x$ from each point of $O_{yz}$ (see Fig.~6). Using (\ref{oddboundaux21}),
  from the relations $\omega_X=X+X_k+Y_k Z_k=X+\bar{X}_k+Y_{k+1} Z_{k}$ one
  finds
  \be\label{oddboundaux23}
  \begin{array}{lll}
  \omx&=&X_k-2\cos\pi\bigl[r_X(4k_0-4k+1)+2r_Z\bigr]=\vspace{0.1cm}\\
  &=&\bar{X}_k-2\cos\pi\bigl[r_X(4k_0-4k-1)+2r_Z\bigr],
 \end{array}
  \eb
  for any $k=0,1,\ldots,N-1$. If the point $(X_k,Y_k,Z_k)$ is good then by Lemma~\ref{xyzcosineslemma}
  \be\label{xrational}
  X_k=2\cos\pi r_{X_k},\qquad\qquad  r_{X_k}\in\Qb.
  \eb
  It can be bad in two cases:\vspace{0.1cm}\\
  (1) The graph of $O_{yz}$ is a line, $(X,Y_k,Z_k)$ corresponds to one of its end vertices and $X_k=X$. Since $N>1$,
  $X_k$ still has the form (\ref{xrational}).\\
  (2) $(X_k,Y_k,Z_k)$ is fixed by the transformations $y$ and $z$. Then from (\ref{yzrecursion}) follows that either
  $X_k=\pm2$  or $Y_k=Z_k=0$. In the latter case, however, the condition $N>1$ is violated since the whole orbit $O$
  consists of only two points $(X,0,0)$ and $(\omx-X,0,0)$.\vspace{0.1cm}\\
  Thus all $X_k$ and $\bar{X}_k$ have the form (\ref{xrational}) and
  the solutions of (\ref{oddboundaux23}) are classified by Lemma~\ref{lemcos}.
    \begin{figure}[!h]
 \begin{center}
 \resizebox{9cm}{!}{
 \includegraphics{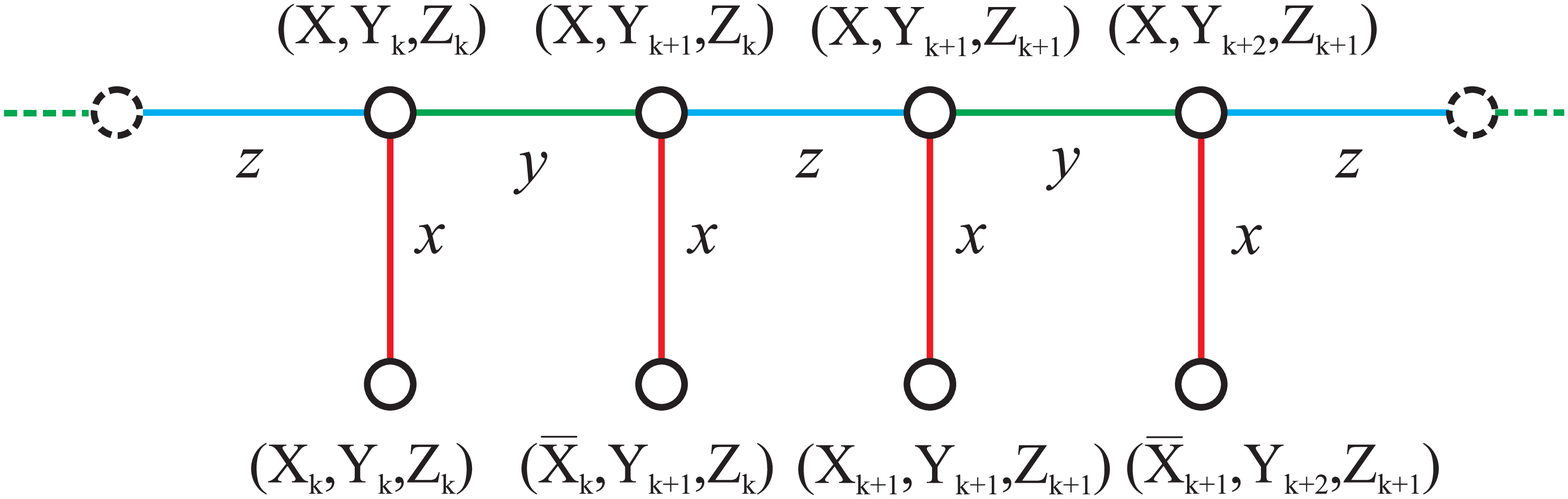}} \\   Fig. 6
 \end{center}
 \end{figure}

  Introduce $2N$ quantities $W_0,\ldots,W_{2N-1}$ defined by
  \ben
  W_{2k}=X_k-\omx,\quad W_{2k+1}=\bar{X}_k-\omx,\qquad
  k=0,\ldots,N-1.
  \ebn
  Obviously, $W_l=2\cos\pi\bigl[r_X(1+4k_0-2l)+2r_Z\bigr]$. We now want to show that
  the number of coinciding $W_l$ cannot exceed 4. Indeed, fix some
  $l$, then $W_{l'}=W_l$ implies that (a)
  $l'-l=0\;\mathrm{mod}\;N$ or (b) $r_X(1+4k_0-l-l')+2r_Z\in\Zb$.
  The former case leads to one compatible $W_{l'}$, while the latter
  gives at most two: if $l'_1$ and $l'_2$ satisfy (b), then
  necessarily $l'_1-l'_2=0\;\mathrm{mod}\;N$.

  In the proof of Propositions~\ref{evenbound} and~\ref{oddbound1}
  we have shown that the maximal number of ordered pairs $(\cos\pi r_1,\cos\pi
  r_2)$,  $r_{1,2}\in\Qb$ such that $\cos\pi r_1+\cos\pi
  r_2=\mathrm{const}\neq0$
  is equal to 6. Hence the number of distinct possible values for
  all $W_l$'s cannot exceed $6$ and the total number of $W_l$'s, equal to
  $2N$, cannot exceed 24.
  \epf

  Let us summarize the results of this subsection. Given a finite orbit $O$, common
  coordinate $X$ of all points of any 2-colored suborbit $O_{yz}\subset O$
  of length $N>1$  has the form $ X=2\cos{\pi n_X}/{N}$, $0<n_X<N$, where $N$ and $n_X$ are coprime.
  Unless $\omx=\omy=\omz=\omi=0$, one has a number of restrictions on possible values of $N$ and $n_X$
  listed in Table~2. These restrictions imply in particular that
  $X$ can take only a finite number of explicitly defined values.
  In the next subsection, we use this observation to construct an
  exhaustive search algorithm giving all finite orbits of (\ref{lxyz}).

  \begin{center}
  \begin{tabular}{c||c|c|}
  & restrictions on $N$, $n_X$ & {$\begin{array}{c}\text{number of} \\ \text{possible
  }X\end{array}$}\\
  \hline\hline
  $\omy^2\neq\omz^2$ & $N\leq10$, $n_X$ odd and even &
  31 \\
  \hline
  $\omy=\omz\neq0$ & $\begin{array}{c} N\leq10\text{, }n_X\text{ odd and even,}\\
  N=11,15,21\text{, }n_X\text{ odd}\end{array}$ & 46 \\
  \hline
  $\begin{array}{c} \omy=\omz=0\text{ with}\\ \omx\neq0\text{ or }\omi\neq0\end{array}$ & $N\leq15$, $n_X\text{ odd and even}$ & 71 \\
  \hline
  \end{tabular}\vspace{0.2cm}\\
  Table 2: Restrictions on possible values of $X$ for $N>1$.
  \end{center}

  \subsection{Search algorithm}
  Let $O\subset \Cb^3$ be a finite orbit of the
  induced $\bar{\Lambda}$~action~(\ref{lxyz}) consisting of more than one point.
  Since we are interested in nonequivalent orbits, it can be asssumed
  that the parameters $\omega_{X,Y,Z,4}$ satisfy one of the following sets of constraints: \\
  (A) $\omx^2\neq\omy^2\neq \omz^2$,\\
  (B) $\omx^2\neq\omy^2$, $\omy=\omz\neq0$,\\
  (C) $\omx\neq0$, $\omy=\omz=0$,\\
  (D) $\omx=\omy=\omz\neq0$,\\
  (E) $\omx=\omy=\omz=0$, $\omi\neq0$,\\
  (F) $\omx=\omy=\omz=\omi=0$.\\
  In what follows, the case (F) will be omitted, as all finite orbits with such
  parameter values
  have already been described above.
  \begin{defn}
  Let $\mathbf{r}=(X,Y,Z)$ be a point in $O$. Its coordinate $X$
  (or $Y$, $Z$) will be called \textit{good} if $\mathbf{r}$ is not fixed
  by at least one of the transformations $y$ and $z$ (resp. $x$ and $z$, $x$ and
  $y$).
  \end{defn}
  \begin{rmk} All coordinates of a good point are good. If
  $\mathbf{r}$ is a bad point, e.g. fixed by $y$ and $z$ but not by
  $x$, then it has good coordinates $Y$ and $Z$.
  \end{rmk}

  Define
  three finite sets of numbers (cf. Table~2):
  \begin{align*}
  \mathcal{S}_1&=\left\{2\cos\frac{\pi n}{N}\,\bigl|\,1<N\leq
  10,\;n\text{ odd and even}\right\},\\
   \mathcal{S}_2&=\left\{2\cos\frac{\pi n}{N}\,\bigl|\,1<N\leq 10,
   \;n\text{ odd and even};\,N=11,15,21,\; n\text{
   odd}\right\},\\
   \mathcal{S}_3&=\left\{2\cos\frac{\pi n}{N}\,\bigl|\,1<N\leq
   15,\;n\text{ odd and even}\right\}.
 \end{align*}
 In all three cases $n$ is supposed to be coprime with $N$ and
 $0<n<N$. Now the results of the previous subsection imply that good
 coordinates of any point $\mathbf{r}\in O$ belong to one of these
 lists according to Table~3.
 \begin{center}
 \begin{tabular}{c||c|c|c|}
  & good $X$ & good $Y$ & good $Z$ \\
 \hline\hline
 (A) & $\mathcal{S}_1$ & $\mathcal{S}_1$ & $\mathcal{S}_1$ \\
 \hline
 (B) & $\mathcal{S}_2$ & $\mathcal{S}_1$ & $\mathcal{S}_1$ \\
 \hline
 (C) & $\mathcal{S}_3$ & $\mathcal{S}_1$ & $\mathcal{S}_1$ \\
 \hline
 (D) & $\mathcal{S}_2$ & $\mathcal{S}_2$ & $\mathcal{S}_2$ \\
 \hline
 (E) & $\mathcal{S}_3$ & $\mathcal{S}_3$ & $\mathcal{S}_3$ \\
 \hline
 \end{tabular}\vspace{0.2cm}\\
 Table~3: Admissible values of good coordinates
 \end{center}

 Any orbit $O$ is completely defined by a point $\mathbf{r}\in O$
 and the parameter triple $\boldsymbol{\omega}=(\omx,\omy,\omz)$.
 Equivalently, instead of $\boldsymbol{\omega}$ one can use three
 points $x(\mathbf{r})$, $y(\mathbf{r})$, $z(\mathbf{r})$ (some of them can coincide with $\mathbf{r}$). Denote
 \be\label{xyzpr}
 {X'}=X(x(\mathbf{r})),\qquad {Y'}=Y(y(\mathbf{r})),\qquad
 {Z'}=Z(z(\mathbf{r})),
 \eb
 then we have
 \be\label{omegas}
 \omx=X+{X'}+YZ,\qquad \omy=Y+{Y'}+XZ,\qquad \omz=Z+{Z'}+XY.
 \eb
 \begin{defn}
 Let $\mathbf{r}$ be a good point in a finite orbit $O$. The set of four
 points $\{\mathbf{r},x(\mathbf{r}),y(\mathbf{r}),z(\mathbf{r})\}$
 will be called a good generating configuration (GGC)  for $O$ if at
 least two of three points  $x(\mathbf{r})$, $y(\mathbf{r})$,
 $z(\mathbf{r})$ are good.
 \end{defn}
 \begin{lem}
 Let $O$ be a finite orbit that does not contain a GGC. Then $\Sigma(O)$ can only be equivalent (up to permutations of colors)
 to one of the four graphs shown in Fig.~7.
 \end{lem}
 \begin{figure}[!h]
 \begin{center}
 \resizebox{10cm}{!}{
 \includegraphics{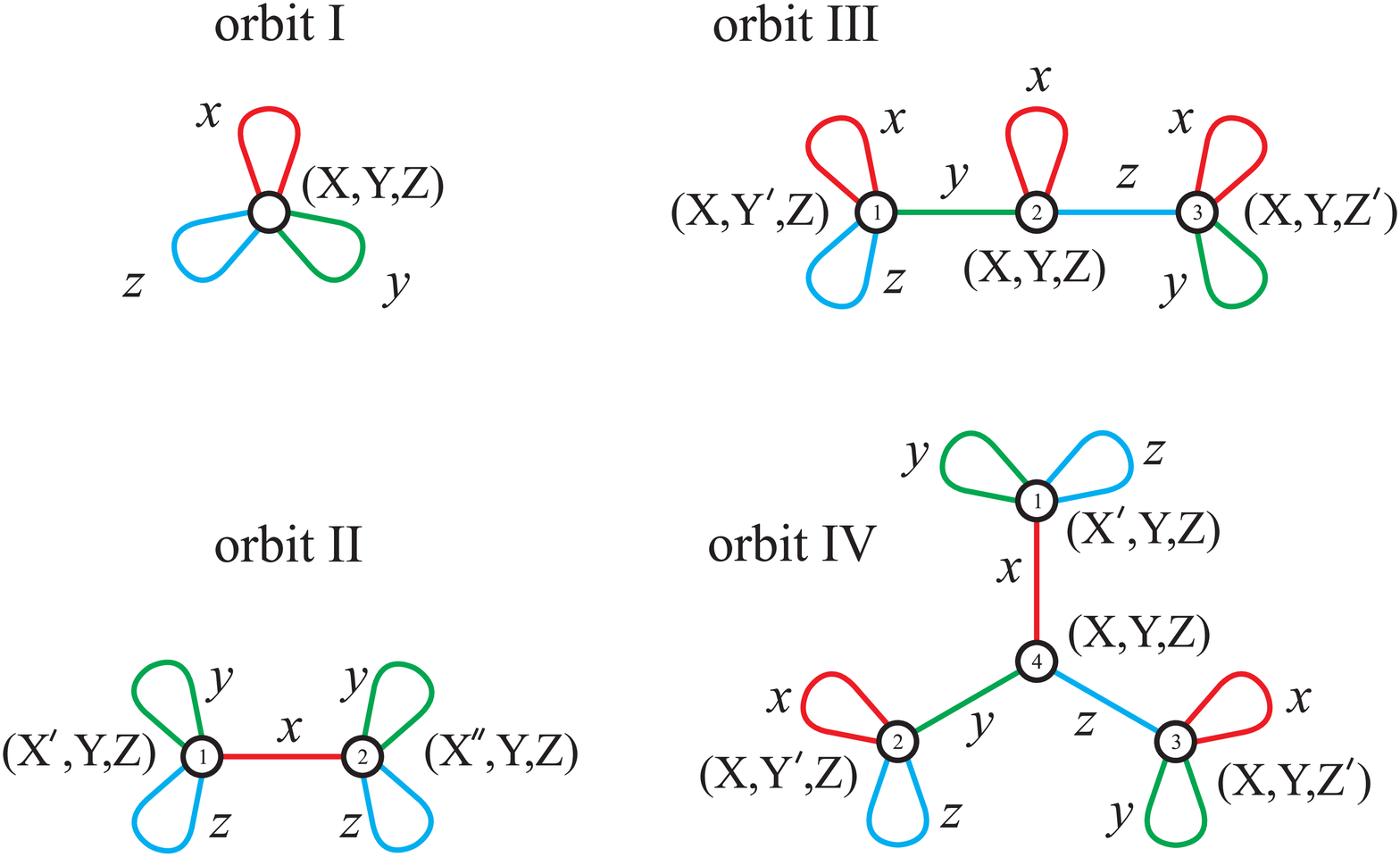}} \\  Fig. 7: Four orbits without GGCs
 \end{center}
 \end{figure}
 \pf
 If $O$ contains more than 2 points, then at least one of them is
 good. Denoting this point by $\mathbf{r}$, we can assume that
 $y(\mathbf{r})$ and $z(\mathbf{r})$ are bad. Now if
 $x(\mathbf{r})=\mathbf{r}$, then one obtains orbit~III. The case
 when $x(\mathbf{r})\neq\mathbf{r}$ is bad corresponds to orbit~IV. Finally, if
 $x(\mathbf{r})\neq \mathbf{r}$ is another good point, then by
 assumptions of the Lemma
 the points  $y(x(\mathbf{r}))$ and  $z(x(\mathbf{r}))$ are bad,
 and $\Sigma(O)$ is given by the 6-vertex graph represented in
 Fig.~8.
  \begin{figure}[!h]
 \begin{center}
 \resizebox{5cm}{!}{
 \includegraphics{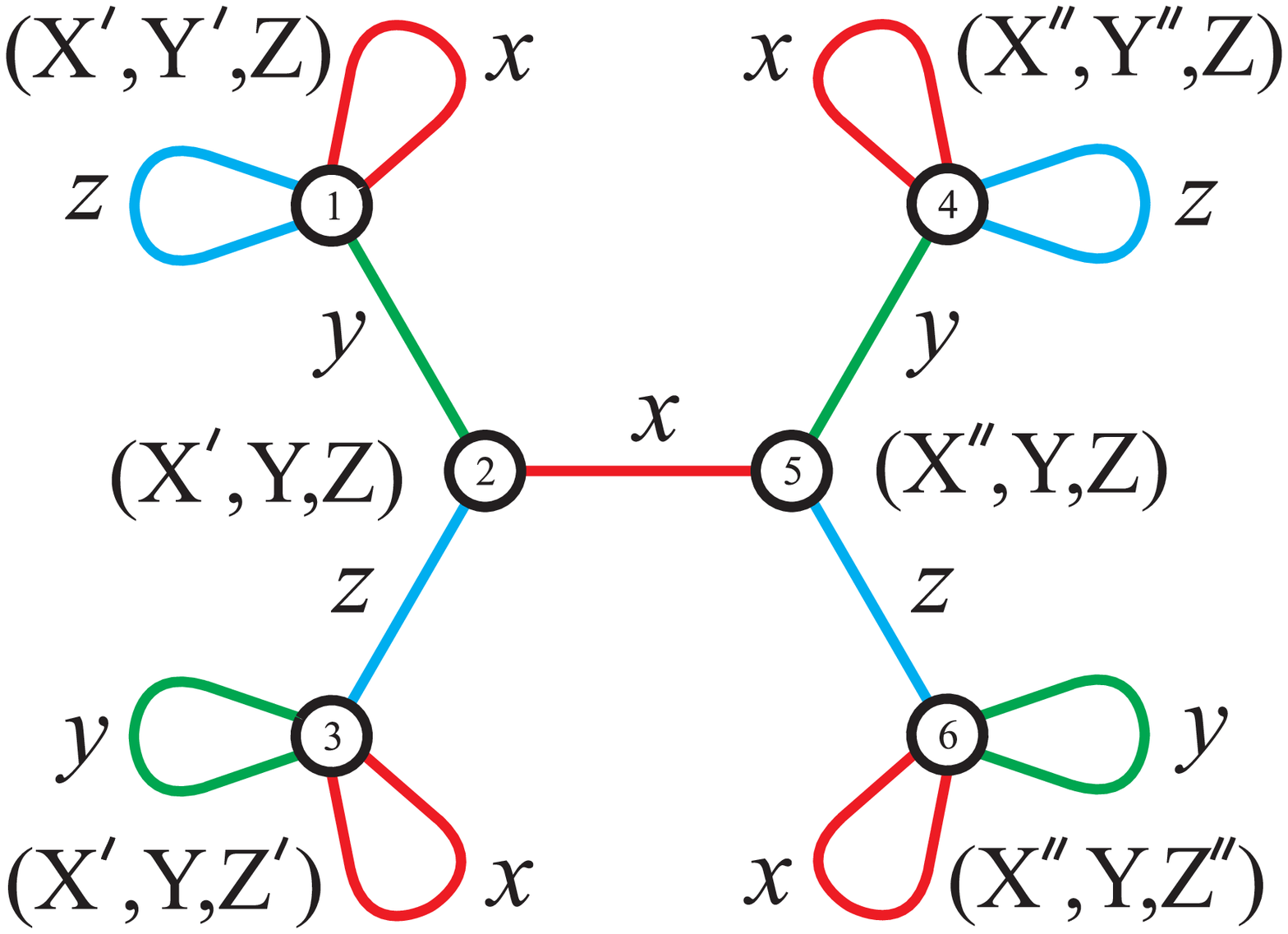}} \\ Fig. 8: 6-vertex graph without GGCs
 \end{center}
 \end{figure}

 It turns out, however, that this last graph is forbidden. To see
 this, note that $yz$-suborbits 1-2-3 and 4-5-6 both have length
 3, therefore $X'$ and $X''$ are equal to $\pm1$. Since $X'\neq
 X''$, one can set $X'=1$, $X''=-1$. Then from the relations
 corresponding to $y$- and $z$-edges,
 \begin{align*}
 \omy&=Y+Y'+X'Z=2Y+X'Z'=2Y+X''Z''=Y+Y''+X''Z,\\
 \omz&=Z+Z'+X'Y=2Z+X'Y'=2Z+X''Y''=Z+Z''+X''Y,
 \end{align*}
 it follows that $Y=-Z'=Z''$ and $Z=-Y'=Y''$. Self-loops of color $x$ at the points 1, 3, 4 and 6
 in turn imply that $\omx=0$, $Y^2=Z^2=2$. However, this is
 incompatible with the $x$-edge 2-5, which gives $\omx=YZ$.
 \epf

 The orbits of (\ref{lxyz}) with graphs~I--IV are completely described by the
 following:
 \begin{lem}\label{orb14} 1.~Orbits of type~I consist of one point
 $(X,Y,Z)\in\Cb^3$.
 The parameters $\omega_{X,Y,Z,4}$
 are given by
 \be\label{tetrahr}
 \omx=2X+YZ,\qquad\omy=2Y+XZ,\qquad \omz=2Z+XY,
 \eb
 \be\label{tetrahr2}
 \omi=4+2XYZ+X^2+Y^2+Z^2.
 \eb
 2.~Any orbit of type~II is equivalent to an orbit consisting of
 2 points $(X',0,0)$ and $(X'',0,0)$, where $X',X''\in \Cb$, $X'\neq X''$
 and $\omx=X'+X''$, $\omy=\omz=0$, $\omi=4+X'X''$.\vspace{0.1cm}\\
 3.~Any orbit of type~III is equivalent to an orbit consisting of
 3 points $(1,0,0)$, $(1,\omega,0)$, $(1,0,\omega)$, where
 $\omega\in\Cb^*$ and  $\omx=2$,
 $\omy=\omz=\omega$, $\omi=5$.\vspace{0.1cm}\\
 4.~Any orbit of type~IV is equivalent to an orbit consisting of
 4 points $(1,1,1)$, $(\omega-2,1,1)$, $(1,\omega-2,1)$,
 $(1,1,\omega-2)$, where $\omega\in\Cb$, $\omega\neq3$ and
 $\omx=\omy=\omz=\omega$, $\omi=3\omega$.
 \end{lem}
 \pf Statement~1 is obvious ($\omi$ is determined from (\ref{jfr})),
  hence we start with orbits of type~II. In this case,
 since $xy$- and $xz$-suborbits 1-2 have length 2, one finds $Y=Z=0$. From the
 relations corresponding to the self-loops then follows $\omy=\omz=0$.

 For orbits of type~III, $xy$-suborbit 1-2 and $xz$-suborbit 2-3 both have length 2, therefore  $Y=Z=0$.
 Similarly, $yz$-suborbit 1-2-3 has length 3 and thus $X=\pm1$. Since the simultaneous change
 of signs of e.g. $\omx$, $\omy$, and also $X$- and $Y$-coordinates of all points leads to an equivalent orbit,
 one can set $X=1$, and then $x$-self-loop at the point 2 gives $\omx=2$.
 At last, $y$- and $z$-edges of the graph imply that $\omy=\omz=Y'=Z'$.

 In graph~IV, $xy$-suborbit 1-4-2, $xz$-suborbit 1-4-3 and $yz$-suborbit 2-4-3 have length 3,
 therefore $X$, $Y$ and $Z$ are equal to $\pm1$. It can be assumed that either
 (a) $X=Y=Z=1$ or (b) $X=Y=Z=-1$. In the case (a), $y$- and $z$-self-loops at the point 1 imply that
 $\omy=\omz=2+X'$, hence by symmetry
 \ben
 \omx=\omy=\omz=2+X'=2+Y'=2+Z',
 \ebn
 and the relations corresponding to the edges 1-4, 2-4 and 3-4 are
 satisfied automatically. In the case (b), one similarly finds
 $\omx=\omy=\omz=-2-X'=-2-Y'=-2-Z'$, but e.g. the relation 1-2
 gives $\omx=X'$. Thus $X'=-1$ and we obtain a contradiction.
 \epf

 Unless $\omx=\omy=\omz=\omi=0$, one has only a finite number of GGCs (and hence only a finite number of finite
 orbits different from I--IV). Indeed, these configurations can be of two types:\vspace{0.1cm}\\
 \underline{Type (i)}. All four points $\mathbf{r},x(\mathbf{r}),y(\mathbf{r}),z(\mathbf{r})\in O$
 are good. In this case six coordinates $X$, $Y$, $Z$, $X'$, $Y'$,
 $Z'$ (defined by (\ref{xyzpr})) are good, hence each of them
 can take only a finite number of values, as specified in
 Table~3.\vspace{0.1cm}\\
 \underline{Type (ii)}. One of three points $x(\mathbf{r}), y(\mathbf{r}),
 z(\mathbf{r})\in O$ is bad. Suppose e.g. that $x(\mathbf{r})$ is
 bad, then $X$, $Y$, $Z$, $Y'$,
 $Z'$ are good coordinates, but $X'$ is not. However, since $x(\mathbf{r})$ is fixed by $y$ and
 $z$, we have the equations
 \be\label{xprbad}
 \begin{cases}
 \omy=2Y+X'Z=Y+Y'+XZ,\\
 \omz=2Z+X'Y=Z+Z'+XY.
 \end{cases}
 \eb
 Unless $Y=Z=0$, one can use (\ref{xprbad}) to express $X'$ in
 terms of good coordinates. Also notice that $Y$, $Z$, $Y'$,
 $Z'$ should satisfy an additional relation
 \be\label{addc}
 Y(Y-Y')=Z(Z-Z').
 \eb
 On the other hand if $Y=Z=0$, then (\ref{xprbad}) implies that
 $\omy=\omz=0$, the orbit $O$ is of type~II and in particular it does not contain a
 GGC.\vspace{0.1cm}\\

 Let us now describe in more detail the sets of good coordinates generating
 all possible candidates for finite orbits, different from
 orbits~I--IV and those of Cayley type:\vspace{0.1cm}\\
 \underline{Class 1 (A-i)}.
 Here one has $31^6\approx 10^9$  GGCs, corresponding to
 all possible $X, Y, Z$, $X', Y', Z'\in \mathcal{S}_1$. Since we are interested in
 nonequivalent orbits, it can be assumed that either $0\leq X\leq Y\leq Z$ or $0\geq X\geq Y\geq Z$,
 and then the above number reduces to $16\cdot 17\cdot 18\cdot 31^3/3-1=48\,618\,911$.
 We do not exclude the remaining equivalent GGCs for simplicity of the algorithm.\vspace{0.1cm}\\
 \underline{Class 2 (A-ii,B-ii,C-ii,D-ii,E-ii)}. In this case it is convenient to deal not only with $\omega_{X,Y,Z}$ satisfying one of the conditions (A)--(E), but also with equivalent parameter triples. One can then
 assume that $x(\mathbf{r})$ is bad and $0\leq Y\leq Z$, $Z>0$. Since $Z'$ can now be determined from
 (\ref{addc}), the whole orbit is completely fixed by four good coordinates $X,Y,Z,Y'$, taking their values in
 the set
 \ben
 \mathcal{S}_4=\left\{2\cos\frac{\pi n}{N}\,\bigl|\,1<N\leq
  15,N=21,\;n\text{ odd and even}\right\},\\
 \ebn
 consisting of $83$ elements. The total number of configurations to be checked is therefore equal to
 $41\cdot 22\cdot 83^2=6\,213\,878$.
 \vspace{0.1cm}\\
 \underline{Class 3 (B-i,C-i)}. Here we use good coordinates
 $X,X'\in\mathcal{S}_4$, $Y,Y',Z\in\mathcal{S}_1$, while $Z'$ is computed from
 \ben
 Z'=\left(Y+Y'+XZ\right)-Z-XY,
 \ebn
 and it can be assumed that $0\leq|Y|\leq Z$. This gives $16^2\cdot 31\cdot 83^2=54\,671\,104$ configurations, from which we should choose
 only those with $Z'\in \mathcal{S}_1$.\vspace{0.1cm}\\
 \underline{Class 4 (D-i,E-i)}. These orbits are completely determined by
 $X,X',Y,Z\in\mathcal{S}_4$.
 Since $x(\mathbf{r})$, $y(\mathbf{r})$, $z(\mathbf{r})$ are good,
 it can be assumed that $X\leq Y\leq Z$, which leads to $83^2\cdot 84\cdot 85/6=8\,197\,910$
 possibilities.\vspace{0.1cm}\\

  In order to check which generating sets do actually lead to finite orbits, one can use the following algorithm:
  \begin{enumerate}
  \item[1.] Consider any generating set from the above as a set $\mathcal{P}$ of known orbit
  \textit{points} and known \textit{adjacency relations} between them. E.g. if it is known by construction
   that   $x(\mathbf{r})=\mathbf{r'}$ for some $\mathbf{r},\mathbf{r'}\in \mathcal{P}$,  we will say that
  $\mathbf{r'}$ is a known $x$-neighbor of $\mathbf{r}$ and vice versa. Thus any point $\mathbf{r}\in\mathcal{P}$
  has  at most 3 known neighbors, corresponding to $x$-, $y$- and $z$-edges originating from $\mathbf{r}$.
  \item[2.] If the set is characterized by $\omx=\omy=\omz=\omi=0$, the algorithm stops (the only finite orbits with such parameters
  are Cayley orbits).
  \item[3.] Using $\mathcal{P}$, construct the set $\mathcal{P}^u$ of points with at least one unknown neighbor.
  \item[4.] Choose an arbitrary point $\mathbf{r}=(X,Y,Z)\in\mathcal{P}^u$. Assume for definiteness
  that its $x$-neighbor $x(\mathbf{r})=(X',Y,Z)$ is unknown. Then compute $X'$ and proceed as follows:
  \begin{enumerate}
  \item[4.1.] If $(X',Y,Z)\in\mathcal{P}^u$, then add the appropriate $x$-adjacency relation to $\mathcal{P}$, update $\mathcal{P}^u$ and go to Step~5, else
      \item[4.2.] If $X'$ has a good value (in practice it is sufficient to require $X'\in\mathcal{S}_4$), then add $(X',Y,Z)$ and
      the appropriate $x$-adjacency relation to $\mathcal{P}$, update $\mathcal{P}^u$ and go to Step~4, else
      \item[4.3.] If $2Y+X'Z=\omy$ and $2Z+X'Y=\omz$, then add $(X',Y,Z)$ and
      the appropriate $x$-, $y$- and $z$-adjacency relations to $\mathcal{P}$, update $\mathcal{P}^u$ and go to Step~5, else the algorithm stops (the orbit cannot be finite).
  \end{enumerate}
  \item[5.] If $\mathcal{P}^u$ is empty, the algorithm stops
  (the orbit is finite and its points are given by $\mathcal{P}$), otherwise go to Step~4.
  \end{enumerate}
 \begin{rmk} It is easy to see that the algorithm stops after
 a finite number of steps. Indeed,
 the total number $N_{g}$ of good points in any finite orbit which is not of Cayley type
 cannot exceed $71^2\cdot2=10\,082$, while the number of bad points cannot exceed $N_g+2$.
 \end{rmk}

 \subsection{List of finite orbits}
 We have implemented the above algorithm with a computer program written in C~language. The check of all generating
 sets took less than 10 minutes on a usual 1.7GHz desktop  computer. It turned out that there are only 45 nonequivalent finite
 exceptional orbits, different from orbits~I--IV and Cayley orbits. We describe these orbits in Table~4 by indicating one of the
 orbit points
 \ben
 (X,Y,Z)=\bigl(2\cos\pi r_X,2\cos\pi r_Y, 2\cos\pi r_Z\bigr),
 \ebn
 and the parameter triple $(\omx,\omy,\omz)$. For further convenience, we also include the value of $4-\omi$, computed
 from the Jimbo-Fricke relation (\ref{jfr}). The graphs of exceptional $\bar{\Lambda}$ orbits are shown in Fig.~9--11 (marked
 vertices correspond to the points listed in Table~4).

 Our results can now be summarized in
 \begin{thm}\label{theorem1}
 The list of all nonequivalent finite orbits of the induced $\bar{\Lambda}$ action (\ref{lxyz}) consists of the following:
 \begin{itemize}
 \item four orbits I--IV, described in Lemma~\ref{orb14};
 \item 45 exceptional orbits listed in Table~4;
 \item Cayley orbits; all of these can be generated from the points
 \ben
 \bigl(-2\cos\pi(r_Y+r_Z),2\cos\pi r_Y,2\cos\pi r_Z\bigr),\qquad r_{Y,Z}\in\Qb
 \ebn
 with $\omx=\omy=\omz=0$ (the relation $\omi=0$ is satisfied automatically).
 \end{itemize}
 \end{thm}
 \begin{rmk}\label{notsplit}
 Note that the graphs of orbits~I--IV
 and of all exceptional orbits except orbits~30, 43--45 contain self-loops. It means in particular that these
 orbits do not split under the action of non-extended modular group~$\Lambda$. In fact the last statement holds
 for orbits~30, 43--45 as well, because in all four cases the orbit graphs contain simple cycles with an odd number
 of edges.
 \end{rmk}

   \newpage
  \scriptsize
  \begin{center}
  \begin{tabular}{|c|c|c|c|}
  \hline  & size & $(\omx,\omy,\omz,4-\omi)$ & $(r_X,r_Y,r_Z)$ \\  \hline
  \hline 1 & 5 & $(0,1,1,0)$ & $(2/3,1/3,1/3)$ \\
  \hline 2 & 5 & $(3,2,2,-3)$ & $(1/3,1/3,1/3)$ \\
  \hline 3 & 6 & $(1,0,0,2)$ & $(1/2,1/3,1/3)$ \\
  \hline 4 & 6 & $(\sqrt{2},0,0,1)$ & $(1/4,1/3,3/4)$ \\
  \hline 5 & 6 & $(3,2\sqrt{2},2\sqrt{2},-4)$ & $(1/2,1/4,1/4)$ \\
  \hline 6 & 6 & $\left(1-\sqrt{5},\frac{3-\sqrt{5}}{2},\frac{3-\sqrt{5}}{2},-2+\sqrt{5}\right)$ & $(4/5,1/3,1/3)$ \\
  \hline 7 & 6 & $\left(1+\sqrt{5},\frac{3+\sqrt{5}}{2},\frac{3+\sqrt{5}}{2},-2-\sqrt{5}\right)$ & $(2/5,1/3,1/3)$ \\
  \hline 8 & 7 & $(1,1,1,0)$ & $(1/2,1/2,1/2)$ \\
  \hline 9 & 8 & $(2,0,0,0)$ & $(0,1/3,2/3)$ \\
  \hline 10 & 8 & $(1,\sqrt{2},\sqrt{2},0)$ & $(1/2,1/2,1/2)$ \\
  \hline 11 & 8 & $\left(\frac{3+\sqrt{5}}{2},1,1,-\frac{\sqrt{5}+1}{2}\right)$ & $(1/3,1/2,1/2)$ \\
  \hline 12 & 8 & $\left(\frac{3-\sqrt{5}}{2},1,1,\frac{\sqrt{5}-1}{2}\right)$ & $(1/3,1/2,1/2)$ \\
  \hline 13 & 9 & $\left(2-\sqrt{5},2-\sqrt{5},2-\sqrt{5},\frac{5\sqrt{5}-7}{2}\right)$ & $(4/5,3/5,3/5)$ \\
  \hline 14 & 9 & $\left(2+\sqrt{5},2+\sqrt{5},2+\sqrt{5},-\frac{5\sqrt{5}+7}{2}\right)$ & $(2/5,1/5,1/5)$\\
  \hline 15 & 10 & $(1,0,0,1)$ & $(1/3,1/3,2/3)$\\
  \hline 16 & 10 & $\left(3-\sqrt{5},3-\sqrt{5},3-\sqrt{5},\frac{7\sqrt{5}-11}{2}\right)$ & $(3/5,3/5,3/5)$\\
  \hline 17 & 10 & $\left(3+\sqrt{5},3+\sqrt{5},3+\sqrt{5},-\frac{7\sqrt{5}+11}{2}\right)$ & $(1/5,1/5,1/5)$\\
  \hline 18 & 10 & $\left(-\frac{\sqrt{5}-1}{2},-\frac{\sqrt{5}-1}{2},-\frac{\sqrt{5}-1}{2},0\right)$ & $(1/2,1/2,1/2)$ \\
  \hline 19 & 10 & $\left(\frac{\sqrt{5}+1}{2},\frac{\sqrt{5}+1}{2},\frac{\sqrt{5}+1}{2},0\right)$ & $(1/2,1/2,1/2)$ \\
  \hline 20 & 12 & $(0,0,0,3)$ & $(2/3,1/4,1/4)$ \\
  \hline 21 & 12 & $(1,0,0,2)$ & $(0,1/4,3/4)$ \\
  \hline 22 & 12 & $(2,\sqrt{5},\sqrt{5},-2)$ & $(1/5,2/5,2/5)$ \\
  \hline 23 & 12 & $\left(\frac{3+\sqrt{5}}{2},\frac{\sqrt{5}+1}{2},\frac{\sqrt{5}+1}{2},-\sqrt{5}\right)$ & $(2/5,2/5,2/5)$ \\
  \hline 24 & 12 & $\left(\frac{3-\sqrt{5}}{2},-\frac{\sqrt{5}-1}{2},-\frac{\sqrt{5}-1}{2},\sqrt{5}\right)$ & $(4/5,4/5,4/5)$ \\
  \hline 25 & 12 & $\left(\frac{\sqrt{5}+1}{2},\frac{\sqrt{5}-1}{2},1,0\right)$ & $(1/2,1/2,1/2)$ \\
  \hline 26 & 15 & $\left(\frac{3-\sqrt{5}}{2},\frac{3-\sqrt{5}}{2},\frac{3-\sqrt{5}}{2},\sqrt{5}-1\right)$ & $(1/2,3/5,3/5)$ \\
  \hline 27 & 15 & $\left(\frac{3+\sqrt{5}}{2},\frac{3+\sqrt{5}}{2},\frac{3+\sqrt{5}}{2},-\sqrt{5}-1\right)$ & $(1/2,1/5,1/5)$ \\
  \hline 28 & 15 & $\left(\frac{5-\sqrt{5}}{2},1-\sqrt{5},1-\sqrt{5},\frac{3\sqrt{5}-5}{2}\right)$ & $(3/5,4/5,4/5)$ \\
  \hline 29 & 15 & $\left(\frac{5+\sqrt{5}}{2},1+\sqrt{5},1+\sqrt{5},-\frac{3\sqrt{5}+5}{2}\right)$ & $(1/5,2/5,2/5)$ \\
  \hline 30 & 16 & $(0,0,0,2)$ & $(2/3,2/3,2/3)$ \\
  \hline 31 & 18 & $(2,2,2,-1)$ & $(0,1/5,3/5)$ \\
  \hline 32 & 18 & $\left(1-2\cos2\pi/7,1-2\cos2\pi/7,1-2\cos2\pi/7,4\cos2\pi/7\right)$ & $(6/7,5/7,5/7)$ \\
  \hline 33 & 18 & $\left(1-2\cos4\pi/7,1-2\cos4\pi/7,1-2\cos4\pi/7,4\cos4\pi/7\right)$ & $(2/7,3/7,3/7)$ \\
  \hline 34 & 18 & $\left(1-2\cos6\pi/7,1-2\cos6\pi/7,1-2\cos6\pi/7,4\cos6\pi/7\right)$ & $(4/7,1/7,1/7)$ \\
  \hline 35 & 20 & $\left(\frac{3-\sqrt{5}}{2},0,0,1+\sqrt{5}\right)$ & $(0,1/3,2/3)$ \\
  \hline 36 & 20 & $\left(\frac{3+\sqrt{5}}{2},0,0,1-\sqrt{5}\right)$ & $(0,1/3,2/3)$ \\
  \hline 37 & 20 & $\left(1,-\frac{\sqrt{5}-1}{2},-\frac{\sqrt{5}-1}{2},\frac{\sqrt{5}+1}{2}\right)$ & $(2/3,3/5,3/5)$ \\
  \hline 38 & 20 & $\left(1,\frac{\sqrt{5}+1}{2},\frac{\sqrt{5}+1}{2},-\frac{\sqrt{5}-1}{2}\right)$ & $(2/3,1/5,1/5)$ \\
  \hline 39 & 24 & $(1,1,1,1)$ & $(1/5,1/2,1/2)$ \\
  \hline 40 & 30 & $\left(-\frac{\sqrt{5}+1}{2},0,0,\frac{3-\sqrt{5}}{2}\right)$ & $(2/3,2/3,2/3)$ \\
  \hline 41 & 30 & $\left(\frac{\sqrt{5}-1}{2},0,0,\frac{3+\sqrt{5}}{2}\right)$ & $(2/3,2/3,2/3)$ \\
  \hline 42 & 36 & $(1,0,0,2)$ & $(0,1/5,4/5)$ \\
  \hline 43 & 40 & $\left(0,0,0,\frac{5-\sqrt{5}}{2}\right)$ & $(2/5,2/5,2/5)$ \\
  \hline 44 & 40 & $\left(0,0,0,\frac{5+\sqrt{5}}{2}\right)$ & $(4/5,4/5,4/5)$ \\
  \hline 45 & 72 & $(0,0,0,3)$ & $(1/2,1/5,2/5)$ \\
  \hline
  \end{tabular}\vspace{0.1cm}\\
  \small
  Table~4: Exceptional finite $\bar{\Lambda}$ orbits
  \end{center}

  \normalsize
  \pagebreak[4]
 \begin{figure}[!h]
 \begin{center}
 $\;$\\ $\;$\\
 \resizebox{15cm}{!}{
 \includegraphics{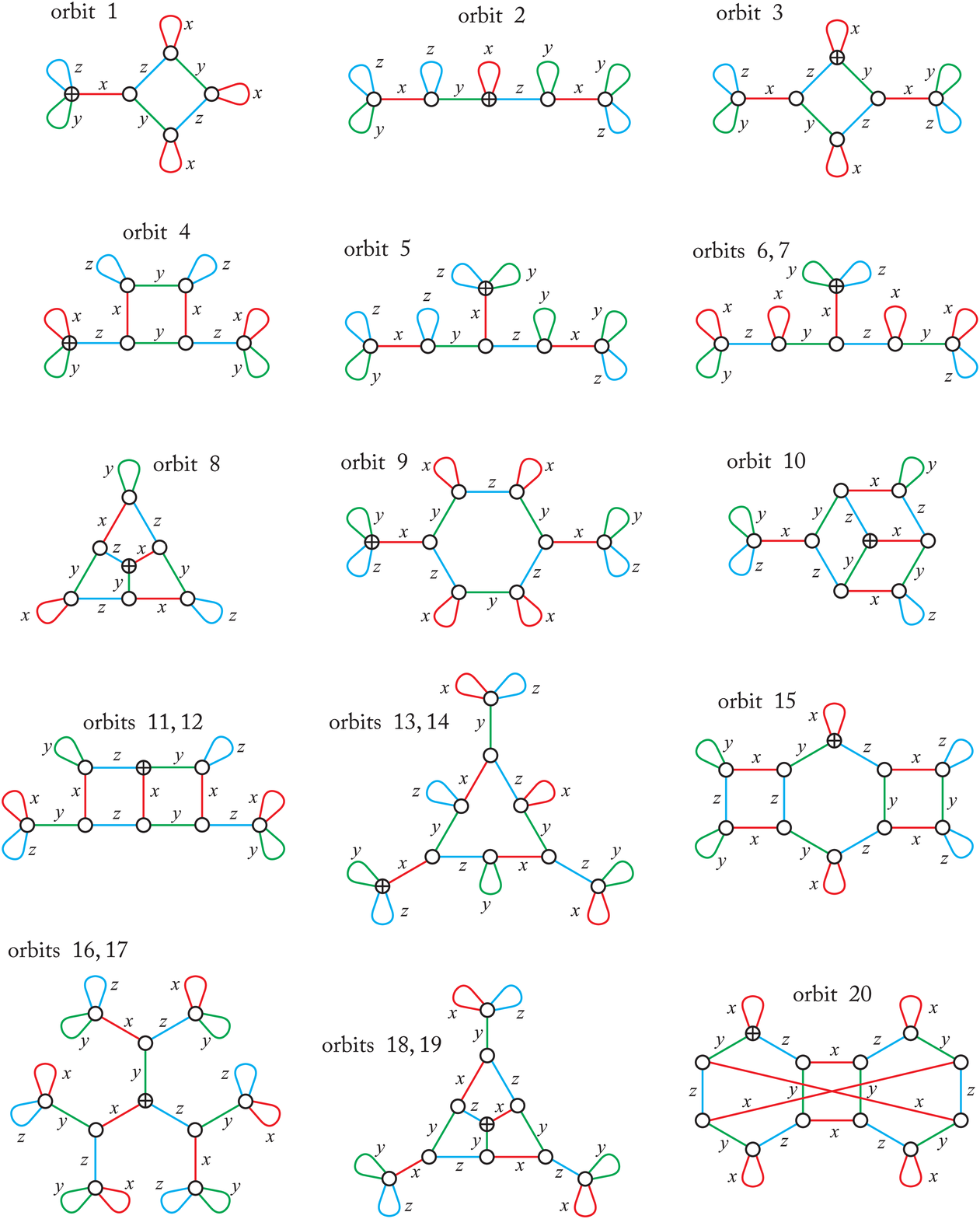}} \\  $\;$ \\ Fig. 9: Graphs of exceptional orbits 1--20 \\
 \end{center}
 \end{figure}
 \pagebreak[4]
  \begin{figure}[!h]
 \begin{center}
 $\;$\\ $\;$\\
 \resizebox{15cm}{!}{
 \includegraphics{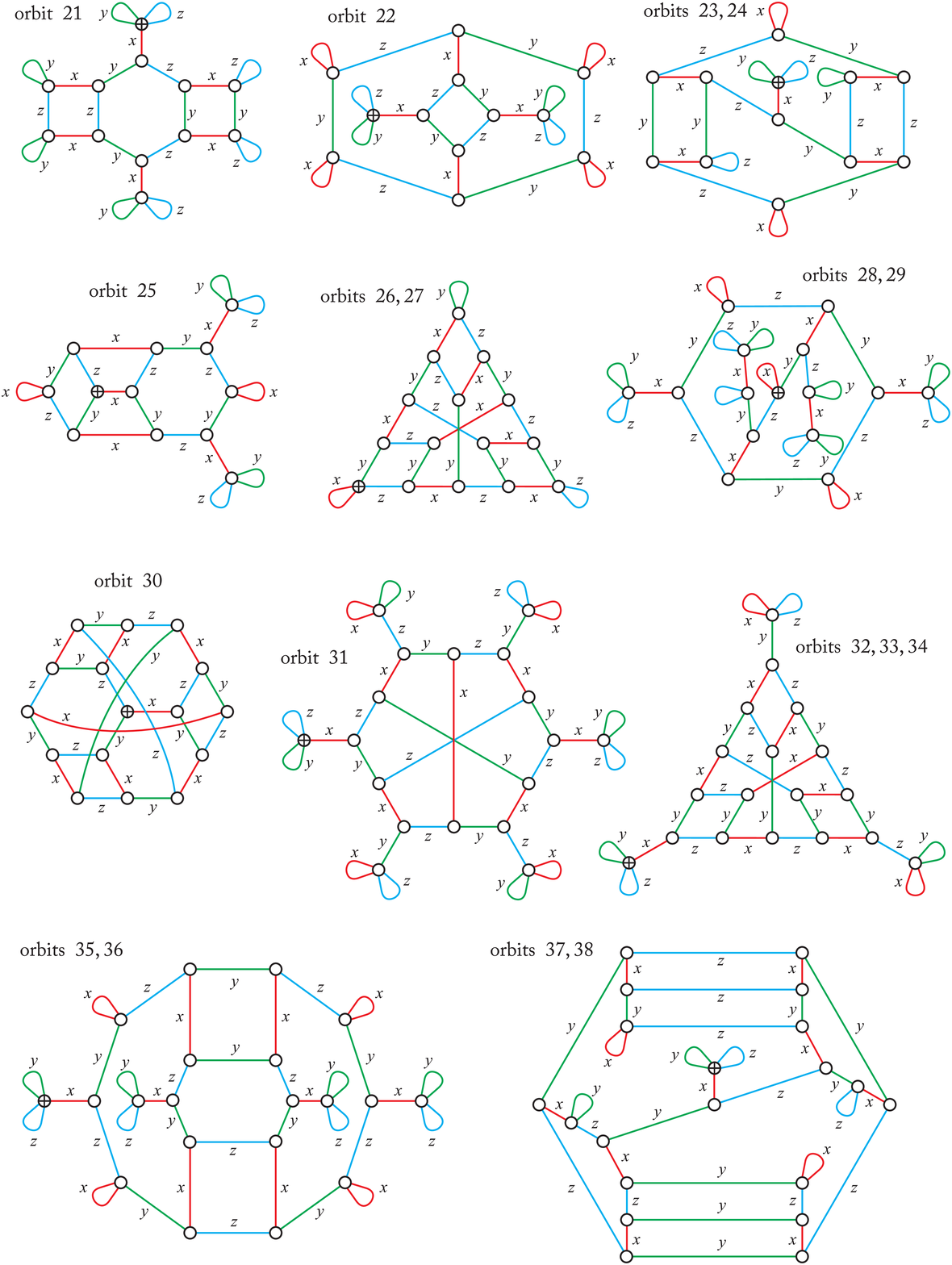}} \\  $\;$ \\ Fig. 10: Graphs of exceptional orbits 21--38 \\
 \end{center}
 \end{figure}
 \newpage
  \begin{figure}[!h]
 \begin{center}
 $\;$\\ $\;$\\
 \resizebox{15cm}{!}{
 \includegraphics{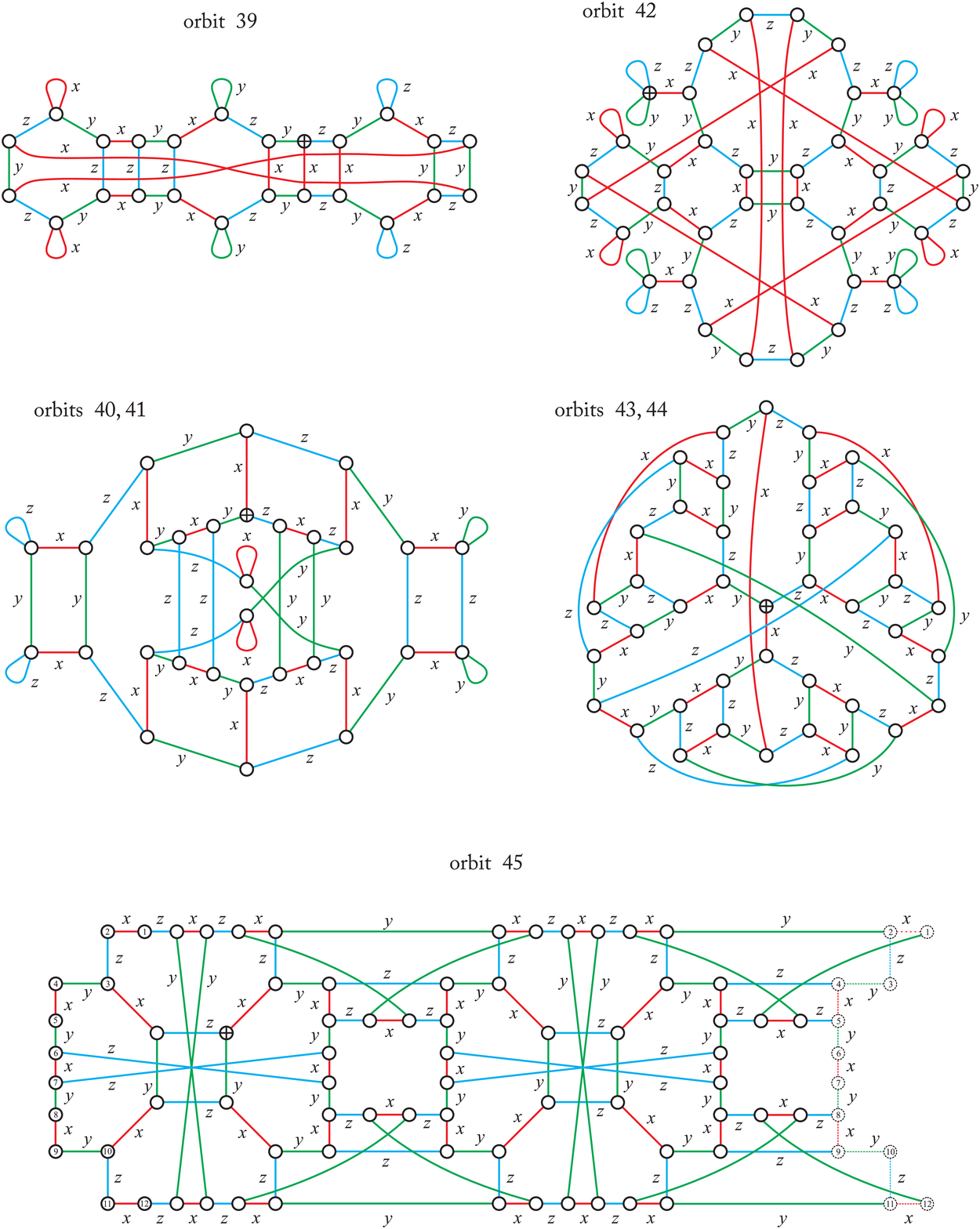}} \\  $\;$ \\ Fig. 11: Graphs of exceptional orbits 39--45 \\
 \end{center}
 \end{figure}
 \newpage

 We now turn to the description of nonequivalent finite orbits of the $\bar{\Lambda}$ action (\ref{lamb_matr}) on $\mathcal{M}$.
 Note that, given $\omega_{X,Y,Z,4}$, the equations (\ref{omxyz})--(\ref{omi})
 have only a finite number of solutions for $\{p_x,p_y,p_z,p_{\infty}\}$. In fact this number cannot exceed 24,
 see proof of Proposition~\ref{inj}, and all such solutions are related by the affine $D_4$ transformations.
 A natural question is therefore: when does the 7-tuple $\mathbf{p}=(p_x,p_y,p_z,p_\infty,X,Y,Z)$
 (see (\ref{pxyzi}), (\ref{xyz})) completely fix the conjugacy class of the triple $(M_x,M_y,M_z)\in G^3$, $G=SL(2,\Cb)$
 in $\mathcal{M}=G^3/G$~?

 Let us first prove an auxiliary result:
 \begin{lem}\label{mabc}
 Let $M^a,M^b,M^c\in G$ be three matrices such that the eigenvalues of at least
 one of them are different from $\pm1$. Then one and only one of the following holds:
 \begin{itemize}
 \item[1.] seven quantities
 \begin{align}\label{tabc1}
 t_a=\mathrm{Tr}\,M^a,\qquad t_b=\mathrm{Tr}\,M^b,\qquad t_c=\mathrm{Tr}\,M^c,\qquad
 t_{abc}=\mathrm{Tr}\left(M^a M^b M^c\right),
 \\ \label{tabc2}
  t_{ab}=\mathrm{Tr}\left(M^a M^b \right),\qquad  t_{ac}=\mathrm{Tr}\left(M^a M^c\right),
  \qquad  t_{bc}=\mathrm{Tr}\left(M^b M^c\right),\qquad
 \end{align}
 completely fix the conjugacy class of the triple $(M^a,M^b,M^c)$ in $\mathcal{M}$;
 \item[2.] $M^a,M^b,M^c$ have a common eigenvector.
 \end{itemize}
 \end{lem}
 \pf Using the same tricks as in the proof of Lemma~\ref{lemxyz}, one
 easily expresses $t_{bac}=\mathrm{Tr}\left(M^b M^a M^c\right)$ in terms of  (\ref{tabc1})--(\ref{tabc2}):
 \ben
 t_{bac}=\mathrm{Tr}\left(\left[t_{ab}\mathbf{1}-(M^a)^{-1}(M^b)^{-1}\right] M^c\right)
 =t_{ab}t_c-\mathrm{Tr}\left((t_a\mathbf{1}-M^a)(t_b\mathbf{1}-M^b)M^c\right)=
 \ebn
 \ben
 =t_{ab}t_c+t_{ac}t_b+t_{bc}t_a-t_a t_b \,t_c-t_{abc}.
 \ebn
 We may therefore assume without loss of generality that the eigenvalues
 of $M^a$ are not equal to $\pm1$;
 in particular, $M^a$ is diagonalizable. It is convenient to transform it into  diagonal form
 $M^a=\mathrm{diag}(\lambda_a , \lambda_a^{-1} )$, with $\lambda_a$ fixed by $t_a$.
 Now the equations for $t_b$ and $t_{ab}$ ($t_c$ and $t_{ac}$) fix ${M}^b_{11}$, ${M}^b_{22}$ and
 $M^b_{12}M^b_{21}$ (resp. ${M}^c_{11}$, ${M}^c_{22}$ and
 $M^c_{12}M^c_{21}$). The equations for $t_{bc}$ and $t_{abc}$, in their turn,
 completely determine $(M^bM^c)_{11}$ and $(M^bM^c)_{22}$, hence the products $M^b_{12}M^c_{21}$ and $M^b_{21}M^c_{12}$
 are also fixed.

 If $M^b_{12}M^b_{21}=M^c_{12}M^c_{21}=M^b_{12}M^c_{21}=M^b_{21}M^c_{12}=0$, then either
 $M^b_{12}=M^c_{12}=0$ or $M^b_{21}=M^c_{21}=0$, i.e. $M^{a,b,c}$ are simultaneously lower or
 upper triangular. On the other hand if at least one of the four products, say $M^b_{12}M^b_{21}$,
 is non-zero, then, using the remaining freedom of conjugation of $M^{a,b,c}$ by any diagonal matrix, one
 can set $M^{b}_{12}=1$ and $M^{b}_{21}(\neq0)$, $M^c_{12}$ and $M^c_{21}$ become completely fixed.
 Moreover, in this case $M^{a,b,c}$ clearly cannot have a common eigenvector.
 \epf
 \begin{lem}\label{fixed}
 Let $M_x,M_y,M_z\in G$. One and only one of the following holds:
 \begin{itemize}
 \item[1.]
 Conjugacy class of the triple $(M_x,M_y,M_z)$ in $\mathcal{M}$ is uniquely fixed by the
 7-tuple $\mathbf{p}=(p_x,p_y,p_z,p_{\infty},X,Y,Z)$, defined by (\ref{pxyzi})--(\ref{xyz}).
 \item[2.] $M_x,M_y,M_z$ have a common eigenvector.
 \end{itemize}
 \end{lem}
 \pf When the eigenvalues of at least one of three matrices $M_{x,y,z}$ are not equal to $\pm1$,
 the statement is equivalent to the previous lemma.

 Similarly, if e.g. the eigenvalues of $M_xM_y$ (or $M_xM_y^{-1}$)
 are different from $\pm1$, we can apply Lemma~\ref{mabc} to the triple $M^a=M_xM_y$,
 $M^b=M_y^{-1}$, $M^c=M_z$ (resp. $M^a=M_xM_y^{-1}$,
 $M^b=M_y$, $M^c=M_z$). Since $t_a$, $t_b$, $t_c$, $t_{ab}$, $t_{ac}$, $t_{bc}$, $t_{abc}$ are clearly
 expressible in terms of $\mathbf{p}$, the conjugacy class of $(M^a,M^b,M^c)$,
 and hence that of $(M_x,M_y,M_z)$, is fixed unless $M^{a,b,c}$
 can be simultaneously brought to lower or upper triangular form.

 Therefore, it is sufficient to prove the Lemma in the case when the eigenvalues
 of $M_{x,y,z}$, $M_xM_y^{\pm1}$, $M_xM_z^{\pm1}$ and $M_yM_z^{\pm1}$ are equal to $\pm1$.
 We can assume without loss of generality that
 $\mathrm{Tr}\,M_x=\mathrm{Tr}\,M_y=\mathrm{Tr}\,M_z=2$, but then
 from the relation $\mathrm{Tr}\left(M_xM_y\right)+\mathrm{Tr}\left(M_xM_y^{-1}\right)=
 \mathrm{Tr}\,M_x\cdot \mathrm{Tr}\,M_y$
 follows that $\mathrm{Tr}\left(M_xM_y\right)=2$.
 Similarly, one has $\mathrm{Tr}\left(M_xM_z\right)=\mathrm{Tr}\left(M_yM_z\right)=2$. Now, if we transform
 $M_x$ into upper triangular form, the relations $\mathrm{Tr}\,M_x=\mathrm{Tr}\,M_y=\mathrm{Tr}\left(M_xM_y\right)=2$
 imply that either $M_x$ is the identity matrix or $M_y$ is also upper triangular.
 Combining with analogous result for $M_x$, $M_z$ we see that all three matrices should have a common eigenvector.\epf
 \begin{lem}
 If three matrices $M_x,M_y,M_z\in G$ have a common eigenvector, then the elements of
 $\mathbf{p}$ satisfy characteristic relations (\ref{tetrahr}) of orbit~I, with $\omega_{X,Y,Z}$
 defined by (\ref{omxyz}).
 \end{lem}
 \pf Transforming $M_{x,y,z}$
 into upper triangular form, we see that $\mathbf{p}$ can be written in terms of the eigenvalues
 of $M_{x,y,z}$. It is sufficient to substitute these expressions into the relations (\ref{tetrahr})
 to check that they are satisfied automatically.\epf

 We now formulate a converse statement:
 \begin{lem}
 Let $M_x,M_y,M_z\in G$ be three matrices with no common eigenvector. If $\mathbf{p}$
  satisfies the relations (\ref{tetrahr}), then
  at least one of four matrices $M_x$, $M_y$, $M_z$,
  $M_z M_y M_x$ is equal to $\pm\mathbf{1}$.
 \end{lem}
 \pf Using (\ref{tetrahr}) and (\ref{jfr}), write $\omi$ in terms of $X,Y,Z$:
 \ben
 \omi=4+2XYZ+X^2+Y^2+Z^2.
 \ebn
 Substituting the expressions for $\omega_{X,Y,Z,4}$ into the cubic equation (\ref{cubicxi})
 for $\xi=p_x^2+p_y^2+p_z^2+p_{\infty}^2$,
 one finds that it has only two solutions:
 (1) $\xi=8+XYZ$ and (2) $\xi=4+X^2+Y^2+Z^2$.\vspace{0.1cm}\\
 \underline{Case (1)}. Let us write $X=2\cos\pi r_X$, $Y=2\cos\pi r_Y$, $Z=2\cos\pi r_Z$.
 It is straightforward to check that $(p_x^0,p_y^0,p_z^0,p_{\infty}^0)$ defined by
 \ben
 \begin{array}{ll}
 p_x^0=2\cos\pi(r_Y+r_Z-r_X)/2,\quad & p_y^0=2\cos\pi(r_X+r_Z-r_Y)/2,\vspace{0.1cm}\\
 p_z^0=2\cos\pi(r_X+r_Y-r_Z)/2,\quad & p_{\infty}^0=2\cos\pi(r_X+r_Z+r_Y)/2,
 \end{array}
 \ebn
 is one of possible solutions for $(p_x,p_y,p_z,p_{\infty})$. All other solutions characterized
 by the same value of $\xi$ have the form (\ref{pxyzid4}), see proof of Proposition~\ref{inj}. However,
 it is not difficult to show that for all such $(p_x,p_y,p_z,p_{\infty})$ one can find infinitely
 many triples $(M_x',M_y',M_z')$ of upper triangular matrices
 with the same $\mathbf{p}$ as $(M_x,M_y,M_z)$. E.g. if $p_{\nu}=p_{\nu}^0$,
 $\nu=x,y,z,\infty$, then we may set
 \begin{align*}
 &M_x'=\left(
 \begin{array}{cc}e^{i\pi(r_Y+r_Z-r_X)/2} & * \\ 0 & e^{-i\pi(r_Y+r_Z-r_X)/2} \end{array}\right),\\
 &M_y'=\left(
 \begin{array}{cc}e^{i\pi(r_X+r_Z-r_Y)/2} & * \\ 0 & e^{-i\pi(r_X+r_Z-r_Y)/2} \end{array}\right),\\
 &M_z'=\left(
 \begin{array}{cc}e^{i\pi(r_X+r_Y-r_Z)/2} & * \\ 0 & e^{-i\pi(r_X+r_Y-r_Z)/2} \end{array}\right).
 \end{align*}
 Now since $\mathbf{p}$ does not fix the conjugacy class of the triple $(M_x,M_y,M_z)$ uniquely,
 by Lemma~\ref{fixed} $M_{x,y,z}$ should have a common eigenvector.\vspace{0.1cm}\\
 \underline{Case (2)}. Here, one possible solution for $(p_x,p_y,p_z,p_{\infty})$ is
 \be\label{1psolution}
 p_x^0=X,\qquad p_y^0=Y,\qquad p_z^0=Z,\qquad p_{\infty}^0=2,
 \eb
 and all the others are given by (\ref{pxyzid4}).
 Consider the solution (\ref{1psolution}) and transform $M_z M_y M_x$ into upper triangular form:
 $M_zM_yM_x=\left(\begin{array}{cc}1 & \alpha \\ 0 & 1 \end{array}\right)$. Since
 \ben
 X=\mathrm{Tr}\left(M_y M_z\right)=
 \mathrm{Tr}\left(M_zM_yM_x\cdot M_x^{-1}\right)=p_x-\alpha \left(M_x\right)_{21},
 \ebn
 the relation $p_x=X$ implies that either $M_zM_yM_x=\mathbf{1}$ or $M_x$ is  upper triangular.
 Repeating  the same procedure with $p_y=Y$, $p_z=Z$ and using the assumption that $M_{x,y,z}$
 have no common eigenvectors, one concludes that $M_zM_yM_x=\mathbf{1}$. Other solutions
 for $(p_x,p_y,p_z,p_{\infty})$ are treated in a similar manner.\epf

 We thus obtain a
 description of all nonequivalent finite orbits of the $\bar{\Lambda}$ action (\ref{lamb_matr}) on $\mathcal{M}$:
 \begin{itemize}
 \item There are two families of nonequivalent orbits that consist of one point.
 They are given by the conjugacy
 classes of triples  (a) $(\mathbf{1},M_y,M_z)$ and (b) $(M_x,M_y,M_x^{-1}M_y^{-1})$,
 where $M_y$, $M_z$ in (a)  and $M_x$, $M_y$ in (b) have no common eigenvectors, $M_{x,y,z}\in G$.
 \item Each finite orbit $O$
 of the induced $\bar{\Lambda}$ action (\ref{lxyz}) that consists of more than one
  point (i.e. each of orbits~II--IV, 1--45 and Cayley orbits of size greater than one)
  generates a finite number of orbits
  of (\ref{lamb_matr}), which have the same size as $O$ and correspond to different 4-tuples $(p_x,p_y,p_z,p_{\infty})$
  solving (\ref{omxyz})--(\ref{omi}).
  (Recall that the parameters $\omega_{X,Y,Z,4}$ for orbits~II--IV and 1--45 are specified by Lemma~\ref{orb14} and Table~4, while for Cayley orbits $\omx=\omy=\omz=\omi=0$.) Once a solution
  for  $(p_x,p_y,p_z,p_{\infty})$ is chosen,  the orbit in $\mathcal{M}$ is completely fixed by the 7-tuple
  $(p_x,p_y,p_z,p_{\infty},X,Y,Z)$, where $(X,Y,Z)$ is any point in~$O$.
  \item All remaining finite orbits of (\ref{lamb_matr}) belong to the space $\mathcal{U}\subset\mathcal{M}$
  of conjugacy classes of triples of upper triangular $SL(2,\Cb)$-matrices.
 \end{itemize}

 \section{Algebraic Painlev\'e VI solutions}
 We are now prepared for the classification of PVI solutions with finite branching up to parameter
 equivalence.
 \begin{defn}
 Let us associate to any PVI solution branch the 7-tuple of monodromy data
 $(\omx,\omy,\omz,\omi,X,Y,Z)\in \Cb^7$  defined by (\ref{pxyzi})--(\ref{omi}). Two finite branch PVI solutions
  will be called
  \begin{itemize}
  \item
  \textit{equivalent}
 if they are related by B\"acklund transformations specified in Table~1;
 \item \textit{parameter equivalent}
 if their analytic continuation leads to equivalent (under $K_4\rtimes S_3$ transformations of Subsection~2.2) orbits in the space of  7-tuples of monodromy data.
 \end{itemize}
 \end{defn}
 \begin{rmk}
 Our parameter equivalence is strictly stronger than that of \cite{boalch52}, and is rather similar to geometric
 equivalence, cf. \cite{boalch52}, Def.~8. In particular, it distinguishes solutions~3, 21 and~42 (see below), whose parameters $\boldsymbol{\theta}=(\theta_x,\theta_y,\theta_z,\theta_{\infty})$ lie in the same orbit of the Okamoto affine~$F_4$ action. Another such example is given by solutions~20 and~45.
 \end{rmk}
 \begin{rmk} In \cite{boalch52}, p.~13 it is stated that the four-branch octahedral PVI solution \cite{hitchin2}
 \be\label{4oct}
 w=\frac{(s-1)^2}{s(s-2)},\qquad t=\frac{(s+1)(s-1)^3}{s^3 (s-2)},
 \eb
 with parameters $\boldsymbol{\theta}=(\vartheta,\vartheta,\vartheta,1-3\vartheta)$ and the four-branch dihedral solution~IV
 below are inequivalent for $\vartheta=\theta=1/6$, although characterized by the same parameters. This seems to be
 incorrect; replacing $s\mapsto{1}/({s+1})$ in (\ref{4oct}) and applying to $w$ affine~$D_4$
 transformation $s_x s_y s_z s_{\delta} s_x s_y s_z$, one finds solution~IV with $\theta=1/2-2\vartheta$.
 Despite the failure of the above counterexample, our parameter equivalence is presumably weaker than the equivalence under B\"acklund transformations.
 \end{rmk}

 Let us now examine one by one all finite orbits listed in Theorem~\ref{theorem1}  (recall that finite orbits which are not
 of Cayley type do not split under the action of $\Lambda$). First consider orbit~I, consisting of a single
 point. In this case all solutions of Painlev\'e~VI can be found explicitly. In particular, for reducible monodromy
 (i.e.  when $M_x,M_y,M_z$ have a common eigenvector)  PVI equation  linearizes
 and one has the following:
 \begin{prop}[Theorem 4.1 in \cite{mazzocco}]
 All solutions of PVI corresponding to reducible monodromy are equivalent to
 the one-parameter family of Riccati solutions
 \be\label{1pfamily}
   w(t)=\frac{(1+\theta_x+\theta_z-t-\theta_z t)u(t)-t(t-1)u'(t)}{(1+\theta_x+\theta_y+\theta_z)u(t)},
 \eb
 realized for $\theta_{\infty}=-(\theta_x+\theta_y+\theta_z)$, where $u(t)=u_1(t)+\nu u_2(t)$ and $u_{1,2}(t)$
 are two linear independent solutions of the following hypergeometric equation:
 \be\label{HGequation}
 t(1-t)u''+[(2+\theta_x+\theta_z)-(4+\theta_x+\theta_y+2\theta_z)t]u'-
 (2+\theta_x+\theta_y+\theta_z)(\theta_z+1)u=0.
 \eb
 \end{prop}
 \begin{rmk} It is well-known that one-parameter family (\ref{1pfamily}) contains solutions
 with a finite number of branches if and only if the parameters of the hypergeometric equation
 (\ref{HGequation}) belong to the Schwarz table, see \cite{schwarz} or Table~1 in \cite{boalch_SL}.
 \end{rmk}
 The solutions of PVI in the case of ``1-smaller monodromy'',
 when one of the matrices $M_x$, $M_y$, $M_z$ or $M_{\infty}=(M_zM_yM_x)^{-1}$ is equal to $\pm\mathbf{1}$,
 have been completely described in
 \cite{mazzocco_garnier}. Any such solution is either i) degenerate ($w=0,1,t,\infty$)
 or ii) equivalent via B\"acklund transformations to a Riccati solution
 or iii) belongs to a set of generalized Chazy solutions, expressible in terms of hypergeometric functions; see Lemma~33 in \cite{mazzocco_garnier}
 for the details.

 Next we consider Cayley orbits. Since in this case $\omx=\omy=\omz=\omi=0$, the 4-tuple $(p_x,p_y,p_z,p_{\infty})$
 can only be $(0,0,0,0)$ or a permutation of $(\pm2,\pm2,\pm2,\mp2)$. This in turn implies that the 4-tuple
 of PVI parameters
 $(\theta_x,\theta_y,\theta_z,\theta_{\infty})$ consists of either i)~1~odd and 3~even integers or ii)~1~even and 3~odd integers or iii) all four $\theta_{x,y,z,\infty}$ have half-integer values. For
 $\theta_x=\theta_y=\theta_z=0$, $\theta_{\infty}=1$ the general solution of Painlev\'e~VI is known:
 \begin{prop}
 All solutions of the sixth Painlev\'e equation with $\theta_x=\theta_y=\theta_z=0$, $\theta_{\infty}=1$
 are given by Picard solutions
 \be\label{picard}
 w(t)=\wp \left(\nu_1 u_1+\nu_2 u_2;u_1,u_2\right)+\frac{t+1}{3},\qquad\qquad \nu_{1,2}\in\Cb,\quad
 0\leq\mathrm{Re}\,\nu_{1,2}<2,
 \eb
 where $\wp(z;u_1,u_2)$ is the Weierstrass elliptic function and
  $u_{1,2}(t)$ are two linearly independent solutions of the following hypergeometric equation:
 \be\label{hgeq2}
 4t(1-t)u''+4(1-2t)u'-u=0.
 \eb
 \end{prop}
 \pf This statement was proved by Fuchs in \cite{fuchs}.\epf\\
 All finite branch solutions corresponding to Cayley orbits are therefore parameter equivalent
 to solutions from the above two-parameter family. Equivalence under B\"acklund transformations
 is slightly more subtle, see e.g.~\cite{mazzocco2}.

 There remain precisely 45 parameter inequivalent finite branch PVI solutions and three families depending
 on continuous parameters, which correspond to orbits~1--45 and II--IV
 (existence of solutions with appropriate monodromy data
 follows from their explicit construction below).
 Surprisingly, each equivalence class contains algebraic representatives that
 have already appeared in the literature \cite{andreev_kitaev,boalch_klein,boalch52,boalch_RH,boalch_HG,
  dubrovin2,dubrovin,hitchin,hitchin2,kitaev,kitaev2}.
 Complete list of these parameter inequivalent algebraic solutions is given below.
 For each solution we specify the 4-tuple of
 PVI parameters $\boldsymbol{\theta}=(\theta_x,\theta_y,\theta_z,\theta_{\infty})$, the number
 of branches and the explicit solution curve. We also give references to original papers
 where the corresponding algebraic solutions have been obtained and correct a few  misprints
 (in solutions~13, 24, 43 and 44).\vspace{0.2cm}\\
 \textit{Solution~II}, 2 branches, $\boldsymbol{\theta}=(\theta_a,\theta_b,\theta_b,1-\theta_a)$:
 \ben
 w(t)=\pm\sqrt{t}.
 \ebn
 In Lemma~\ref{orb14}, $X'=2\cos2\pi\theta_b$, $X''=-2\cos2\pi\theta_a$.\\
 \textit{Solution~III}, 3 branches, $\boldsymbol{\theta}=(2\theta,\theta,\theta,2/3)$:
 \ben
 w=\frac{(s-1)(s+2)}{s(s+1)},\qquad t=\frac{(s-1)^2(s+2)}{(s+1)^2(s-2)},
 \ebn
 first obtained in \cite{dubrovin2}, (E.31);
 in the above form it appeared in \cite{hitchin2}. In Lemma~\ref{orb14}, $\omega=2\cos3\pi\theta$.\\
 \textit{Solution~IV}, 4 branches, $\boldsymbol{\theta}=({\theta},{\theta},{\theta},1/2)$:
 \ben
 w=\frac{s^2(s+2)}{
s^2+s+1},\qquad t=\frac{s^3(s+2)}{2s+1},
 \ebn
 first obtained in \cite{dubrovin2}, (E.29);
 in the above form it appeared in \cite{hitchin}. In Lemma~\ref{orb14}, $\omega=4\cos^2\pi{\theta}$.
 \\
 \textit{Solution~1}, 5 branches, $\boldsymbol{\theta}=(2/5,1/5,1/3,2/3)$:
 \ben
 w=\frac{2(s^2+s+7)(5s-2)}{s(s+5)(4s^2-5s+10)}, \qquad t=\frac{27(5s-2)^2}{(s+5)(4s^2-5s+10)^2},
 \ebn
 solution~20 in \cite{boalch52}, p.~21.\\
 \textit{Solution~2}, 5 branches, $\boldsymbol{\theta}=(1/5,2/5,1/5,2/5)$:
 \ben
 w=\frac{s^2(s-1)}{3(s-2)(s+3)}, \qquad t=\frac{2s^3(s^2-5)}{(s-2)^2(s+3)^3},
 \ebn
 first found in \cite{kitaev}, Eq. (3.3).\\
 \textit{Solution 3}, 6 branches, $\boldsymbol{\theta}=(1/2,1/3,1/3,1/2)$:
 \ben
 w=-
\frac{s (s + 1) (s - 3)^2}{
3 (s + 3) (s - 1)^2},\qquad t=-
\frac{(s + 1)^3 (s - 3)^3}{
(s - 1)^3 (s + 3)^3},
 \ebn
 first found in \cite{andreev_kitaev}, equivalent to solution~4.1.1A; in the above form in
 \cite{boalch_RH}, tetrahedral solution~6, p.~9.\\
 \textit{Solution 4}, 6 branches, $\boldsymbol{\theta}=(1/2,1/4,1/2,2/3)$:
 \ben
 w=\frac{9s (2 s^3 - 3 s + 4)}{
4 (s + 1) (s - 1)^2 (2 s^2 + 6 s + 1)},\qquad t=\frac{27s^2}{
4 (s^2 - 1)^3},
 \ebn
 octahedral solution~7 in \cite{boalch_RH}, p.~12.\\
 \textit{Solution 5}, 6 branches, $\boldsymbol{\theta}=(1/4,1/4,1/3,1/3)$:
 \ben
 w=\frac{(3s - 1)(2s - 1)(s + 1)^3}{
4s(3s^2 - 1)(s^2 + 1)},\qquad t=\frac{(s + 1)^4(2s - 1)^2}{
8s^3(3s^2 - 1)},
 \ebn
 first found in \cite{kitaev} 3.3.3, p.~22.\\
 \textit{Solution 6}, 6 branches, $\boldsymbol{\theta}=(2/5,1/5,2/5,2/3)$:
 \ben
 w=\frac{18s(s-3)}{(s-4)(s+1)(s^2+5)},\qquad t=\frac{432s}{(s+5)(s+1)^3(s-4)^2},
 \ebn
 solution~23 in \cite{boalch52}, p.~23. \\
  \textit{Solution 7}, 6 branches, $\boldsymbol{\theta}=(1/5,2/5,1/5,1/3)$:
 \ben
 w=\frac{-54s(s-7)}{(s-4)(s+1)(s^4-20s^2-35)},\qquad t=t_6,
 \ebn
 solution~22 in \cite{boalch52}, p.~23. \\
 \textit{Solution 8}, 7 branches, $\boldsymbol{\theta}=(2/7,2/7,2/7,4/7)$:
 \ben
 w=-\frac{(5s^2-8s+5)(7s^2-7s+4)}{s(s-2)(s+1)(2s-1)(4s^2-7s+7)},\qquad t=\frac{(7s^2-7s+4)^2}{s^3(4s^2-7s+7)^2}.
 \ebn
 Klein solution of \cite{boalch_klein}, p.~26.\\
 \textit{Solution 9}, 8 branches, $\boldsymbol{\theta}=(1/4,1/2,1/4,1/2)$:
 \ben
 w=-
\frac{(s^2 - 2s + 2)(s^2 + 2)^2}{
4(s + 1)(s^2 - 4s - 2)(s - 1)^2},\qquad t=\frac{(s^2 - 2)(s^2 + 2)^3}{
16(s + 1)^3(s - 1)^3},
 \ebn
 first found in \cite{kitaev} 3.3.5, p.~23.\\
 \textit{Solution 10}, 8 branches, $\boldsymbol{\theta}=(1/3,1/2,1/4,2/3)$:
 \ben
 w=\frac{s^3 (2 s^2 - 4 s + 3) (s^2 - 2 s + 2)}{
(2 s^2 - 2 s + 1) (3 s^2 - 4 s + 2)},\qquad t=\left(\frac{s^2 (2 s^2 - 4 s + 3)}{
3 s^2 - 4 s + 2}\right)^2,
 \ebn
 octahedral solution~9 in \cite{boalch_RH}, p. 12.\\
   \textit{Solution 11}, 8 branches, $\boldsymbol{\theta}=(1/2,1/5,2/5,4/5)$:
 \ben
 w=\frac{s(s+4)(3s^4-2s^3-2s^2+8s+8)}{8(s-1)(s+1)^2(s^2+4)},\qquad t=\frac{s^5(s+4)^3}{4(s-1)(s+1)^3(s^2+4)^2},
 \ebn
 solution~24 in \cite{boalch52}, p.~21. \\
    \textit{Solution 12}, 8 branches, $\boldsymbol{\theta}=(2/5,1/2,2/5,4/5)$:
 \ben
 w=\frac{s^2(s+4)^2(5s^3+2s^2-4s-8)}{4(s-1)(s+1)^2(s^2+4)(s^2+3s+6)},\qquad t=t_{11},
 \ebn
 solution~25 in \cite{boalch52}, p.~21. \\
 \textit{Solution 13}, 9 branches, $\boldsymbol{\theta}=(2/5,2/5,2/5,2/3)$:
 \begin{align*}
 &w=\frac{1}{2}+\frac{350s^3+63s^2-6s-2}{30s(2s+1)u},\\
 &t=\frac{1}{2}+\frac{(25s^4+170s^3+42s^2+8s-2)u}{54s^3(5s+4)^2},\\
 &u^2=s(8s+1)(5s+4),
 \end{align*}
 solution~27 in \cite{boalch52}, p.~23 (parameters in \cite{boalch52} are defined with a misprint, which
 is corrected by interchanging $\theta_3\leftrightarrow\theta_4$).\\
 \textit{Solution 14}, 9 branches, $\boldsymbol{\theta}=(1/5,1/5,1/5,1/3)$:
 \begin{align*}
 &w=\frac{1}{2}-\frac{(s - 1)\bigl(5(s^6 + 1) + 58(s^5 + s) + 1771(s^4 + s^2) + 8620s^3\bigr)u}{8s(s+1)(5s^3 + 25s^2 + 95s + 3)(3s^3 + 95s^2 + 25s + 5)},\\
 &t=\frac{1}{2}-\frac{(s - 1)\bigl(25(s^8 + 1) + 760(s^7 + s) + 4924(s^6 + s^2) + 75464(s^5 + s^3) + 329174s^4)}{
 2048s(s + 1)^5 u},\\
 &u^2=s(5s^2+118s+5),
 \end{align*}
 first found in \cite{kitaev}, p.~11.\\
 \textit{Solution 15}, 10 branches, $\boldsymbol{\theta}=(1/2,1/5,1/2,3/5)$:
 \begin{align*}
 &
 w=\frac{(s^2-5)(s^2+5)(s^5+5s^4-20s^3+75s+75)}{(s+1)^2 (s+5)(s^2-4s+5)(s^4+6s^2-75)},\\
 &t=\frac{2(s^2+5)^3(s^2-5)^2}{(s+1)^3(s+5)^3(s^2-4s+5)^2},
 \end{align*}
 solution~28 in \cite{boalch52}, p.~21.\\
 \textit{Solution 16}, 10 branches, $\boldsymbol{\theta}=(0,0,0,-4/5)$:
 \begin{align*}
 &w=\frac{(s-1)^2(3s+1)^2(s^2+4s-1)(119s^8-588s^6+314s^4-108s^2+7)^2}{(s+1)^3(3s-1)P(s)},\\
 &t=\frac{(s-1)^5(3s+1)^3(s^2+4s-1)}{(s+1)^5(3s-1)^3(s^2-4s-1)},\\
 &P(s)=42483s^{18}-719271s^{16}+5963724s^{14}+13758708s^{12}-7616646s^{10}\\
 &\qquad\quad+1642878s^8-259044s^6+34308s^4-2133s^2+49,
 \end{align*}
 first obtained in \cite{dubrovin2}, the above parametrization corresponds to
 icosahedron solution ($H_3$) in \cite{dubrovin}, p.~76.\\
 \textit{Solution 17}, 10 branches, $\boldsymbol{\theta}=(0,0,0,-2/5)$:
 \begin{align*}
 &w=\frac{(s-1)^4(3s+1)^2(s^2+4s-1)(11s^4-30s^2+3)^2}{(s+1)(3s-1)(3s^2+1) P(s)},\\
 &t=\frac{(s-1)^5(3s+1)^3(s^2+4s-1)}{(s+1)^5(3s-1)^3(s^2-4s-1)},\\
 &P(s)=121s^{12}-1942s^{10}+63015s^8-28852s^6+4855s^4-342s^2+9,
 \end{align*}
 great icosahedron solution ($H_3$)$'$ in \cite{dubrovin}, p.~77.\\
 \textit{Solution 18}, 10 branches, $\boldsymbol{\theta}=(1/3,1/3,1/3,4/5)$:
 \ben
 w=\frac{s^2(s+2)(s^2+1)(2s^2+3s+3)}{2(s^2+s+1)(3s^2+3s+2)},\qquad
 t=\frac{s^5(s+2)(2s^2+3s+3)^2}{(2s+1)(3s^2+3s+2)^2},
 \ebn
 solution~29 in \cite{boalch52}, p.~23.\\
  \textit{Solution 19}, 10 branches, $\boldsymbol{\theta}=(1/3,1/3,1/3,2/5)$:
 \ben
 w=\frac{s^4(s+2)(2s^2+3s+3)(7s^2+10s+7)}{(3s^2+3s+2)\bigl(4(s^6+1)+12(s^5+s)+15(s^4+s^2)+10s^3\bigr)},\qquad
 t=t_{18},
 \ebn
 solution~30 in \cite{boalch52}, p.~23.\\
 \textit{Solution 20}, 12 branches, $\boldsymbol{\theta}=(1/2,1/2,1/2,2/3)$:
 \begin{align*}
 &w=\frac{1}{2}+\frac{
45 s^6 + 20 s^5 + 95 s^4 + 92 s^3 + 39 s^2 - 3}{
4 (5 s^2 + 1) (s + 1)^2 u},\\
&t=\frac{1}{2}+
\frac{s (2 s + 1)^2 (27 s^4 + 28 s^3 + 26 s^2 + 12 s + 3)}{
(s + 1)^3 u^3},\\
&u^2 = (2 s + 1)
(9 s^2 + 2 s + 1),
 \end{align*}
 octahedral solution~12 in \cite{boalch_RH}, p.~13.\\
 \textit{Solution 21}, 12 branches, $\boldsymbol{\theta}=(1/3,1/2,1/2,2/3)$:
 \begin{align*}
 &w=\frac{4(s + 1)(3s^2-4s+2) (7 s^4 + 16 s^3 + 4 s^2 - 4) }{
s^3 (s - 2)(s^2+4s+6) (s^4 - 4 s^2 + 32 s - 28)},\\
&t=\frac{16(s+1)^4(3s^2-4s+2)^2}{s^4(s-2)^4(s^2+4s+6)^2},
 \end{align*}
 octahedral solution~11 in \cite{boalch_RH}, p.~12.\\
 \textit{Solution 22}, 12 branches, $\boldsymbol{\theta}=(1/3,1/3,1/5,2/5)$:
 \begin{align*}
 &w=\frac{1}{2}
+
\frac{140 s^6 + 1029 s^5 - 1023 s^4 + 360 s^3 - 288 s^2 + 27 s + 27}{
18 u (s + 1) (7 s^3 - 3 s^2 - s + 1)},\\
&t=\frac{1}{2}+\frac{
40 s^6 + 540 s^5 - 765 s^4 + 540 s^3 - 270 s^2 + 27}{
6 u (8 s^2 - 9 s + 3) (s + 1)^2},\\
&u^2 = 3(5s+1)(8s^2-9s+3),
 \end{align*}
 solution~36 in \cite{boalch52}, p.~22.\\
 \textit{Solution 23}, 12 branches, $\boldsymbol{\theta}=(1/5,1/5,1/3,1/2)$:
 \begin{align*}
 &w=\frac{1}{2}+
\frac{(3 s + 5) (8 s^4 - 10 s^3 + 12 s^2 - 13 s + 11)}{
2 (2 s^3 - 15 s + 5) u},\\
&t=\frac{1}{2} -
\frac{8 s^6 + 20 s^3 - 15 s^2 + 66 s - 15}{
2 (8 s^2 - 5 s + 5) u},\\
& u^2 = (3s + 5)(8s^2 - 5s + 5),
 \end{align*}
 solution~34 in \cite{boalch52}, p.~21.\\
  \textit{Solution 24}, 12 branches, $\boldsymbol{\theta}=(2/5,2/5,1/3,1/2)$:
 \begin{align*}
 &w=\frac{1}{2}-
\frac{(3 s + 5) (16 s^5 - 8 s^4 + 18 s^3 - 8 s^2 + 115 s + 3)}{
2 (26 s^3 + 60 s^2 + 15 s + 35) u},\\
&t=t_{23},\qquad u = u_{23},
 \end{align*}
 solution~35 in \cite{boalch52}, p.~22 (in \cite{boalch52} there is a sign misprint in the formula for $w$).\\
 \textit{Solution 25}, 12 branches, $\boldsymbol{\theta}=(2/5,1/3,1/2,4/5)$:
 \begin{align*}
 &w=-\frac{9s(s^2 + 1)(3s - 4)(15s^4 - 5s^3 + 3s^2 - 3s + 2)}{
 (2s - 1)^2(9s^2 + 4)(9s^2 + 3s + 10)},\\
 & t=\frac{27s^5(s^2 + 1)^2(3s - 4)^3}{
4(2s - 1)^3(9s^2 + 4)^2},
 \end{align*}
 solution~33 (generic icosahedral solution) in \cite{boalch52}, Th.~B, p.~4.\\
  \textit{Solution 26}, 15 branches, $\boldsymbol{\theta}=(1/3,1/3,1/3,3/5)$:
 \begin{align*}
 &w =\frac{1}{2}-
\frac{250 s^6 + 500 s^5 + 518 s^4 + 261 s^3 + 76 s^2 + 13 s + 2}{
2 (s + 2) (5 s + 1) (5 s^3 + 6 s^2 + 3 s + 1) u},\\
&t =\frac{1}{2}-
\frac{3 (500 s^7 + 925 s^6 + 1164 s^5 + 830 s^4 + 340 s^3 + 105 s^2 + 20 s + 4)}{
2 (s + 2)^2 (5 s+ 1) u^3},\\
& u^2 =
(4 s^2 + s + 1)
(5 s + 1),
 \end{align*}
 solution~38 in \cite{boalch52}, p.~26.\\
 \textit{Solution 27}, 15 branches, $\boldsymbol{\theta}=(1/3,1/3,1/3,1/5)$:
 \begin{align*}
 &w =\frac{1}{2}-
\frac{1000 s^8 + 2425 s^7 + 4171 s^6 + 3805 s^5 + 1999 s^4 + 874 s^3 + 244 s^2 + 58 s + 4}{
4 (s + 2) (25 s^6 + 135 s^5 + 111 s^4 + 91 s^3 + 36 s^2 + 6 s + 1) u},\\
 &t=t_{26},\qquad u=u_{26},
 \end{align*}
 solution~37 in \cite{boalch52}, p.~26.\\
 \textit{Solution 28}, 15 branches, $\boldsymbol{\theta}=(3/5,3/5,2/3,2/3)$:
 \begin{align*}
 &w=\frac{1}{2}-
\frac{2 s^9 + 20 s^8 + 53 s^7 - 89 s^6 - 605 s^5 - 851 s^4 - 1389 s^3 - 5775 s^2 - 10125 s - 5625}{
2 (s^2 - 5) (s^2 - 6 s - 15) (s^2 + 4 s + 5) u},\\
&t =\frac{1}{2}-\frac{
(2 s^7 + 10 s^6 - 90 s^4 - 135 s^3 + 297 s^2 + 945 s + 675) u}{
18 (4 s^2 + 15 s + 15)^2 (s^2 - 5)},\\
& u^2 = 3 (s + 5) (4 s^2 + 15 s + 15),
 \end{align*}
 solution~40 in \cite{boalch52}, p.~22.\\
 \textit{Solution 29}, 15 branches, $\boldsymbol{\theta}=(1/3,1/3,4/5,4/5)$:
 \begin{align*}
 &w=\frac{1}{2}+\frac{
14 s^5 + 61 s^4 - 66 s^3 - 660 s^2 - 900 s - 225}{
6 (s + 1) (s^2 - 5) u},\\
 &t=t_{28},\qquad u=u_{28},
 \end{align*}
 solution~39 in \cite{boalch52}, p.~22.\\
 \textit{Solution 30}, 16 branches, $\boldsymbol{\theta}=(1/2,1/2,1/2,3/4)$:
 \begin{align*}
 &w=-
\frac{(1 + i) (s^2 - 1) (s^2 + 2 is + 1) (s^2 - 2 is + 1)^2 P(s)}{
4 s (s^2 + i) (s^2 - i)^2 (s^2 + (1 + i)s - i) Q(s)},\\
&t=\frac{(s^2 - 1)^2 (s^4 + 6 s^2 + 1)^3}{
32 s^2 (s^4 + 1)^3},\\
&P(s)=s^8-(2-2i) s^7-(6+2i) s^6+(10+2i) s^5+4 is^4+(10-2i) s^3\\ &\qquad\quad +(6-2i) s^2-(2+2i) s-1,\\
&Q(s)=s^6 - (3 + 3i) s^5 + 3 is^4 + (4 - 4i) s^3 + 3 s^2 + (3 + 3i) s + i,
 \end{align*}
 octahedral solution~13 in \cite{boalch_RH}, p.~13.\\
 \textit{Solution 31}, 18 branches, $\boldsymbol{\theta}=(1/3,1/3,1/3,1/3)$:
 \begin{align*}
 &w=\frac{1}{2}-\frac{8s^7-28s^6+75s^5+31s^4-269s^3+318s^2-166s+56}{18u(s-1)(3s^3-4s^2+4s+2)},\\
 &t=\frac{1}{2}+\frac{(s+1)\bigl(32(s^8+1)-320(s^7+s)+1112(s^6+s^2)-2420(s^5+s^3)+3167s^4\bigr)}{54u^3s(s-1)},\\
 &u^2=s(8s^2-11s+8).
 \end{align*}
 A solution with equivalent parameters was first obtained
 in \cite{dubrovin} (great dodecahedron solution ($H_3$)$''$, see pp.~78--87 in the preprint
 version of \cite{dubrovin} for the explicit form), the above elliptic
 parameterization was produced in \cite{boalch52}, Th.~C, p.~4.\\
 \textit{Solution 32}, 18 branches, $\boldsymbol{\theta}=(4/7,4/7,4/7,1/3)$:
 \begin{align*}
 &w =\frac{1}{2}-
\frac{
P(s)u}{
Q(s)},\qquad
t =\frac{1}{2}-
\frac{R(s) u}{
432 s (s + 1)^2 (s^2 + s + 7)^2},\qquad u^2 = s (s^2 + s + 7),\\
 &P(s)=s^{10} + 5 s^9 + 24 s^8 + 20 s^7 - 266 s^6 - 2874 s^5 - 14812 s^4 \\
 &\qquad\quad  - 40316 s^3 - 85359 s^2 - 100067 s - 67396,\\
 &Q(s)=16 (s + 1) (s^2 + s + 7) (5 s^6 + 63 s^5 + 252 s^4 + 854 s^3 + 1449 s^2 + 1827 s + 2030),\\
 &R(s)=s^9 - 84 s^6 - 378 s^5 - 1512 s^4 - 5208 s^3 - 7236 s^2 - 8127 s - 784,
 \end{align*}
 first appeared in \cite{boalch_RH}, p.~22.\\
 \textit{Solution 33}, 18 branches, $\boldsymbol{\theta}=(1/3,1/7/,1/7,6/7)$:
 \begin{align*}
 &w=1 +
\frac{(3s - 2)(s^2 - 2s + 4)^2}{
4(s + 2)(s - 1)^2(s^2 - s + 1)(3s^2 - 4s + 4)}\times\\
&\qquad\quad \times\frac{-14s^5 + 25s^4 - 20s^3 - 8s^2 + 16s - 8 - 8(s - 1)(s^2 - s + 1)u}{(2s + 1)(3s^3 - 10s^2 + 6s - 2) - 14(s - 1)u},
 \\
 &t =\frac{1}{2} -
\frac{14s^9 - 105s^8 + 252s^7 - 392s^6 + 420s^5 - 336s^4 + 112s^3 + 72s^2 - 96s + 32}{
16(s + 2)^2(s - 1)^3(s^2 - s + 1)u},\\
& u^2=(2s + 1)(1 - s)(s^2 - s + 1),
 \end{align*}
 solution (3.16)--(3.17) in \cite{kitaev2}, p. 15.
 \\
 \textit{Solution 34}, 18 branches, $\boldsymbol{\theta}=(2/7,2/7,2/7,1/3)$:
 \begin{align*}
 &w=\frac{1}{2}-
\frac{(3 s^8 - 2 s^7 - 4 s^6 - 204 s^5 - 536 s^4 - 1738 s^3 - 5064 s^2 - 4808 s - 3199) u}{
4 (s + 1)(s^2 + s + 7)(s^6 + 196 s^3 + 189 s^2 + 756 s + 154)  },\\
 &t=t_{32},\qquad u=u_{32},
 \end{align*}
 first appeared in \cite{boalch_RH}, p.~17, Eq.~(12).\\
 \textit{Solution 35}, 20 branches, $\boldsymbol{\theta}=(0,0,1/10,9/10)$:
 \begin{align*}
 &w=\frac{1}{2}-\frac{9 s^5-49 s^4-822 s^3+238 s^2-1699 s+1299}{2 (3 s-7) (s^2-2 s+17) u},\\
 &t=\frac{1}{2}-\frac{P(s)}{Q(s)u^3},\qquad\qquad u^2=(9 s^2 - 2 s + 9) (s^2 - 2 s + 17)\\
 &P(s)= 27 s^{10}-630 s^9+4055 s^8+30520 s^7-174970s^6+258492 s^5-724490 s^4\\
 &\qquad\quad+600760 s^3-1097825s^2+186570 s-131085,\\
 &Q(s)=2(s^2-2 s+17) (s^2-18 s+1),
 \end{align*}
 solution~45 of \cite{boalch52}, first obtained explicitly in \cite{boalch_HG}, p.~7.\\
 \textit{Solution 36}, 20 branches, $\boldsymbol{\theta}=(0,0,3/10,7/10)$:
 \begin{align*}
 &w=\frac{1}{2}-\frac{(s+3)(9s^4-100 s^3+118 s^2-228 s-55)}{(6 s^3-42s^2-30 s-62)u},\\
 &t=t_{35},\qquad u=u_{35},
 \end{align*}
  solution~44 of \cite{boalch52}, first obtained explicitly in \cite{boalch_HG}, p.~8.\\
 \textit{Solution 37}, 20 branches, $\boldsymbol{\theta}=(1/3,1/3,1/2,2/5)$:
 \begin{align*}
 &w=\frac{1}{2}+
\frac{(s + 3)P(s)}{
18 (s^2 + 1) (s^6 - 7 s^4 + 42 s^3 - 45 s^2 + 34 s + 7) u},\\
&t=\frac{1}{2} -
\frac{(s + 3) Q(s)}{
2 (s^2 + 1)^2 u^3},\qquad u^2 = 3 (s + 3) (8 s^2 - 13 s + 17),
 \end{align*}
 \begin{align*}
 &P(s) = 28 s^9-235 s^8+556 s^7-1334 s^6+2174 s^5-3854 s^4\\
 &\qquad\quad+4360 s^3-4738 s^2+2362 s-1047,\\
 &Q(s)=8 s^{10} + 100 s^7 - 135 s^6 + 834 s^5 - 1205 s^4 + 2280 s^3\\
 &\qquad\quad - 1365 s^2 + 890 s + 321,
 \end{align*}
 solution~43 in \cite{boalch52}, p. 24.\\
 \textit{Solution 38}, 20 branches, $\boldsymbol{\theta}=(1/3,1/3,1/2,4/5)$:
 \begin{align*}
  & w=\frac{1}{2}+\frac{(s + 3)
(8 s^6 - 28 s^5 + 85 s^4 - 196 s^3 + 214 s^2 - 196 s + 41) }{
6 (s^2 + 1) (3 s^2 - 4 s + 5) u},\\
& t=t_{37},\qquad u=u_{37},
 \end{align*}
 solution~42 in \cite{boalch52}, p. 24.\\
 \textit{Solution 39}, 24 branches, $\boldsymbol{\theta}=(1/3,1/3,1/3,1/2)$:
 \ben
 w=\frac{1}{2}-\frac{P(s)}{
 R(s)u},\quad
 t =\frac{1}{2}+
\frac{(s^2 + 4 s - 2)Q(s)}{
2 (s + 2) (3 s^2 - 2 s + 2)^2 u^3},\quad
u^2 =
(8 s^2 - 7 s + 2)
(s + 2),
 \ebn
 \begin{align*}
 &P(s)=16 s^{11}+72 s^{10}+50 s^9-242 s^8-3143 s^7+6562 s^6-8312 s^5\\
 &\qquad\quad+9760 s^4-9836 s^3+6216 s^2-2288 s+416,\\
 &Q(s)=8 s^{10} + 16 s^9 + 24 s^8 - 84 s^7 + 429 s^6 - 312 s^5 + 258 s^4 - 288 s^3 + 288 s^2 - 128 s+ 32,\\
 &R(s)= 2 (3 s^2 - 2 s + 2)(26 s^6 + 18 s^5 - 75 s^4 + 50 s^3 + 270 s^2 - 312 s + 104),
 \end{align*}
 solution~46 in \cite{boalch52}, p.~27.\\
 \textit{Solution 40}, 30 branches, $\boldsymbol{\theta}=(1/15,1/15,7/30,23/30)$:
 \begin{align*}
 &w=\frac{1}{2}+\frac{(s+1)(s^8+8 s^7+90 s^6+348 s^5+972 s^4+1296 s^3+4374 s^2+8748 s+19683)}{
 2(s+3)^2(s^4-4 s^3-6 s^2+81)u},\\
 &t=\frac{1}{2}+\frac{(s+1)^2(s+9)^2 P(s)}{2(s-3)^2(s+3)^5(s^2+9)u^3},\qquad u^2=(s+1)(s+9) (s^2+9)(s^2+4 s+9),\\
 &P(s)=s^{14}+10 s^{13}+63 s^{12}+180 s^{11}+621 s^{10}+3942 s^9+26595 s^8+99576 s^7+239355 s^6\\
 &\qquad\quad+319302 s^5+452709 s^4+1180980 s^3+3720087 s^2+5314410 s+4782969,
 \end{align*}
 solution~47 of \cite{boalch52}, first obtained explicitly in \cite{boalch_HG}, p.~9.\\
 \textit{Solution 41}, 30 branches, $\boldsymbol{\theta}=(2/15,2/15,1/30,29/30)$:
 \begin{align*}
 &w=\frac{1}{2}+\frac{(s+9)Q(s)}{2(s-3)(s+3)^4(s^2+9)u},\qquad t=t_{40},\qquad u=u_{40},\\
 &Q(s)=s^9+7 s^8+36s^7+36s^6+126s^5+1170s^4+8100s^3+18468s^2+24057s-6561,
 \end{align*}
 solution~48 of \cite{boalch52}, first obtained explicitly in \cite{boalch_HG}, p.~9.\\
 \textit{Solution 42}, 36 branches, $\boldsymbol{\theta}=(0,0,1/6,5/6)$:
 \begin{align*}
 &w=\frac{1}{2}-\frac{4s^9-24 s^8+84 s^7-240 s^6+96 s^5+1401 s^4-6396 s^3+11136s^2-8160 s-401}{
 2(2s^2-2s+5)(s^3-3s^2+3s-11)u},\\
 &t=\frac{1}{2}-\frac{(s-2)(s+4)P(s)}{4(s^2-7s+1)(s^2-4s+13)(2s^2-2s+5) u^3},\\
 &u^2=(s^2-4 s+13) (2 s^2-2s+5)(2 s^4+2 s^3-3 s^2-58 s+107),\\
 &P(s)=32s^{16}-640 s^{15}+6432 s^{14}-46016 s^{13}+266968s^{12}-1228152s^{11}
 +4546772 s^{10}\\
 &\qquad\quad -13723024 s^9+34628427s^8-74456536s^7+139564088 s^6-224784264s^5\\
 &\qquad\quad +300342142s^4-299494736 s^3+197723868s^2-68764168s+17918807,
 \end{align*}
 solution~49 of \cite{boalch52}, first obtained explicitly in \cite{boalch_HG}, p.~10.\\
 \textit{Solution 43}, 40 branches, $\boldsymbol{\theta}=(3/20,3/20,3/20,17/20)$:
 \begin{align*}
 &w=\frac{1}{2}+\frac{(s^2-18 s+1)(s^2-2s+17)\left(u_{35}\right)^2+8 (s+1)(3s^3-21 s^2-15 s-31)uv}{
 32(s^3+57s^2-69s+75)(s^2-1) v},\\
 &t=\frac{1}{2}+\frac{P_{35}(s) \,u}{1024(s-9)^2(s^2-1)^3(5s^2-2s+13)},\\
 &u^2=2(s-9)(s^2-1),\qquad v^2=-(s-1)(s-9)(5s^2-2 s+13),
 \end{align*}
 solution~50 of \cite{boalch52}, first obtained explicitly in \cite{boalch_HG}, p.~9. (The formula (6)
 for $v$ in \cite{boalch_HG}, p.~8 is incorrect and should be replaced with $v^2=-2(j+1)(5j^2-2j+13)$.
 This is undoubtedly a typing error, because the Maple file accompanying Arxiv version
 of \cite{boalch_HG} contains correct expressions, which yield a solution equivalent to the above).\\
 \textit{Solution 44}, 40 branches, $\boldsymbol{\theta}=(1/20,1/20,1/20,19/20)$:
 \begin{align*}
 &w=\frac{1}{2}+\frac{(s^2-18 s+1)\left(u_{35}\right)^2+4 (s-1)(3s-7)uv}{64(s+3)(s+1)^2 v},\\
 &t=t_{43},\qquad u=u_{43},\qquad v=v_{43},
 \end{align*}
 solution~51 of \cite{boalch52}, in explicit form first obtained in \cite{boalch_HG}, p.~8 (with
 the same misprints as solution~43).\\
 \textit{Solution 45}, 72 branches, $\boldsymbol{\theta}=(1/12,1/12,1/12,11/12)$:
 \begin{align*}
 &w=\frac{1}{2}+\frac{2(s^2-4s+13)(s^2-7s+1)\left(u_{42}\right)^2+9(s-1)(s^3+27s^2-57s+79)uv}{
 6(2s-7)^2 (s^2-1)(2s^2+s+17)(s^3-3s^2+3s-11) v},\\
 &t=\frac{1}{2}+\frac{(s-2)(s+4)P_{42}(s)}{54(2s-7)(s^2-1)(s^2-2s+6) u^3},\\
 &u^2=(2 s-7) (s^2-1)(2 s^2+s+17)(4 s^2-13 s+19),\\
 &v^2=-(s+1)(s^2-2s+6)(4s^2-13s+19),
 \end{align*}
 solution~52 of \cite{boalch52}, in explicit form first obtained  in \cite{boalch_HG}, p.~10--11.
 
 \end{document}